\documentclass[a4paper,11pt]{article}
\usepackage{geometry, graphicx, amsmath, amsthm, amssymb, mathrsfs, subfigure, color, tikz}
\def\x{{\boldsymbol x}}
\newcommand{\bb}[1]{\boldsymbol{\mathbf{#1}}}
\begin{document}
\title{Multiwavelet troubled-cell indicator for discontinuity detection of discontinuous Galerkin schemes}
\author{Mathea J. Vuik\footnote{Email:  M.J.Vuik@tudelft.nl.  Delft Institute of Applied Mathematics, Delft University of Technology, Mekelweg 4, 2628CD Delft, The Netherlands.}
\renewcommand{\thefootnote}{\fnsymbol{footnote}}
and
 \hspace{.05in}  Jennifer~K. Ryan\footnote{Corresponding Author. Email:  Jennifer.Ryan@uea.ac.uk. Telephone: +44 (0)1603 592586.   School of Mathematics, University of East Anglia, Norwich NR4 7TJ, United Kingdom.  Supported by the Air Force Office of Scientific Research, Air Force Material Command, USAF, under grant number FA8655-09-1-3055.}
}

\maketitle

\begin{abstract}
In this paper, we introduce a new {\it global} troubled-cell indicator for the discontinuous Galerkin (DG) method in one- and two-dimensions.  This is done by taking advantage of the global expression of the DG method and re-expanding it in terms of a multiwavelet basis, which is a sum of the global average and finer details on different levels.  Examining the higher level difference coefficients acts as a troubled-cell indicator, thus avoiding unnecessary increased computational cost of a new expansion. In two-dimensions the multiwavelet
decomposition uses combinations of scaling functions and multiwavelets in the $x-$ and $y-$directions for improved troubled-cell indication.  By using such a troubled-cell indicator, we are able to reduce the computational cost by avoiding limiting in smooth regions.  We present numerical examples in one- and two-dimensions and compare our troubled-cell indicator to the subcell resolution technique of Harten (1989) and the shock detector of Krivodonova, Xin, Remacle, Chevaugeon, and Flaherty (2004), which were previously investigated by Qiu and Shu (2005).
\end{abstract}



\section{Introduction and Motivation}
Nonlinear hyperbolic partial differential equations are often solved using the Runge-Kutta discontinuous Galerkin (DG) method \cite{Coc89S, Coc89LS, Coc90HS, Coc98S}. In practical applications, initial conditions may contain discontinuities, or the solution of a nonlinear equation may develop a shock at a certain time. To efficiently apply DG in case of discontinuous solutions, limiting techniques are used to reduce the spurious oscillations, that are developed in the discontinuous regions. Examples of these limiters are the minmod-based TVB limiter \cite{Coc89S}, TVD limiters \cite{Coc98}, WENO \cite{Shu88O, Shu89O}, and the moment limiter \cite{Kri07}.  Unfortunately, most of the limiters do not work well for higher-order approximations (limiting smooth extrema), or multidimensional cases. In order to limit the correct elements, a troubled-cell indicator can be used. This procedure detects discontinuous regions, where the use of a limiter is necessary. A limiter is then applied only in the identified troubled cells. In general, this leads to more accurate results in smooth regions, and reduces the computational cost significantly. 

There are a variety of troubled-cell indicators, some that are tied to the limiting procedure and others that are separate from this procedure. A few of the important methods of troubled-cell indication are minmod \cite{Coc89S}, Harten's subcell resolution \cite{Har89}, moment limiters \cite{Kri07}, monotonicity preserving limiters \cite{Sur97H}, and the shock detector of Krivodonova et al. (KXRCF) \cite{Kri04XRCF}.  These methods for indicating troubled cells were explored and compared by Qiu et al. in \cite{Qiu05S_2}.  They did this in order to improve the performance of a WENO-based limiter for DG.  They found that there is no universally better performing method for every problem.  However, they did find that the minmod-based limiter with a suitably chosen parameter, Harten's method and the KXRCF shock detector performed better than other methods. 

In this paper, we introduce a new troubled-cell indicator using ideas from a multiwavelet formulation.  We explain the relation between the multiwavelet expansion and the DG formulation, \cite{Alp02BGV, Arc11FS}. This multiwavelet expansion is decomposed into a sum of a global average and finer details on different levels. The absolute averages of the highest decomposition level act as a troubled-cell indicator, which suddenly increase in the neighborhood of a discontinuity \cite{Mal98}. In two-dimensions, the multiwavelet decomposition uses combinations of scaling functions and multiwavelets in the $x-$, and $y-$direction. This is the reason why we are able to detect the exact locations of discontinuities in the $x-$, or $y-$direction, or in one of the diagonal directions \cite{Mal98}.

This multiwavelet troubled-cell indicator takes a different tack than most troubled-cell indicators.  Instead of only considering local information, this technique uses {\it global} information to detect the troubled cells.  This technique performs well, even in the vicinity of a strong shock with weaker local shocks such as the double Mach reflection problem. It was recently pointed out by Zaide et al. \cite{Zai11R} that for systems, using local information one will find three different shock locations in each of the conserved variables.  However, by using global information, we obtain one location for the shock(s).  This allows us to implement a limiter in a smaller region, thereby reducing the time for computation. We demonstrate the robust performance of our indicator on a variety of test problems, using the moment limiter in the identified troubled cells \cite{Kri07}. The results using our new troubled-cell indicator will be compared with the method of Harten and the KXRCF shock detector and we show that it reliably performs better and more efficiently.

The outline of this paper is as follows:  in Section \ref{sec:background} we present the relevant background information in discontinuous Galerkin methods and multiwavelets.  In Section \ref{sec:shock_detection} we introduce our new global multiwavelet troubled-cell indicator.  The effectivity of this new method compared with existing methods is presented in Section \ref{sec:results} for standard numerical examples.  We conclude with a discussion of our method and future work in Section \ref{sec:conclude}.

\section{Background}
\label{sec:background}
In this section, relevant background information regarding discontinuous Galerkin methods and multiwavelets is presented, which can be found in \cite{Coc98,Nai11LL, Alp02BGV}. To begin, an explanation of the discontinuous Galerkin method in two dimensions is given.

\subsection{The discontinuous Galerkin method}
\label{sec:dg}
In order to describe the discontinuous Galerkin (DG) method, consider  the following partial differential equation on a rectangular domain $\Omega \in \mathbb{R}^2$:
\begin{subequations}
 \begin{align}
  u_t + \nabla \cdot {\bf f} (u) = 0, & \ \ \x=(x,y) \in \Omega, \ t \geq 0; \label{eqn:linadv}\\
  u(\x,0) = u^0(\x), & \ \ \x \in \Omega, \label{eqn:linadv_ic}
 \end{align}\label{eqn:linadvsystem}
\end{subequations}
where $u = u(\x,t)$, and ${\bf f}(u) = (f(u),g(u))^\top$ is the flux function.

To discretize in space, $\Omega$ is divided into $(N_x+1) \times (N_y+1)$ rectangular elements, given by,
\[
 I_{ij} = \{ (x,y): x \in (x_{i-\frac{1}{2}},x_{i+\frac{1}{2}}], y \in (y_{j-\frac{1}{2}},y_{j+\frac{1}{2}}] \}, \ i=0,\ldots,N_x, \ j = 0,\ldots,N_y.
\]
On each element, the chosen approximation space is defined as,
\[
 V_h(I_{ij}) = \left\{ v \in \mathbb{Q}^k(I_{ij}) \right\}, i=0,\ldots,N_x, j=0,\ldots,N_y.
\]
Here, $\mathbb{Q}^k$ is the space of polynomials, $\mathbb{Q}^k = \mbox{span} \{x^m y^n: 0 \leq m,n \leq k \}$. 

For simplicity, the basis of $\mathbb{Q}^k$ is constructed using a tensor product of the scaled Legendre polynomials, $\phi_{\ell_x}(x) \phi_{\ell_y}(y), \ \ell_x,\ell_y \in \{ 0,\ldots,k \}$. These functions are defined as,
\begin{equation}
 \phi_\ell(x) = \sqrt{\ell + \frac{1}{2}} P^{(\ell)}(x),\label{eqn:Legendre}
\end{equation}
where $P^{(\ell)}$ is the Legendre polynomial of degree $\ell \in \mathbb{N}$. We note that these functions are pairwise orthonormal:
\begin{equation}
 \langle \phi_\ell, \phi_m \rangle = \int_{-1}^1 \phi_\ell(x) \phi_m(x)dx = \delta_{\ell m}. \label{eqn:orthogonality}
\end{equation}
This choice of basis functions provides ease when pairing the discontinuous Galerkin method with a multiwavelet approximation.

The weak formulation of the DG method is obtained by multiplying equation (\ref{eqn:linadv}) by an arbitrary, smooth function $v \in V_h (I_{ij})$, and integrating over $I_{ij},$ \mbox{$i \in \{0,\ldots,N_x \},$} $ j \in \{0,\ldots,N_y \}$. Using the divergence theorem along with replacing $u$ by $u_h$ we obtain
\begin{align}
 \iint\limits_{I_{ij}} & u_{h,t}  \phi_{m_x}(\xi) \phi_{m_y}(\eta) dxdy =  \nonumber \\
 & \hspace{-0.1cm} \iint\limits_{I_{ij}} \left\{ f(u_h) \frac{d}{dx} \phi_{m_x}(\xi) \phi_{m_y}(\eta) + g(u_h) \phi_{m_x}(\xi) \frac{d}{dy} \phi_{m_y}(\eta) \right\} dxdy \nonumber \\
& -\int_{y_{j-\frac{1}{2}}}^{y_{j+\frac{1}{2}}} \left\{ \hat{f}_{i+\frac{1}{2},j}(\eta) \phi_{m_x}(1) \phi_{m_y}(\eta) - \hat{f}_{i-\frac{1}{2},j}(\eta) \phi_{m_x}(-1) \phi_{m_y}(\eta) \right\} dy  \nonumber\\
 & - \int_{x_{i-\frac{1}{2}}}^{x_{i+\frac{1}{2}}} \left\{ \hat{g}_{i,j+\frac{1}{2}}(\xi) \phi_{m_x}(\xi) \phi_{m_y}(1) - 
 \hat{g}_{i,j-\frac{1}{2}}(\xi)  \phi_{m_x}(\xi) \phi_{m_y}(-1) \right\} dx \label{eqn:weakformflux}
\end{align}
in local coordinates. The numerical flux functions at the boundaries are approximated using a monotone flux.  In our implementation the local Lax Friedrichs flux is used \cite{LeV02}.  This is defined as,
\[
 \hat{f}_{i+\frac{1}{2},j}(\eta) = \frac{f_{i+\frac{1}{2}, j}^+(\eta) + f_{i+\frac{1}{2},j}^-(\eta)}{2} - \frac{\lambda_{i+\frac{1}{2},j}(\eta)}{2} (u_{i+\frac{1}{2},j}^+(\eta) - u_{i+\frac{1}{2}}^-(\eta)), \\
\]
where,
\begin{equation}
 \lambda_{i+\frac{1}{2},j}(\eta) = \max(|f^\prime(u_h)|) \mbox{ over all $u_h$ between $u_h(x_{i+\frac{1}{2}},y)^-$ and $u_h(x_{i+\frac{1}{2}},y)^+$}, \nonumber
\end{equation}
where, $ y \in (y_{j-1/2}, y_{j+1/2}].$ If $f$ is convex, this reduces to,
\begin{equation}
 \lambda_{i+\frac{1}{2},j}(\eta) = \max (|f^\prime(u_h(x_{i+\frac{1}{2}},y)^-)|,|f^\prime(u_h(x_{i+\frac{1}{2}},y)^+)|).
\end{equation}\label{eqn:LLF}
The fluxes $\hat{g}_{i,j+1/2}(\xi)$ are computed similarly.

For time evolution, we choose the third order strong stability preserving Runge-Kutta scheme  \cite{Got98S}.  Note that we could use other strong stability preserving time-stepping methods, \cite{Got01ST,Ket09MG,Shu88}, this is only a choice that is made.

\subsection{Multiwavelets} 
\label{sec:multiwavelets}

In this section, a brief description of the theory of multiwavelets  \cite{Alp02BGV,Arc11FS} and the relation to the DG approximation is given. This is done for the one-dimensional case in Sections \ref{sec:scaling} (scaling function space), and \ref{sec:multi} (multiwavelets). An extension to two-dimensions is given in Section \ref{sec:2dmulti}.  Multiwavelets will be used as a troubled-cell indicator for the discontinuous Galerkin approximation.  Although these details are needed in order to put the DG approximation in the context of a multiwavelet basis, computationally it can be unnecessarily expensive.  In practice only portions of the multiwavelet expansion of the DG solution will be used in order to indicate troubled cells and therefore the  computational cost overall will not increase.

\subsubsection{Scaling function space} 
\label{sec:scaling}

\begin{figure}[ht!]
\centering
 \begin{center}
  \begin{tikzpicture}[scale=1]
\draw (0,10) -- (4,10);
\draw (0,10) node{$|$};
\draw (2,9.7) node{$I_0^0$};
\draw (4,10) node{$|$};
\draw (5,10) node{$V_0^{k+1}$};
\draw (7,10) node{Level 0};

\draw (0,9) -- (4,9);
\draw (0,9) node{$|$};
\draw (1,8.7) node{$I_0^1$};
\draw (2,9) node{$|$};
\draw (3,8.7) node{$I_1^1$};
\draw (4,9) node{$|$};
\draw (5,9) node{$V_1^{k+1}$};
\draw (7,9) node{Level 1};

\draw (0,8) -- (4,8);
\draw (0,8) node{$|$};
\draw (0.5,7.7) node{$I_0^2$};
\draw (1,8) node{$|$};
\draw (1.5,7.7) node{$I_1^2$};
\draw (2,8) node{$|$};
\draw (2.5,7.7) node{$I_2^2$};
\draw (3,8) node{$|$};
\draw (3.5,7.7) node{$I_3^2$};
\draw (4,8) node{$|$};
\draw (5,8) node{$V_2^{k+1}$};
\draw (7,8) node{Level 2};

\draw (2,7) node{$\vdots$};
\end{tikzpicture}
\caption{Visualization of the intervals in $V_0^{k+1},V_1^{k+1},\ldots$.}\label{fig:Nested_Vnp}
 \end{center}
\end{figure}
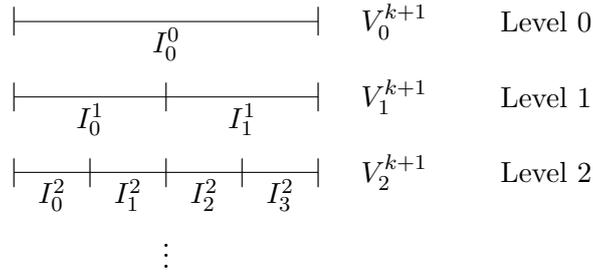

To begin defining multiwavelets in one-dimension, scaling functions defined on $[-1,1]$ are introduced.  We note that although we are using the multiwavelet decomposition in \cite{Alp02BGV,Arc11FS}, we must modify the definitions to accommodate this interval, which helps with use of the DG coefficients for the multiwavelet expansion.  Using this interval, the scaling function space is a space of piecewise polynomial functions, $V_n^{k+1}$, defined as,
\begin{equation}
 V_n^{k+1} = \{f: f \in \mathbb{P}^{k+1}(I_j^n), j=0,\ldots,2^n-1\}, 
 \label{eqn:V_n^p}
\end{equation}
where,
\begin{equation}
 I_j^n = (-1+2^{-n+1}j, -1+2^{-n+1}(j+1)], 
 \label{eqn:intervals}
\end{equation} 
and $\mathbb{P}^{k+1}(I_j^n)$ is the space of all polynomials of degree less than $k+1$ on interval $I_j^n$.
A visualization of the intervals in \mbox{$V_0^{k+1},V_1^{k+1},\ldots$} is demonstrated in Figure \ref{fig:Nested_Vnp}.  Notice that the space $ V_n^{k+1}$ has dimension $2^n(k+1)$ and the following nested property holds,
\[
 V_0^{k+1} \subset V_1^{k+1} \subset \cdots \subset V_n^{k+1} \subset \cdots.
\]

%In particular,
%\[
% V_0^p = \{ f \mbox{ is a polynomial of degree less than } p+1 \mbox{ on the interval } (-1,1] \}.
%\]

The scaled Legendre polynomials $\phi_0,\ldots,\phi_{k}$ used in the DG method (Section \ref{sec:dg}) are chosen to be the orthonormal basis for $V^{k+1}_0$. Next, the space $V_m^{k+1}, m \in \{0,\ldots,n \}$ is spanned by $2^m (k+1)$ functions which 
are obtained from $\phi_0,\ldots,\phi_{k}$ by dilation and translation,
\begin{equation}
 \phi_{\ell j}^m(x) = 2^{m/2} \phi_\ell(2^m(x+1) -2j -1), \ell = 0, \ldots, k, \ j = 0, \ldots, 2^m-1, \label{eqn:phi_{jl}^n}
\end{equation}
where the coefficient $j$ belongs to the various intervals $I_j^m$ \cite{Kei04}, and $x \in I_j^m$. The factor $2^{m/2}$ makes this an orthonormal basis for $V_m^{k+1}$.   The functions $\phi_\ell, \ell=0,\ldots,k,$ are called scaling functions.

Because the DG approximation and the scaling function approximation are composed of the same basis functions, there is a direct relation between the DG approximation and the scaling function approximation. In general, the DG method is applied on an interval $[a,b]$. If the number of elements in $[a,b]$ is chosen such that $N+1=2^n$, then (using definition (\ref{eqn:intervals})), the elements are given by,
\[
 I_j = \left( \left. a+ \frac{b-a}{2^n} j, a + \frac{b-a}{2^n}(j+1)\right.\right], \quad j = 0,\ldots,2^n-1.
 \]
 Noting that $\Delta x = (b-a)/2^n$, and $x_j = a + (j+1/2)\Delta x$,
the {\it global} DG approximation of the solution can be written as
\begin{equation}
 u_h(x,t) = 2^{-\frac{n}{2}}\sum_{j=0}^N\sum_{\ell=0}^k u_j^{(\ell)} \phi_{\ell j}^n(y), 
 \label{eqn:u_h1}
\end{equation}
where $y=-1+2 (x-a)/(b-a).$  However, exploiting the fact that  $u_h(x,t)$ is a piecewise polynomial of degree $k$ and transforming to a reference element, then the DG approximation projected onto the scaling function basis can be written as  
\begin{equation}
  u_h(x,t) = P_n^{k+1} u_h(x,t) = \sum_{j=0}^{2^n-1} \sum_{\ell=0}^k s_{\ell j}^n \phi_{\ell j}^n (y). \label{eqn:u_h2}
\end{equation}
From equations (\ref{eqn:u_h1}) and (\ref{eqn:u_h2}) it follows that for every $\ell = 0,\ldots,k, j=0,\ldots,2^n-1$,
\begin{equation}
 2^{-\frac{n}{2}} u_j^{(\ell)} = s_{\ell j}^n,
 \label{eqn:application}
\end{equation}
thus giving a relation between the coefficients of the DG approximation and the scaling function approximation.

\subsubsection{Multiwavelets} 
\label{sec:multi}

In addition to the scaling function space, $V_m^{k+1},$ a multiwavelet subspace is also needed in order to define the multiwavelet expansion.  Define the multiwavelet subspace $W_m^{k+1},$ to be the orthogonal complement of $V_m^{k+1}$ in $V_{m+1}^{k+1}$:
\begin{equation}
 V_m^{k+1} \oplus W_m^{k+1} = V_{m+1}^{k+1}, \ W_m^{k+1} \perp V_m^{k+1}, \ W_m^{k+1} \subset V_{m+1}^{k+1}, \ m = 0,\ldots, n-1. \label{eqn:orthogonal_complement}
\end{equation}
Note that $V_n^{k+1}$ can be split into $n+1$ orthogonal subspaces:
\[
 V_n^{k+1} = V_0^{k+1} \oplus W_0^{k+1} \oplus W_1^{k+1} \oplus \cdots \oplus W_{n-1}^{k+1}. \label{eqn:orthogonal_subspaces}
\]

By definition (\ref{eqn:orthogonal_complement}), the orthonormal basis for $W_0^{k+1}$ is given by $k+1$ piecewise polynomials, $\psi_0,\ldots,\psi_k$ (polynomials on $(-1,0]$ and $(0,1]$), which are the so-called multiwavelets. 
The term multiwavelet refers to the fact that the bases for $V_0^{k+1}$ and $W_0^{k+1}$ contain multiple elements.  The multiwavelet basis that belongs to the scaled Legendre polynomials, was developed by Alpert \cite{Alp93}. A more throrough explanation is given in \cite{Hov10MS}. Similar to the basis for $V_m^{k+1}$, the space $W_m^{k+1}$ is spanned by the functions,
\[
 \psi_{\ell j}^m(x) = 2^{m/2} \psi_\ell(2^m(x+1) -2j -1), \ell = 0, \ldots, k, \ j = 0, \ldots, 2^m-1, x \in I_j^m.
\]
Note that in general, a linear combination of $\psi_0,\ldots,\psi_k$, is continuous on $(-1,0]$ and $(0,1]$. On level $m$, multiwavelet $\psi_{\ell j}^m$ is continuous on 
\begin{align*}
  (-1 + 2^{-m+1}j, -1 + 2^{-m+1}(j+\frac{1}{2})] &= I_{2j}^{m+1}, \mbox{ and } \\
  (-1 + 2^{-m+1}(j+\frac{1}{2}), -1 + 2^{-m+1}(j+1)] &= I_{2j+1}^{m+1}.
\end{align*}

\bigskip
The multiwavelet expansion of a function $f \in L^2(-1,1)$ in level $m+1$ is given by
\begin{equation}
  Q_m^{k+1}f(x) = P_{m+1}^{k+1}f(x) - P_m^{k+1}f(x) = \sum_{j=0}^{2^m-1} \sum_{\ell=0}^{k}d_{\ell j}^m\psi_{\ell j}^m(x), \label{eqn:Q_n^p}
\end{equation}
which uses the multiwavelets $\psi_{\ell j}^m, \ell=0,\ldots,k, j=0,\ldots,2^m-1$. The coefficients are defined to be
\begin{equation}
 d_{\ell j}^m = \langle f,\psi_{\ell j}^m \rangle = \int_{-1+2^{-m+1}j}^{-1+2^{-m+1}(j+1)} f(x) \psi_{\ell j}^m(x)dx. \label{eqn:dljn}
\end{equation}
Using equation (\ref{eqn:Q_n^p}) recursively, a relation between expansions at different levels (starting from level $n$) is given as
\begin{equation}
  P_n^{k+1}f(x) = S^0(x) + \sum_{m=0}^{n-1} \sum_{j=0}^{2^m-1} D_j^m(x) = S^0(x) + \sum_{m=0}^{n-1}D^m(x).
\end{equation}\label{eqn:multiscale_transformation}

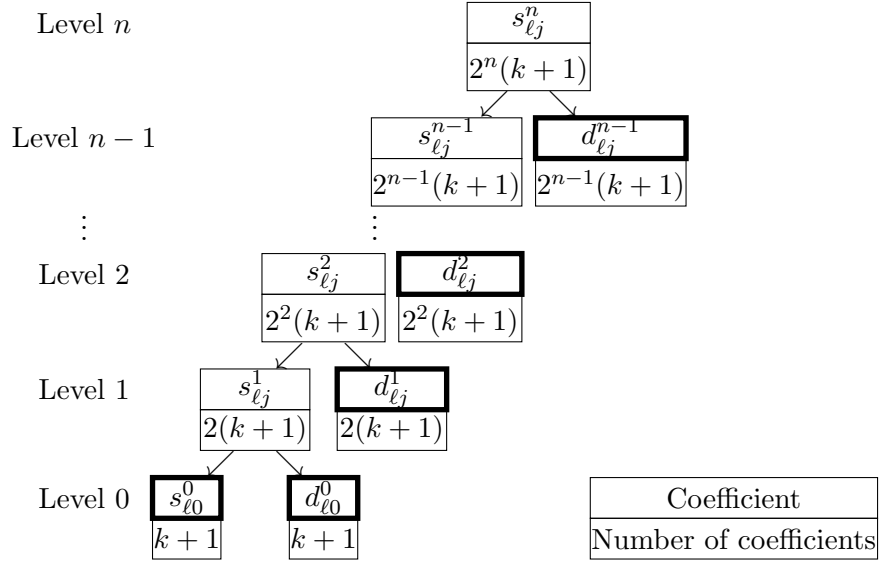
\begin{figure}[ht!]
\centering
\begin{tikzpicture}[scale=0.9]
% Level 0
\draw (0,0) node{$s_{\ell 0}^0$};
\draw (0,-0.6) node{$k+1$};
\draw (-0.5,-0.9) rectangle (0.5,-0.3);
\draw[line width = 2pt] (-0.5,-0.3) rectangle (0.5,0.3);

\draw (2,0) node{$d_{\ell 0}^0$};
\draw (2,-0.6) node{$k+1$};
\draw (1.5,-0.9) rectangle (2.5,-0.3);
\draw[line width = 2pt] (1.5,-0.3) rectangle (2.5,0.3);

\draw (-1.5,0) node{Level 0};

\draw (0.5,0.5) node{$\swarrow$};
\draw (1.5,0.5) node{$\searrow$};

% Level 1
\draw (1,1.6) node{$s_{\ell j}^1$};
\draw (0.2,1.3) rectangle (1.8,1.9);
\draw (1,1) node{$2(k+1)$};
\draw (0.2, 0.7) rectangle (1.8, 1.3);

\draw (3,1.6) node{$d_{\ell j}^1$};
\draw (3,1) node{$2(k+1)$};
\draw (2.2, 0.7) rectangle (3.8, 1.3);
\draw[line width = 2pt] (2.2,1.3) rectangle (3.8,1.9);

\draw (-1.5,1.6) node{Level 1};

\draw (1.5,2.1) node{$\swarrow$};
\draw (2.5,2.1) node{$\searrow$};

% Level 2
\draw (2,3.3) node{$s_{\ell j}^2$};
\draw (1.1,3) rectangle (2.9,3.6);
\draw (2,2.6) node{$2^2(k+1)$};
\draw (1.1, 2.3) rectangle (2.9, 3);

\draw (4,3.3) node{$d_{\ell j}^2$};
\draw (4,2.6) node{$2^2(k+1)$};
\draw (3.1, 2.3) rectangle (4.9, 3);
\draw[line width = 2pt] (3.1,3) rectangle (4.9,3.6);

\draw (-1.5,3.3) node{Level 2};

\draw (2.75,4.1) node{$\vdots$};

\draw (-1.5,4.1) node{$\vdots$};

% Level n - 1
\draw (3.8,5.3) node{$s_{\ell j}^{n-1}$};
\draw (2.7,5) rectangle (4.9,5.6);
\draw (3.8,4.6) node{$2^{n-1}(k+1)$};
\draw (2.7, 4.3) rectangle (4.9, 5);

\draw (6.2,5.3) node{$d_{\ell j}^{n-1}$};
\draw (6.2,4.6) node{$2^{n-1}(k+1)$};
\draw (5.1,4.3) rectangle (7.3,5);
\draw[line width = 2pt] (5.1,5) rectangle (7.3,5.6);

\draw (-1.5,5.3) node{Level $n-1$};

\draw (4.5,5.8) node{$\swarrow$};
\draw (5.5,5.8) node{$\searrow$};

% Level n
\draw (5,7) node{$s_{\ell j}^{n}$};
\draw (4.1,6.7) rectangle (5.9,7.3);
\draw (5,6.3) node{$2^n(k+1)$};
\draw (4.1, 6) rectangle (5.9, 6.7);

\draw (-1.5,7) node{Level $n$};

\draw (8, 0) node{Coefficient};
\draw (8,-0.6) node{Number of coefficients};
\draw (5.9, -0.3) rectangle (10.1, 0.3);
\draw (5.9, -0.9) rectangle (10.1, -0.3);
\end{tikzpicture}
\caption{Multiwavelet decomposition. Marked coefficients together carry the same information as $s_{\ell j}^n, \ell=0,\ldots, k, j=0,\ldots,2^n-1$.}\label{fig:multiscale_transformation}
\end{figure}

This representation of $P_n^{k+1}f(x)$ is called the multiscale decomposition, where
\[ 
 S^0(x) = \sum_{\ell=0}^{k}s_{\ell0}^0 \phi_\ell(x), \ \ D_j^m(x) = \sum_{\ell=0}^{k}d_{\ell j}^m\psi_{\ell j}^m(x).
\]
The coefficients $\{ s_{\ell0}^0\}_{\ell=0}^{k}$ represent the approximate solution on the coarsest level $n=0.$ 
The coefficients $\{ d_{\ell j}^m\}$ carry the multiscale information. These detail coefficients can be seen as carriers of individual fluctuations of the solution, which, if added to the lowest-resolution information, enrich it up to the level $n$
of resolution \cite{Iac11MMS}. The multiwavelet decomposition can be seen in Figure \ref{fig:multiscale_transformation}, which can be computed using quadrature mirror filter coefficients.

\begin{figure}[ht!]
 \begin{center}
  \includegraphics[clip=true, viewport=3cm 20cm 16cm 26.5cm, scale = 0.6]{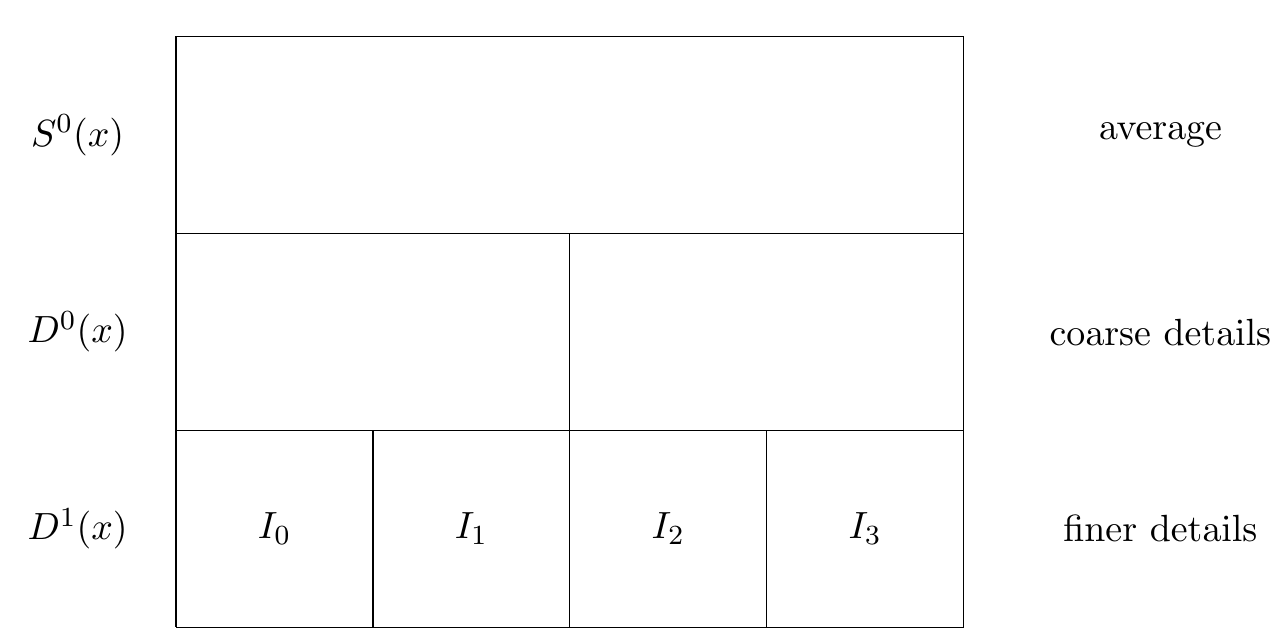}
\caption{Regions where multiwavelet contributions are continuous for $n=2$. Note that the highest level $D^1(x)$ is continuous in exactly the same regions as $u_h(x)$ and correspond to the DG elements.}\label{fig:continuousregions}
 \end{center}
\end{figure}

For $n=2$, the regions where the multiwavelet contributions are continuous are given in Figure \ref{fig:continuousregions}. It should be noticed that the DG approximation, $u_h(x)$, and the multiwavelet contribution of the highest level, $D^1(x)$, are both continuous in DG elements $I_0,\ldots,I_3$ and discontinuous on its boundaries.  In general, contribution $D^{n-1}(x)$ is continuous in exactly the same regions as $u_h(x)$, and level $n-1$ contains the most detailed information.

\subsubsection{Two-dimensional multiwavelet expansions} \label{sec:2dmulti}
% % % 1_2016.pdf, 2_2016.pdf
The two-dimensional multiscale decomposition of the discontinuous Galerkin approximation $u_h$ is more complex than the one-dimensional case. The one-step decomposition is given by
\begin{align*}
  u_h(x,y) = \hspace{-0.2cm} \sum_{i=0}^{2^{n_x-1}-1}  \sum_{j=0}^{2^{n_y-1}-1} \sum_{\ell_x,\ell_y=0}^k & \left\{ s_{\bb{\ell j}}^{\bb{n}-1} \phi_{\ell_x,i}^{n_x-1}(x) \phi_{\ell_y,j}^{n_y-1}(y) +  d_{\bb{\ell j}}^{\alpha, \bb{n}-1} \phi_{\ell_x,i}^{n_x-1}(x) \psi_{\ell_y,j}^{n_y-1}(y) \right.  \\
&  \left.  +  d_{\bb{\ell j}}^{\beta, \bb{n}-1} \psi_{\ell_x,i}^{n_x-1}(x) \phi_{\ell_y,j}^{n_y-1}(y) +  d_{\bb{\ell j}}^{\gamma, \bb{n}-1} \psi_{\ell_x,i}^{n_x-1}(x) \psi_{\ell_y,j}^{n_y-1}(y) \right\},
\end{align*}
and the full decomposition can be written as
 \[
 u_h(x,y) = S^0(x,y) + \sum_{m_y = 0}^{n_y - 1} D^{\alpha, m_y}(x,y) + \sum_{m_x=0}^{n_x-1}  D^{\beta, m_x}(x,y) + \sum_{m_x=0}^{n_x-1} \sum_{m_y = 0}^{n_y - 1} D^{\gamma, {\bf m}}(x,y), 
\]
where, 
\begin{align*}
S^0(x,y) &= \sum_{\ell_x,\ell_y=0}^k s_{{\bf \boldsymbol{\ell} 0}}^{\mathbf{0}} \phi_{\ell_x}(x) \phi_{\ell_y}(y), \\
 D^{\alpha, m_y}(x,y) & = \sum_{j=0}^{2^{m_y} - 1} \sum_{\ell_x,\ell_y=0}^k d_{\boldsymbol{\ell},(0,j)}^{\alpha, (0,m_y)} \phi_{\ell_x}(x) \psi_{\ell_y, j}^{m_y}(y), \\
 D^{\beta, m_x}(x,y) & = \sum_{i=0}^{2^{m_x} - 1}  \sum_{\ell_x,\ell_y=0}^k d_{\boldsymbol{\ell},( i,0)}^{\beta, (m_x,0)} \psi_{\ell_x, i}^{m_x}(x) \phi_{\ell_y}(y), \\
 D^{\gamma, {\bf m}}(x,y) & = \sum_{i=0}^{2^{m_x} - 1} \sum_{j=0}^{2^{m_y} - 1} \sum_{\ell_x,\ell_y=0}^k d_{{\bf \boldsymbol{\ell} j}}^{\gamma, {\bf m}} \psi_{\ell_x, i}^{m_x}(x) \psi_{\ell_y, j}^{m_y}(y),
\end{align*}
and $\boldsymbol{\ell} = (\ell_x, \ell_y)^\top, \ {\bf j} = (i,j)^\top, {\bf m} = (m_x, m_y)^\top$. 

Thus, modes $\alpha$, $\beta$ and $\gamma$ detect troubled cells which are oriented in the $x-$, $y-$, and $xy-$directions, respectively \cite{Mal98}. 

Using the quadrature mirror filter coefficients \cite{Alp02BGV}, the lower-level multiwavelet coefficients can be computed using the relations
\begin{subequations}
 \begin{align}
 s_{{\bf \boldsymbol{\ell} j}}^{{\bf m} - {\bf 1}} &=   \sum_{\tilde{i},\tilde{j}=0}^1 \sum_{r_x,r_y=0}^k h_{\ell_x, r_x}^{(\tilde{i})}h_{\ell_y, r_y}^{(\tilde{j})} s_{{\bf r},2{\bf j}+\tilde{{\bf j}}}^{\bf m}; \\
d_{{\bf \boldsymbol{\ell} j}}^{\alpha, {\bf m} - {\bf 1}} & = \sum_{\tilde{i},\tilde{j} = 0}^1 \sum_{r_x,r_y=0}^k h_{\ell_x,r_x}^{(\tilde{i})}g_{\ell_y,r_y}^{(\tilde{j})} s_{{\bf r},2{\bf j}+\tilde{{\bf j}}}^{\bf m}; \\
d_{{\bf \boldsymbol{\ell} j}}^{\beta, {\bf m} - {\bf 1}} & = \sum_{\tilde{i},\tilde{j} = 0}^1 \sum_{r_x,r_y=0}^k g_{\ell_x,r_x}^{(\tilde{i})}h_{\ell_y,r_y}^{(\tilde{j})} s_{{\bf r},2{\bf j}+\tilde{{\bf j}}}^{\bf m}; \\
 d_{{\bf \boldsymbol{\ell} j}}^{\gamma, {\bf m} - {\bf 1}} & = \sum_{\tilde{i},\tilde{j} = 0}^1 \sum_{r_x,r_y=0}^k g_{\ell_x,r_x}^{(\tilde{i})}g_{\ell_y,r_y}^{(\tilde{j})} s_{{\bf r},2{\bf j}+\tilde{{\bf j}}}^{\bf m}.
\end{align}
\end{subequations}
Analogous to the one-dimensional case,  it holds that
\begin{equation}
 2^{-\frac{n_x + n_y}{2}} u_{ij}^{(\ell_x,\ell_y)} = s_{{\bf \boldsymbol{\ell} j}}^{\bf n}.
\end{equation}

\subsection{Limiting and troubled-cell indication for DG}

\subsubsection{Troubled-cell indicators}

In this section, we look at the Harten troubled-cell indicator (developed by Qiu et al. \cite{Qiu05S_2}), and a shock detection technique by Krivodonova et al. \cite{Kri04XRCF}.  Note that these were the methods deemed to be the most reasonable by Qiu et al. \cite{Qiu05S_2}.

The Harten indicator is based on Harten's subcell resolution idea \cite{Har89}. In one-dimension, we define,
\[
F_i(z) = \frac{1}{\Delta x} \left\{ \int_{x_{i-\frac{1}{2}}}^z u_h|_{I_{i-1}} dx +  \int_z^{x_{i+\frac{1}{2}}} u_h|_{I_{i+1}} dx \right\} - \bar{u}_{ij},
\]
where $u_h|_{I_{i-1}}$ and $u_h|_{I_{i+1}}$ are extensions of the DG approximation in $I_{i-1}$ and $I_{i+1}$ into element $I_i$.   Element $I_i$ is marked as a troubled cell if 
\[
 F_i(x_{i-\frac{1}{2}}) F_i(x_{i+\frac{1}{2}}) \leq 0 \mbox{ and } |u_i^{(k)}| > \alpha |u_{i-1}^{(k)}|, \ |u_i^{(k)}| > \alpha |u_{i+1}^{(k)}|.
\]
Here, $\alpha$ is a parameter, which is chosen to be equal to 1.5 in \cite{Qiu05S_2}. Note that the choice of $\alpha$ also depends on the choice of limiter that is applied in the troubled cells.

The KXRCF indicator \cite{Kri04XRCF} uses inflow boundaries to detect troubled cells. Here, element $I_i$ is detected if 
\[
 \frac{\left| \int_{\partial I_i^-} (u_h|_{I_i} - u_h|_{I_{n_i}}) ds \right|}{h^{\frac{k+1}{2}} |\partial I_i^-| ||u_h|_{I_i}||} > 1. 
\]
Here, $\partial I_i^-$ is the inflow boundary, $u_h|_{I_{n_i}}$ is the DG approximation in the neighbor of $I_i$ on the side of $\partial I_i^-$, $h$ is the radius of the circumscribed circle in $I_i$, and the norm is based on the average in one-dimension and the maximum norm in quadrature points in two-dimensions. 

In Section \ref{sec:results}, the methods of Harten and KXRCF are compared against our multiwavelet troubled-cell indicator. Additionally, these methods will be used in combination with the moment limiter of Krivodonova \cite{Kri07}, which is outlined in the next subsection.

\subsubsection{Limiting the DG approximation}

In general, the solution of a nonlinear PDE develops shocks and discontinuities in time. Higher order methods tend to introduce spurious oscillations in discontinuous regions, thereby loosing accuracy. One way to get rid of these oscillations is by applying a limiter, which reduces the approximation to a low order in discontinuous regions, and maintains the high order if the approximation is smooth enough. Examples of these limiters are the minmod-based TVB limiter \cite{Coc89S}, TVD limiters \cite{Coc98}, WENO \cite{Shu88O, Shu89O}, and the moment limiter \cite{Kri07}.   Each of these limiters has its own mechanism to control which regions should be limited. This often results in limiting smooth local extrema, which makes the approximation too diffusive. Our approach however, is to use a new global multiwavelet troubled-cell indicator as a switch to control where the limiter is applied.  

We have chosen to apply our indicator in combination with the moment limiter \cite{Kri07}, which limits DG coefficients, starting at the highest level $k$. For each element $I_j, j=0,\ldots,N$, we compute
\[
 \widetilde{u}_j^{(k)} = \mbox{minmod} \left( u_j^{(k)},\beta_k\left(u_{j+1}^{(k-1)} - u_j^{(k-1)}\right),\beta_k\left(u_j^{(k-1)}-u_{j-1}^{(k-1)}\right) \right),
\]
where the minmod function is given by,
\[
 \mathrm{minmod}(a_1,\ldots,a_q) = \left\{ \begin{array}{ll} \mbox{sign}(a_1) \cdot \min_{1 \leq r \leq q} |a_r|, & \mbox{sign }(a_1) = \cdots = \mbox{sign} (a_q), \\ 0, &  \mbox{otherwise.} \end{array} \right.
\]
and $\beta_k = (\sqrt{k-1/2})/(\sqrt{k+1/2}).$   If $\widetilde{u}_j^{(k)} = u_j^{(k)},$ the limiting procedure is cut off for this element $I_j$. Otherwise, $u_j^{(k-1)}$ is limited, continuing until $u_j^{(1)}$ is limited ($u_j^{(0)}$ remains the same, such that the average $\bar{u}_j$ is preserved), or stopping the first time $\widetilde{u}_j^{(\ell)} = u_j^{(\ell)}$ for some $\ell = k-1,\ldots,1.$  For a system of equations, the moment limiter is applied to the characteristic variables $\mathbf{w}_j^{(\ell)} = R^{-1}\mathbf{u}_j^{(\ell)}$. If we obtain negative values for density, pressure or energy (due to the characteristic approach), then we should set all higher order coefficients equal to zero, and limit the linear term. If negative values are still found, then the linear coefficient is also set equal to zero.  In two-dimensions, the moment limiter uses the neighboring elements both in the $x$-, and in the $y$-direction, \cite{Kri07}.

\section{Multiwavelet troubled-cell indicator}
\label{sec:shock_detection}
In this section a troubled-cell indicator that exploits information from the multiwavelet expansion is introduced.  This uses the relation between the {\it global} discontinuous Galerkin approximation and the multiwavelet expansion.  It is important to note that most troubled-cell indicators only use local information.  However, by using the relation between the DG approximation and multiwavelet expansion, the detection of troubled cells is exact, even in the presence of local extrema.

In the neighborhood of a discontinuity in the DG approximation, the multiwavelet contribution of the higher levels will suddenly become large with respect to this contribution in continuous regions. For various types of discontinuities, the same behavior is found: the multiwavelet contribution of, for example, level $n-1$, is large in the discontinuous region, compared to the continuous regions. In this paper, the authors propose to use the contribution $D^{n-1}(x)$, which belongs to level $n-1$, for troubled-cell indication, \cite{Mal98}.

The multiwavelet contribution in element $I_j^{n-1}$ is given by
\[
  D_j^{n-1}(x) = \sum_{\ell=0}^k d_{\ell j}^{n-1} \psi_{\ell j}^{n-1}(x), j \in \{ 0,\ldots,2^{n-1}-1 \}, x \in I_j^{n-1}. 
\]
Because $D_j^{n-1}$ is continuous on $I_{2j}^n$ and $I_{2j+1}^n$ (section \ref{sec:multiwavelets}), we propose to indicate troubled cells using the absolute averages of $D_j^{n-1}$ on elements $I_{2j}^n$ and $I_{2j+1}^n$. This corresponds to computing the averages on each element $I_i$ of the discretization, given by
\begin{equation}
 \bar{D}_i^{n-1} = \frac{1}{\Delta x}\int_{I_i}  \left| D^{n-1}(x)\right| dx, i=0,\ldots,2^n-1, \label{eqn:L1}
\end{equation}
which is the weighted $L^1$-norm on element $I_i$, generally used for discontinuity detection \cite{Har95}. The element where the average (\ref{eqn:L1}) is maximal, is assumed to be the element where the strongest shock occurs. 

Due to the computational cost of integral evaluation, the authors choose to implement a three point trapezoidal rule. This is done in place of exact integral evaluation, because finding the roots of the absolute multiwavelet decomposition is not easy. This discrete average of $D^{n-1}$ is easy and fast to compute, and gives a good approximation of the continuous average. 

Combining the shock detector with a limiting strategy, better numerical results are expected.  This is because the DG coefficients are limited only in the neighborhood of the shock.  There is no limiting occurring in continuous regions. As mentioned earlier, the choice made for the limiting step is the moment limiter.  This is applied in regions where the average (\ref{eqn:L1}) is large enough, that is,
\[
 \bar{D}_i^{n-1} > C \cdot \max \{ \bar{D}_i^{n-1}, i=0,\ldots,2^n-1 \}, \ C \in [0,1]. 
\]
If $C=1$, then no element will be detected. In this way, the value of $C$ is a useful tool to prescribe the strictness of the limiter. The lower the value of C, the more cells are limited.  

The complexity of extending these ideas to two dimensions does not increase considerably (Section \ref{sec:results}). The main difference is that now there are separate detail coefficients for the $x,\ y,$ and $xy-$directions.  This means that three different approximations must be computed:  $D^{\alpha,{\bf n} - {\bf 1}},\ D^{\beta, {\bf n} - {\bf 1}},\ D^{\gamma, {\bf n} - {\bf 1}}.$  

To compute the averages of $D^{\alpha, {\bf n} - {\bf 1}},$ one bases the computation on its construction through the functions $\phi_{\ell_x,i}^{n_x-1}$ (polynomial on element $I_i^{n_x-1}$ ($x$-direction)) and $\psi_{\ell_y,j}^{n_y-1}$ (piecewise polynomial on elements $I_{2j}^{n_y}$ and $I_{2j+1}^{n_y}$) for $\ell_x,\ell_y \in \{0,\ldots,k\}$, and where $i \in \{0,\ldots,2^{n_x - 1}-1\}, j \in \{0,\ldots,2^{n_y - 1}-1\}$. Similar to the one-dimensional approach, we compute,
\[
 \bar{D}_{ij}^{\alpha, {\bf n} - {\bf 1}} \equiv \bar{D}^{\alpha, {\bf n} - {\bf 1}}(I_i^{n_x - 1} \times I_j^{n_y}), i = 0,\ldots, 2^{n_x - 1} - 1, j = 0, \ldots, 2^{n_y} - 1,
\]
resulting in $2^{n_x - 1} \cdot 2^{n_y}$ averages.

For the $\beta$ mode, 
\[
 \bar{D}_{ij}^{\beta, {\bf n} - {\bf 1}} \equiv \bar{D}^{\beta, {\bf n} - {\bf 1}}(I_i^{n_x} \times I_j^{n_y - 1}), i = 0,\ldots, 2^{n_x} - 1, j = 0, \ldots, 2^{n_y - 1} - 1,
\]
are computed because multiwavelet $\psi_{\ell_x, i}^{n_x - 1}(x)$ is used in the $x-$direction and scaling function $\phi_{\ell_y, j}^{n_y - 1}(y)$ in the $y-$direction ($2^{n_x} \cdot 2^{n_y - 1}$ averages). 

In mode $\gamma$, multiwavelets both in the $x-$ and in the $y-$direction are used, such that  
\[
 \bar{D}_{ij}^{\gamma, {\bf n} - {\bf 1}} \equiv \bar{D}^{\gamma, {\bf n} - {\bf 1}} (I_i^{n_x} \times I_j^{n_y}), i=0,\ldots,2^{n_x}-1, j=0,\ldots,2^{n_y}-1
\]
is found ($2^{n_x} \cdot 2^{n_y}$ averages).

Analogous to the one-dimensional case, the element $I_i^{n_x-1} \times I_j^{n_y}$ is indicated to be a troubled cell in the $\alpha$ mode if
\[
 \bar{D}_{ij}^{\alpha, {\bf n} - {\bf 1}}  >  C^\alpha \cdot \max \{ \bar{D}_{ij}^{\alpha, {\bf n} - {\bf 1}}, i = 0,\ldots, 2^{n_x - 1} - 1, j = 0, \ldots, 2^{n_y} - 1 \}, \ C^\alpha \in [0,1].
\]
Shock detection in the $\beta$, and $\gamma$ mode is done in the same manner, using the constants $C^\beta$ and $C^\gamma$ to determine the strictness of the troubled-cell indicator. Note that this gives us three parameters to choose.  
Similar to the one-dimensional case, more elements are detected if the values of $C^\alpha, C^\beta,$ and $C^\gamma$ are smaller.  Using this approach,  the $\alpha$ mode detects discontinuities in the $y-$direction (because multiwavelets are used in the $y-$direction), and the $\beta$ mode detects discontinuities in the $x-$direction (multiwavelets in $x$).  The $\gamma$ mode is used for diagonal shock detection, \cite{Mal98}.

\section{Numerical results}
\label{sec:results}

In this section, we look at various examples in order to investigate the effectiveness of the multiwavelet troubled-cell indicator applied to the discontinuous Galerkin approximation.  We compare the results with the subcell resolution method of Harten \cite{Har89} and the shock detection method of Krivodonova et al. \cite{Kri04XRCF}.  These results demonstrate that the multiwavelet troubled-cell indicator performs well using only a moderate computational cost.  

\subsection{One-dimensional Euler equations}

We begin by investigating the performance of the multiwavelet troubled-cell indicator for a nonlinear system of equations (Euler equations) and comparing this to existing troubled-cell indicators. These equations describe density, $\rho$, velocity, $u$, pressure, $p$, and energy, $E$, to form a system of conservation laws. Introducing $\mathbf{u} = (\rho, \rho u, E)^\top = (u^{(1)}, u^{(2)}, u^{(3)})^\top$, these equations are given by
\begin{subequations}
\begin{equation}
 \mathbf{u}_t + \mathbf{f(u)}_x = 0,
\end{equation}
where
\begin{equation}
\mathbf{f(u)} = (\rho u, \rho u^2+p, u(E+p))^\top.
\end{equation}
For simplicity, we use the equation of state for an ideal polytropic gas:
\begin{equation}
 E = \frac{p}{\gamma -1} + \frac{1}{2} \rho u^2.
\end{equation} \label{eqn:Euler}
\end{subequations}

Below the results and comparisons are given using four different sets of initial conditions: the shock tubes of Sod and Lax, the blast waves, and the Shu-Osher problem.  We compare the cells that are detected by our multiwavelet indicator with the KXRCF and Harten's troubled-cell indicator.  For the KXRCF and Harten's indicator, a combination of density and energy, or density and entropy is used in the literature \cite{Kri04XRCF, Qiu05S_2}. We have used density and entropy as indicator values.  In the multiwavelet approach, however, only density is used in the indicator. In Section \ref{sec:sod}, it can be seen that the multiwavelet approach using density detects exactly the same elements as the combination of density and entropy does.

\subsubsection{Sod's shock tube}\label{sec:sod}
Sod's shock tube problem models the situation where a diaphragm halfway inside a tube separates two regions which have different densities and pressures  \cite{Sod78}.  These two regions have constant states, with both fluids are initially at rest. The following initial condition is used:
\begin{subequations}
\begin{equation}
 \rho(x,0)  = \left\{ \begin{array}{ll} 1, & \mbox{for } x <0, \\ 0.125, & \mbox{for } x \geq 0 \end{array} \right. \quad
 p(x,0)  = \left\{ \begin{array}{ll} 1, & \mbox{for } x <0, \\ 0.1, & \mbox{for } x \geq 0 \end{array} \right.
\end{equation}
and
\begin{equation}
 u(x,0) \equiv 0.
\end{equation}
\label{eqn:Sod}
\end{subequations}
At time $t > 0$ the diaphragm is broken.  The physical domain is assumed to be essentially infinite. The computational domain, however,  is set equal to $[-5,5]$. Because we do not compute long enough for the waves to reach the boundaries, initial constant states are used as boundary conditions.  From the literature \cite{Smo83}, the detector should be able to identify the   shock, contact discontinuity and rarefaction wave.

\begin{figure}[ht!]
\centering
\subfigure[$k=1$]{\includegraphics[scale = 0.28]{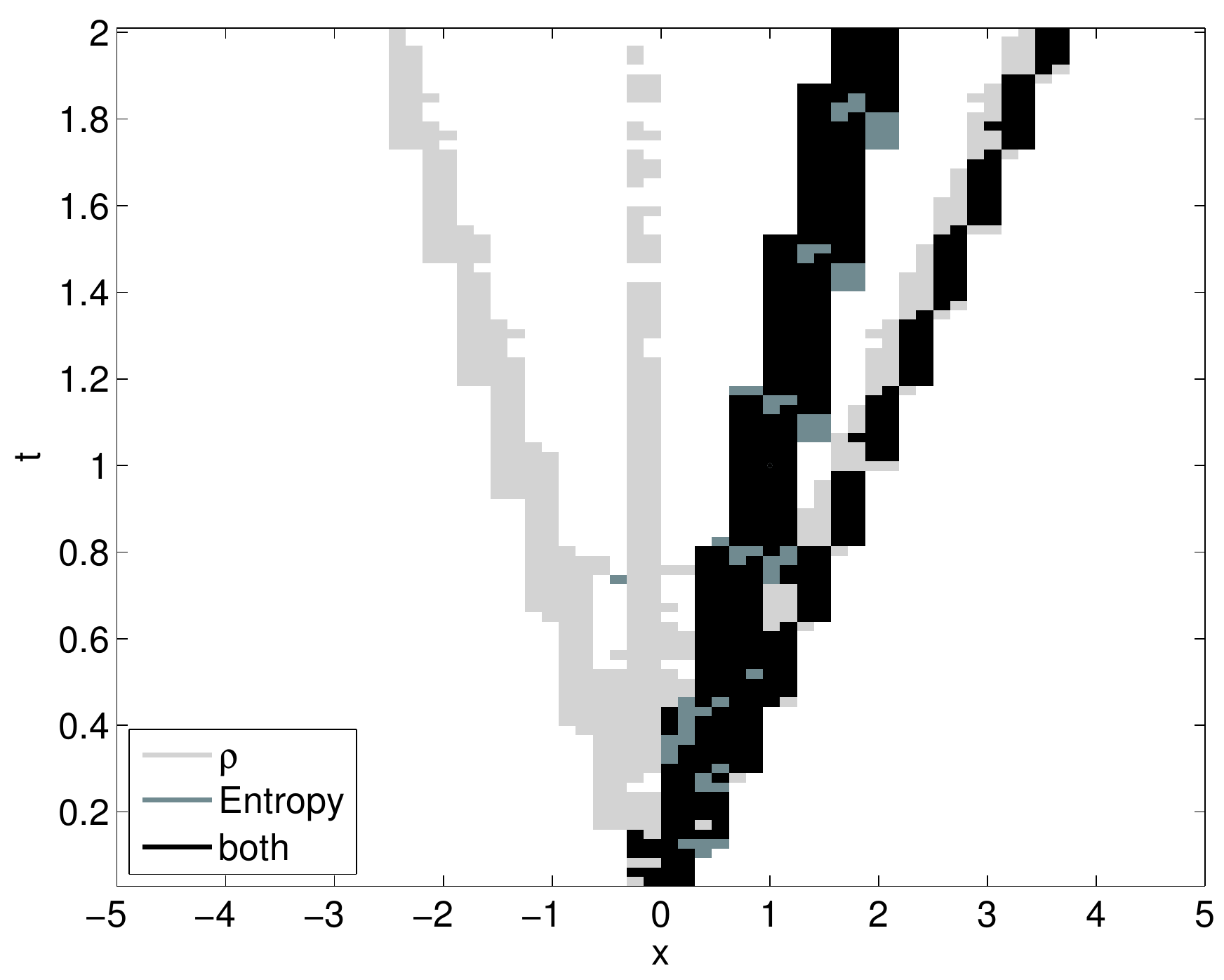}}
\subfigure[$k=2$]{\includegraphics[scale = 0.28]{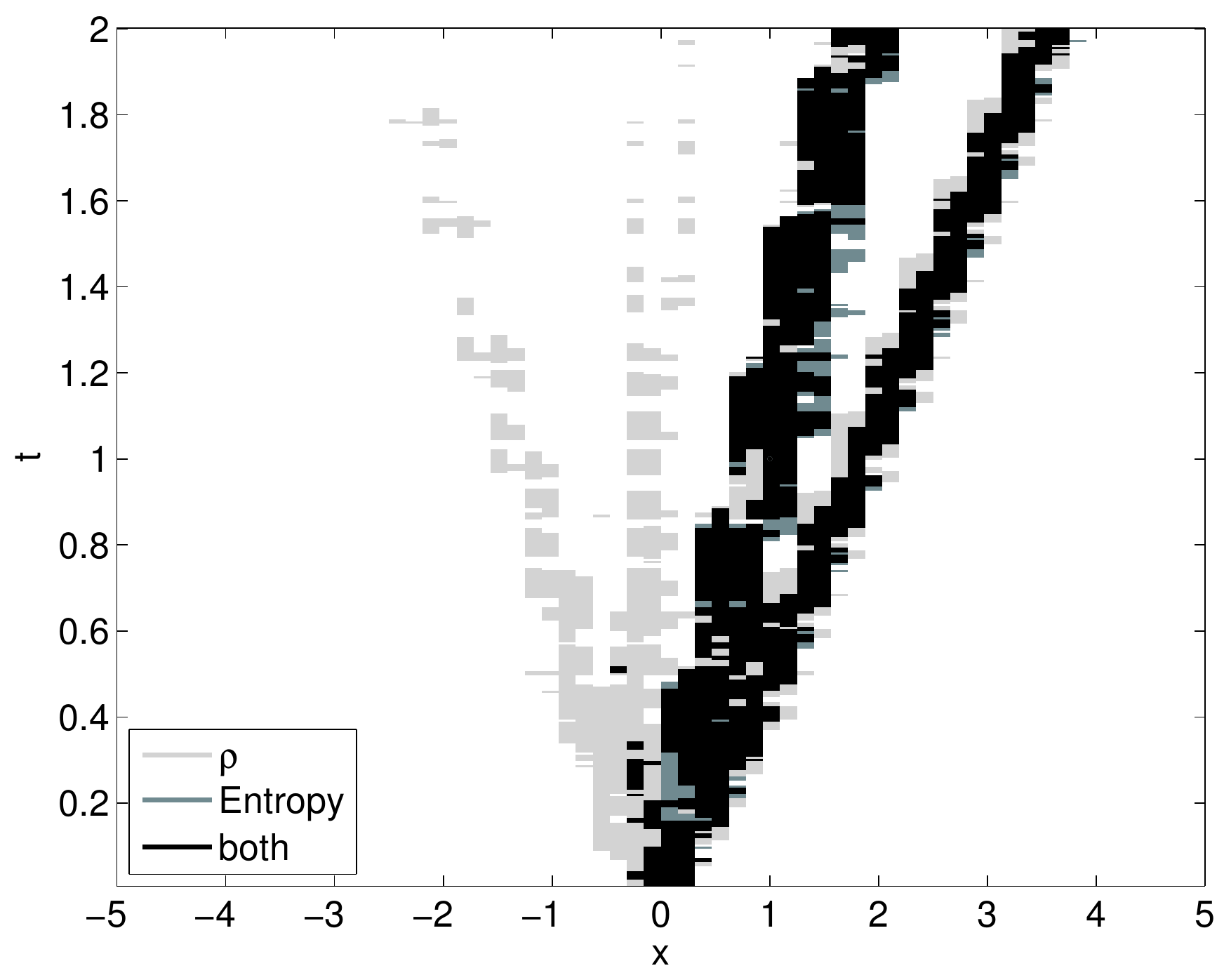}}
\caption{Time history plot of detected troubled cells using the multiwavelet troubled-cell indicator, $C=0.1$ (density and entropy), Sod, 64 elements, $k=1$. Compare to Figure \ref{fig:SodC}.}\label{fig:Soddensityentropy}
\end{figure}

\begin{figure}[ht!]
\centering
 \subfigure[$C=0.9$]{\includegraphics[scale = 0.22]{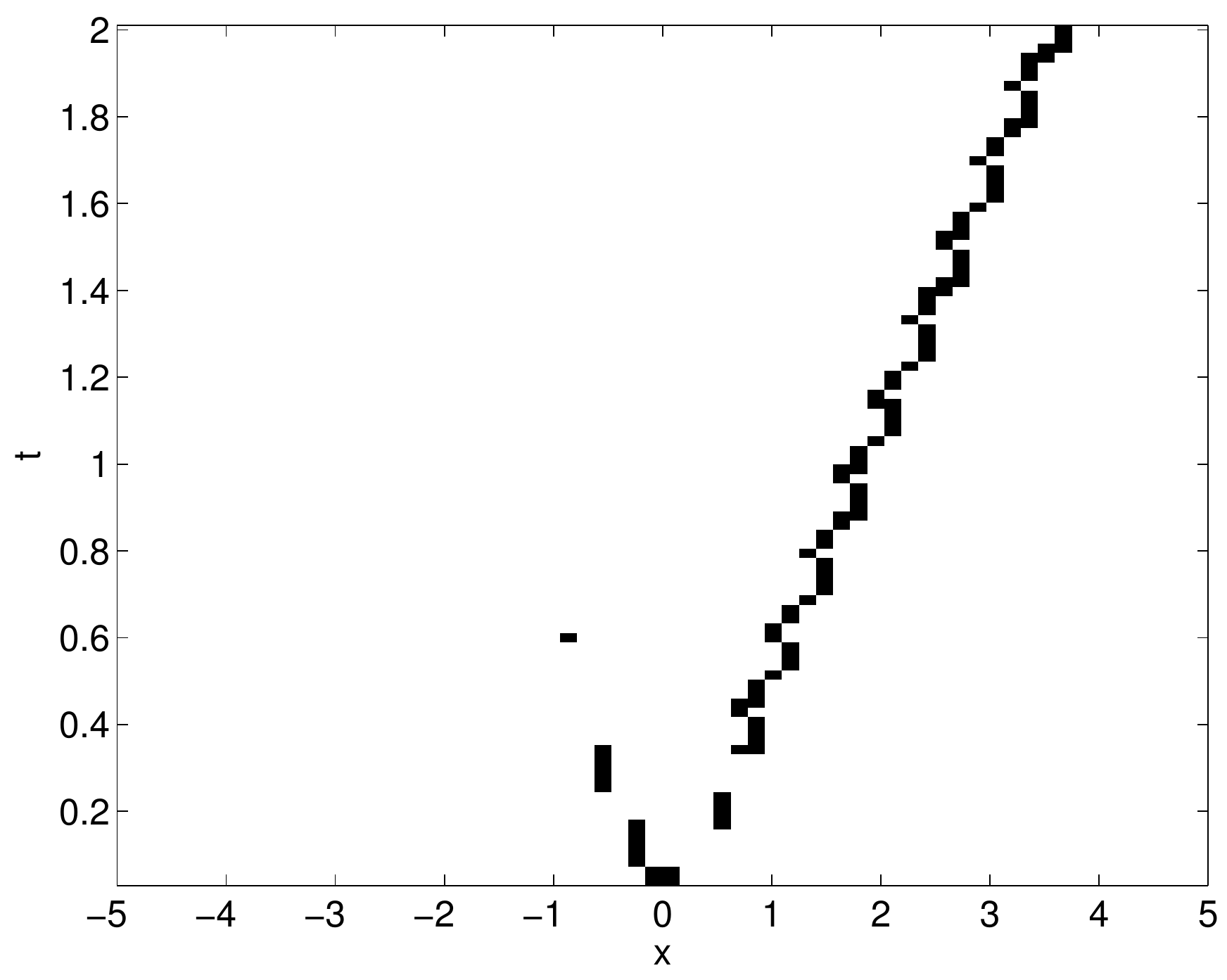}}
 \subfigure[$C=0.5$]{\includegraphics[scale = 0.22]{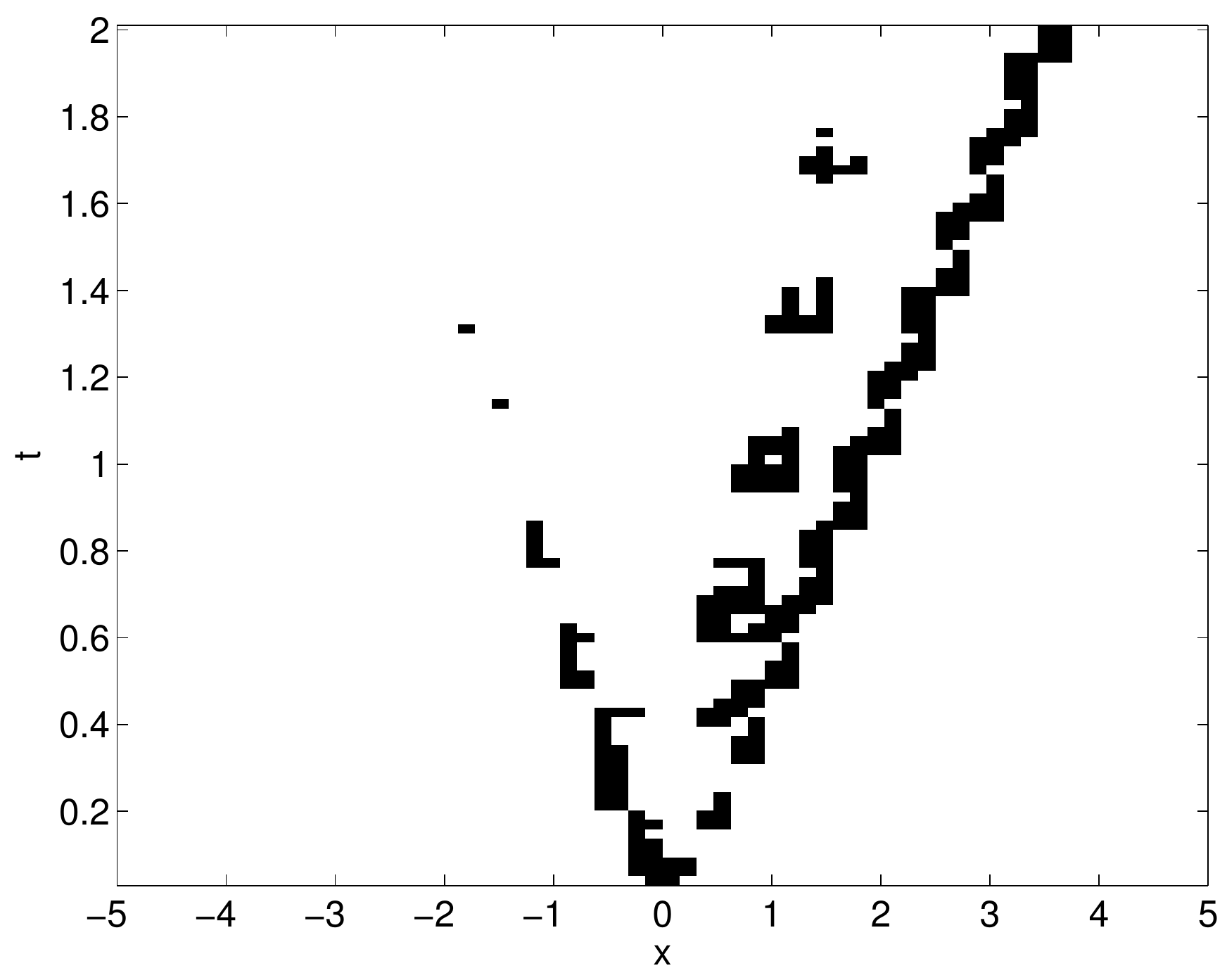}}
 \subfigure[$C=0.1$]{\includegraphics[scale = 0.22]{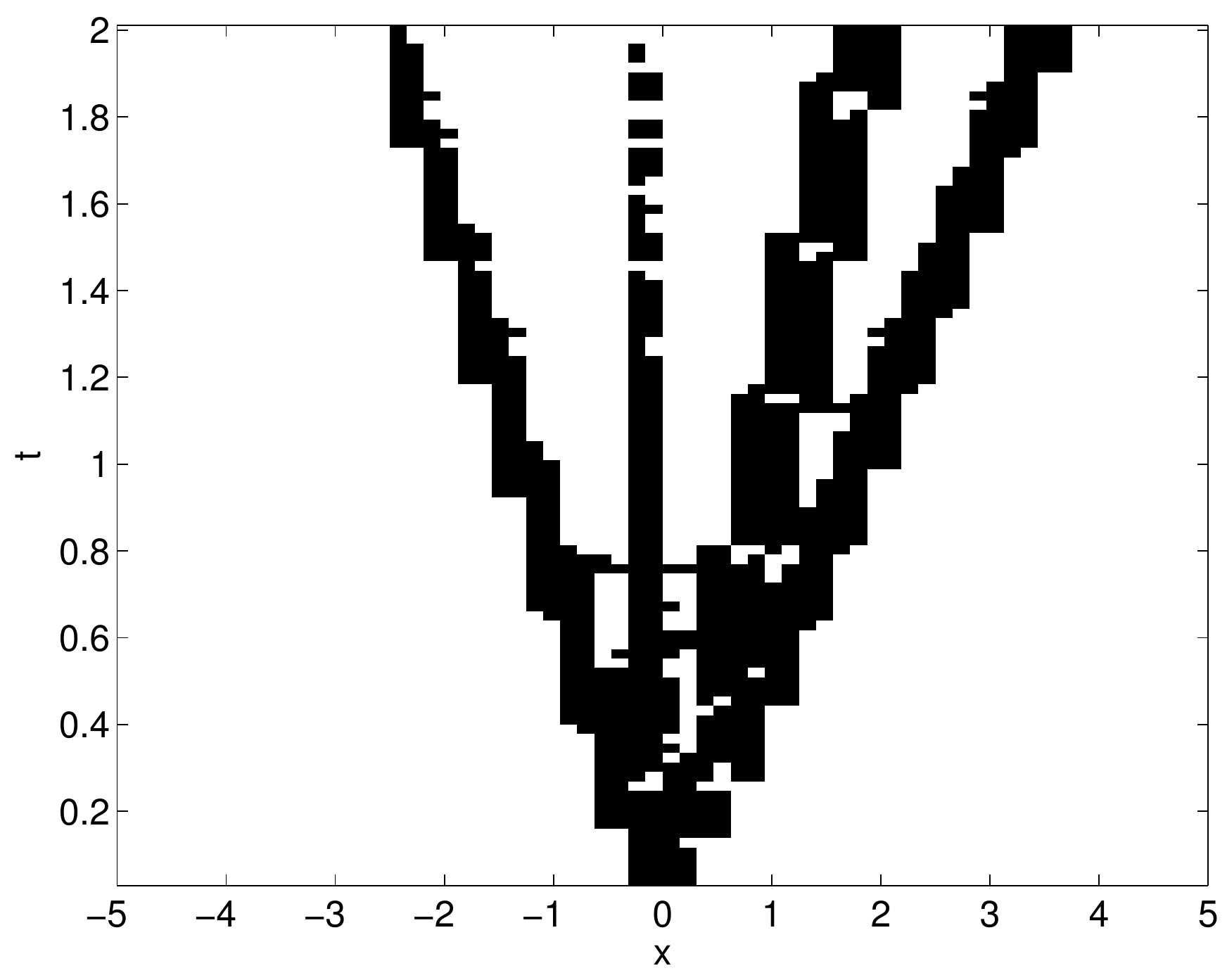}\label{fig:SodC01k1}} \\
\vspace{-0.2cm}
 \subfigure[$C=0.9$]{\includegraphics[scale = 0.22]{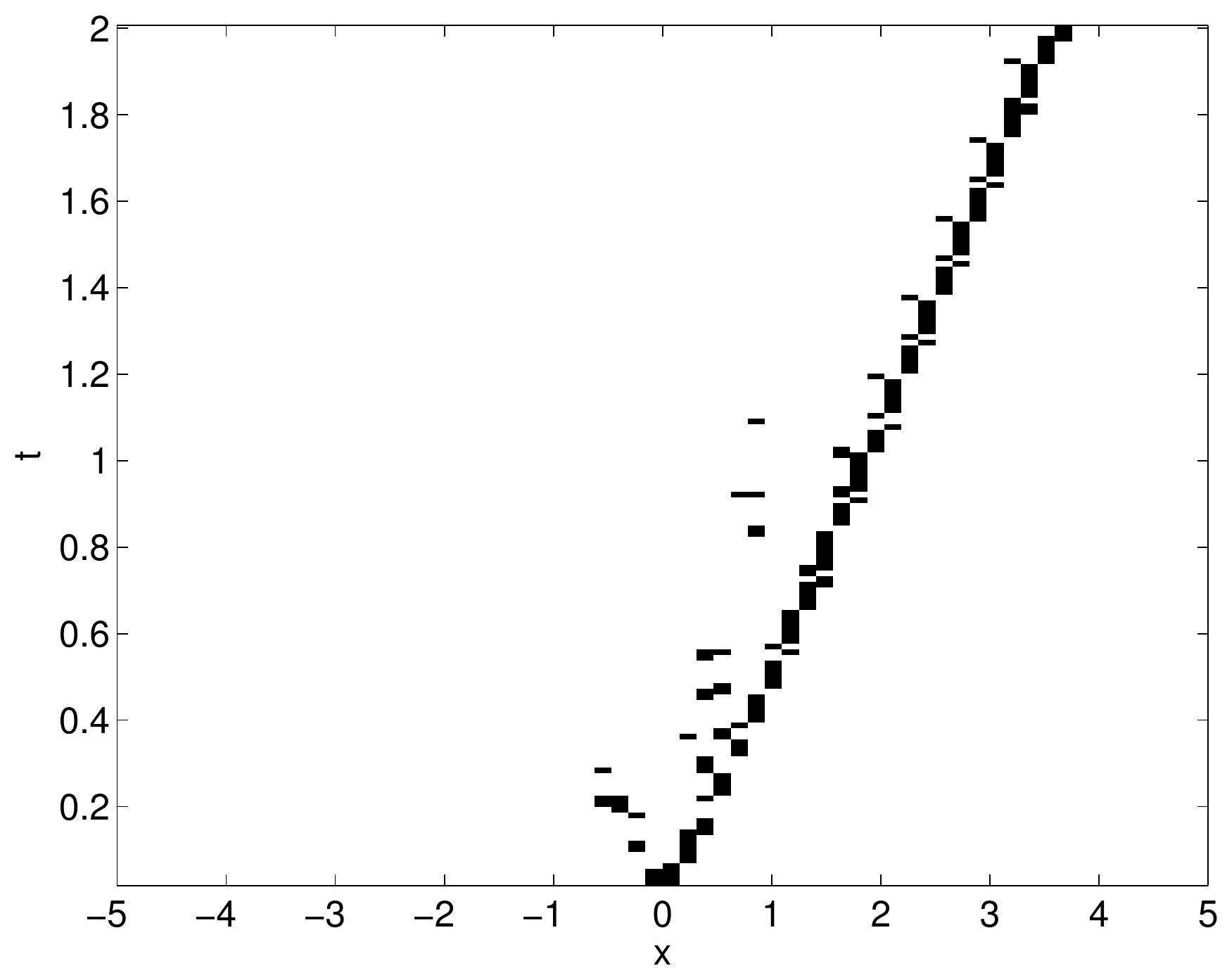}} 
 \subfigure[$C=0.5$]{\includegraphics[scale = 0.22]{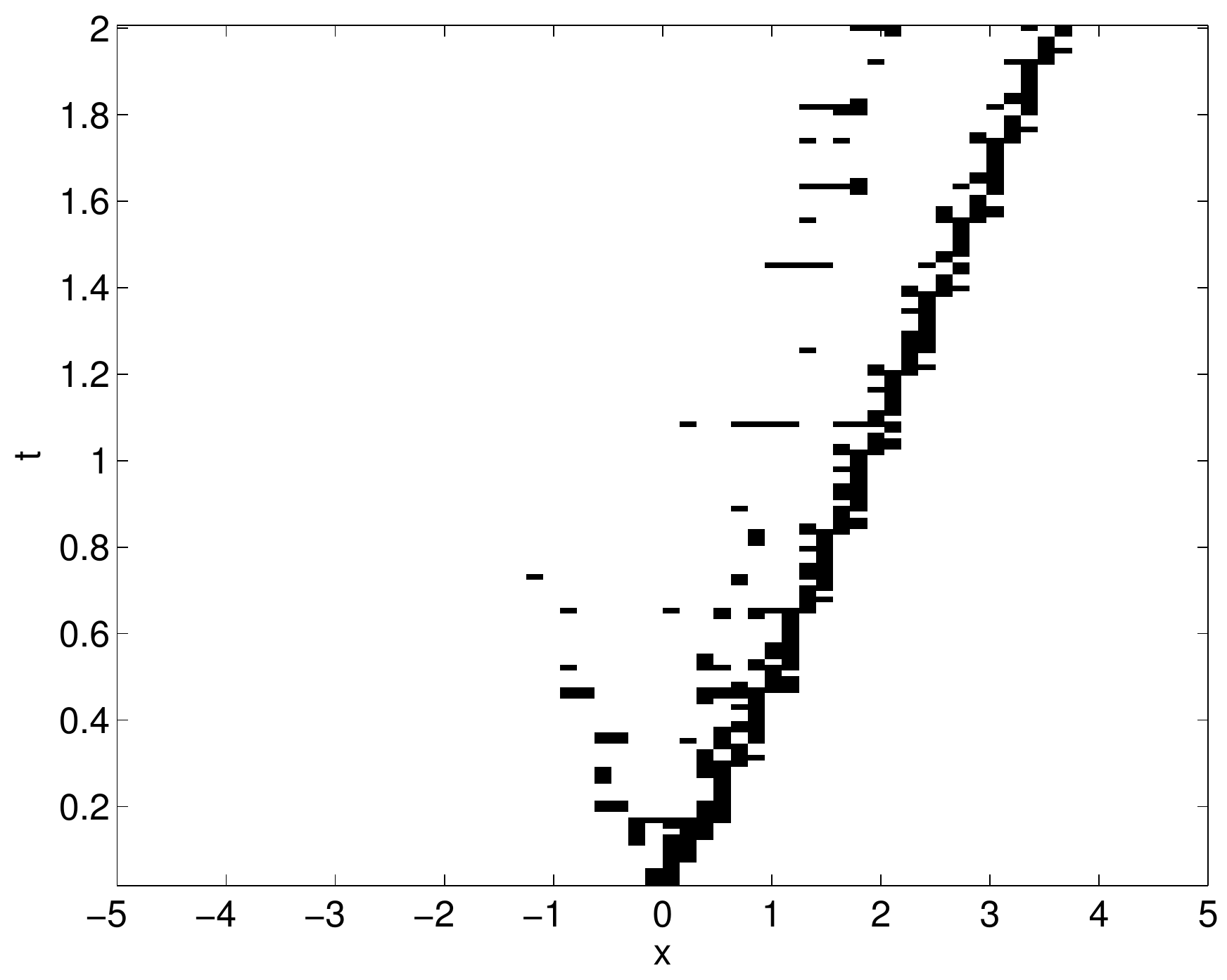}} 
 \subfigure[$C=0.1$]{\includegraphics[scale = 0.22]{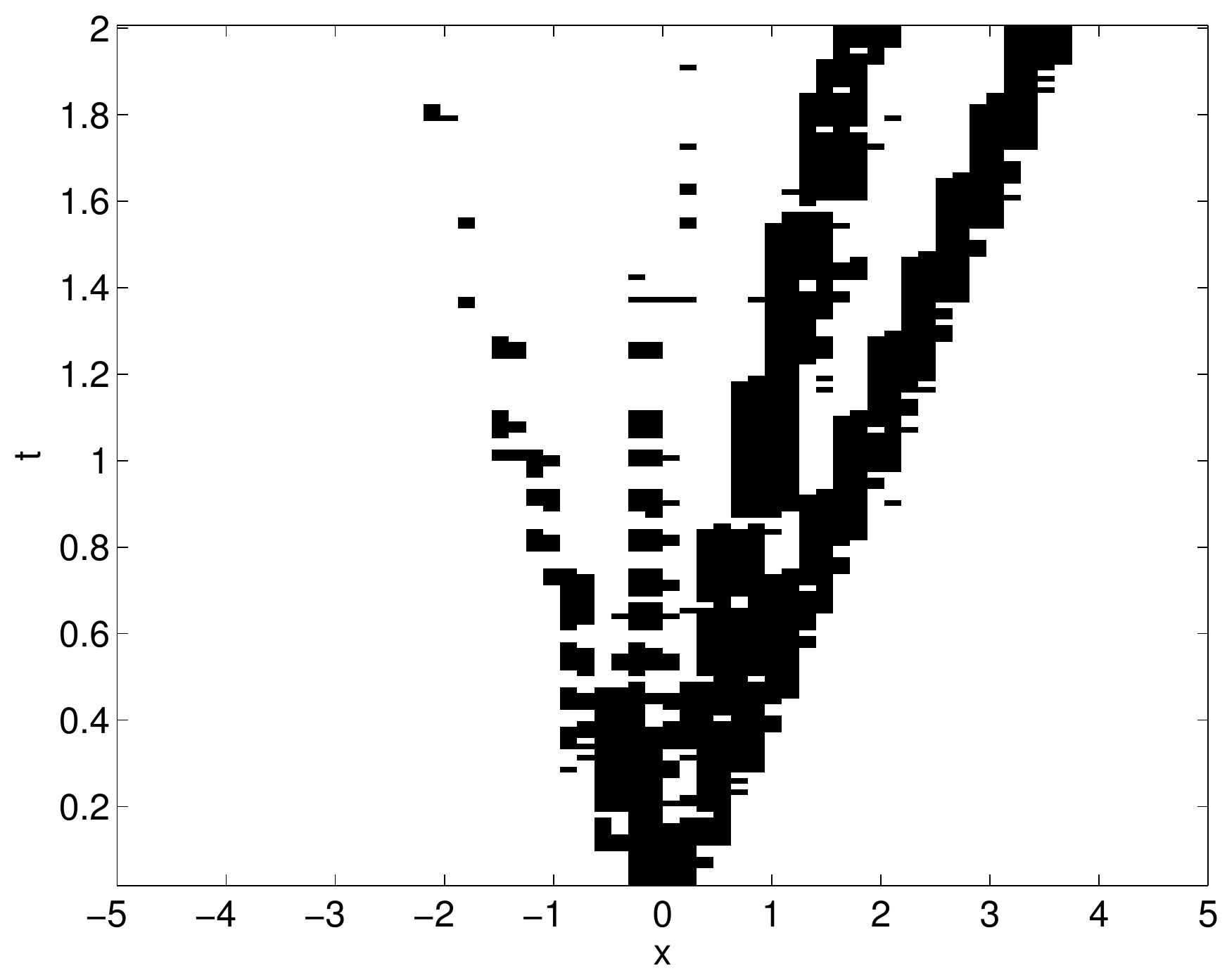}} \\
\vspace{-0.2cm}
\caption{Time history plot of detected troubled cells, multiwavelet troubled-cell indicator (density), Sod, 64 elements. First row: $k=1$, second row: $k=2$.}\label{fig:SodC}
\end{figure}

\begin{figure}[h!]
\centering
 \subfigure[$C=0.9$]{\includegraphics[scale = 0.22]{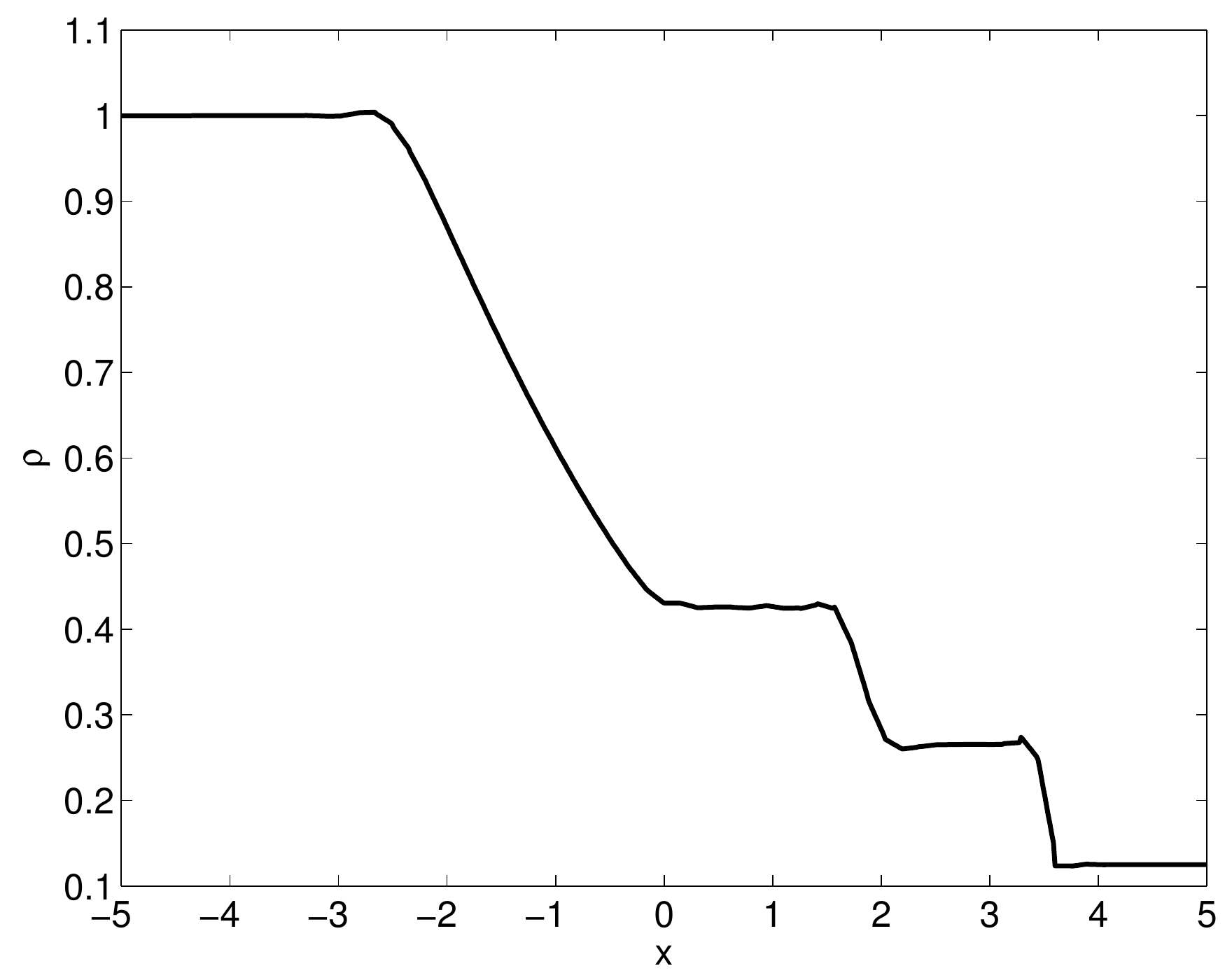}}
 \subfigure[$C=0.5$]{\includegraphics[scale = 0.22]{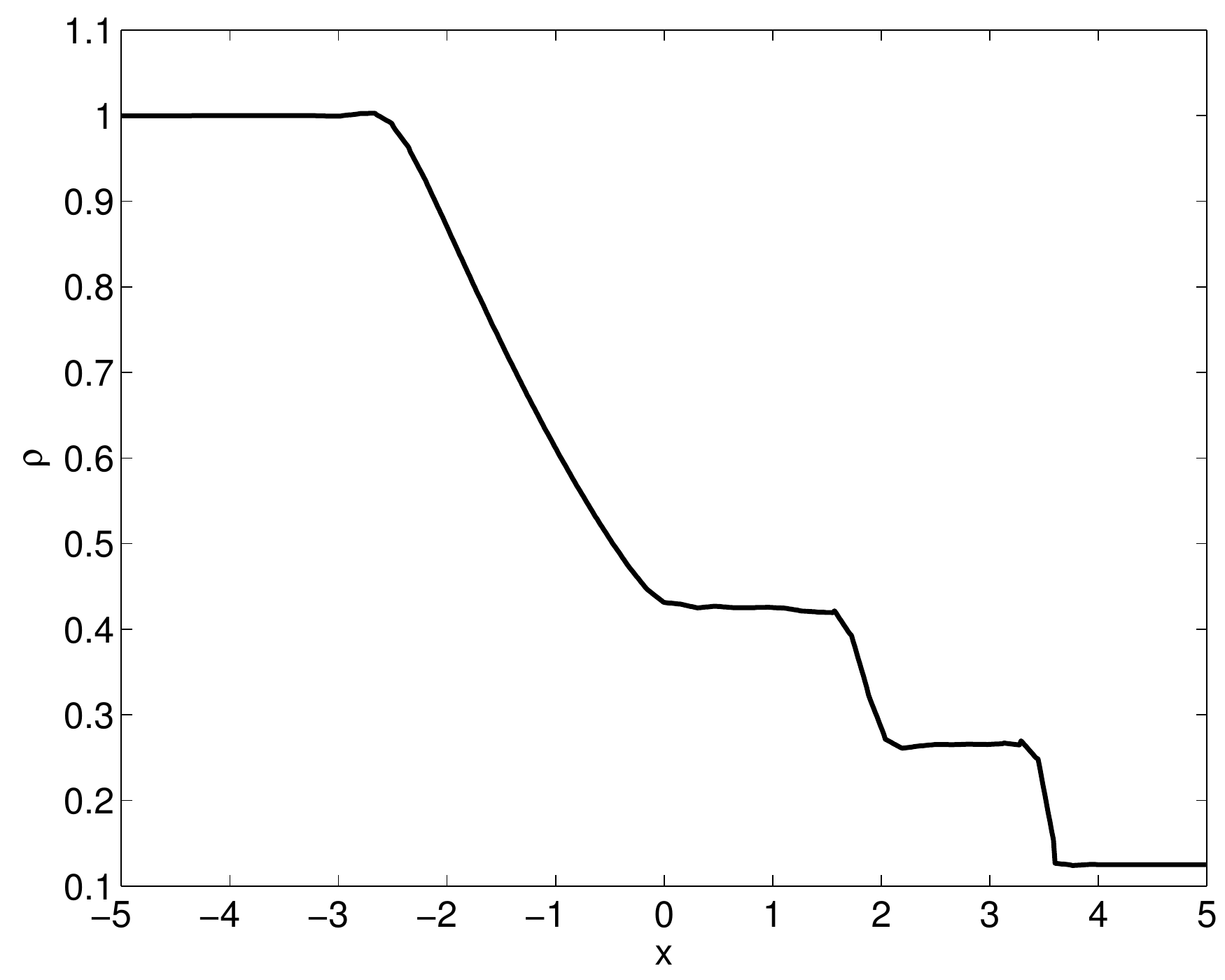}}
 \subfigure[$C=0.1$]{\includegraphics[scale = 0.22]{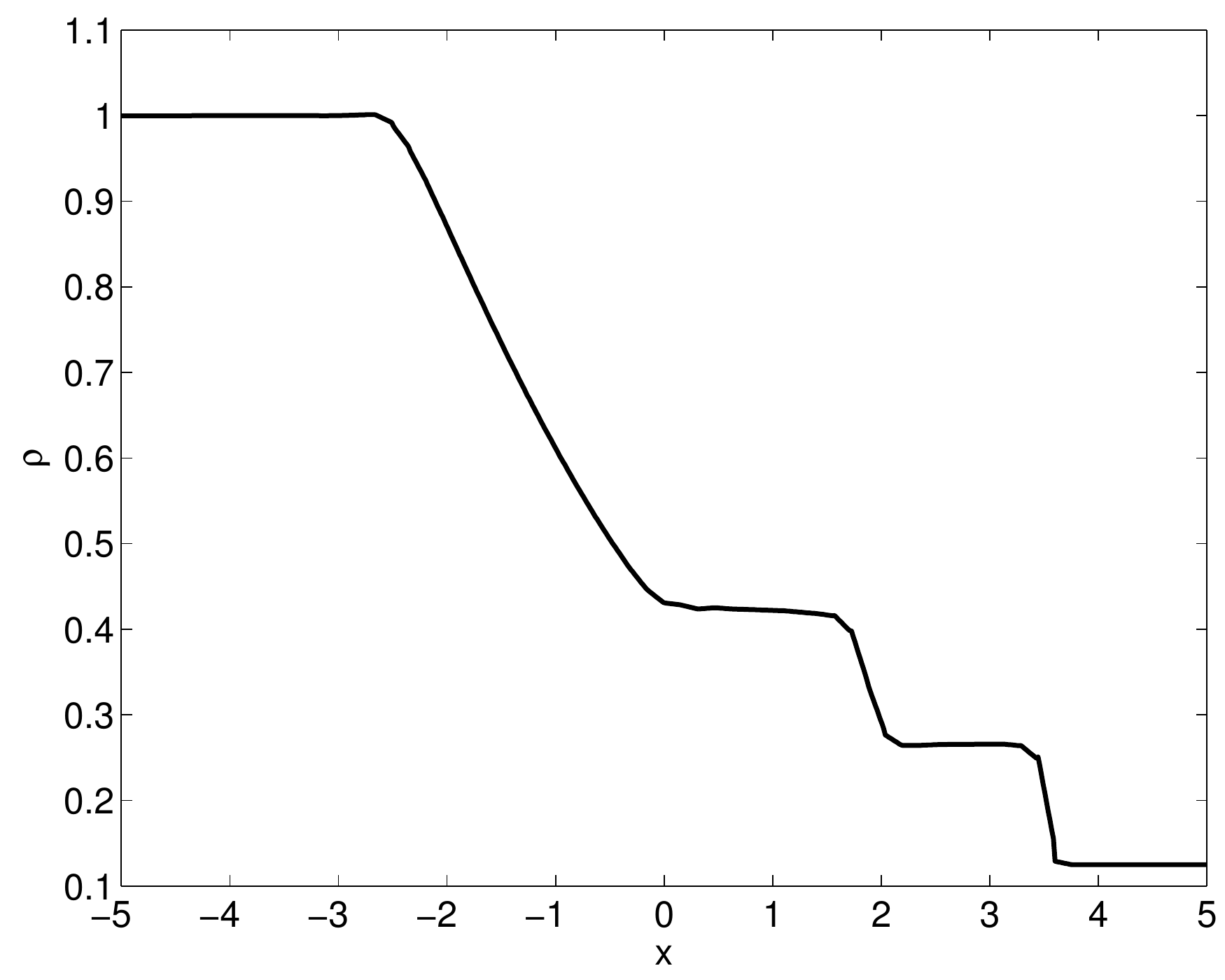}} \\
\vspace{-0.2cm}
 \subfigure[$C=0.9$]{\includegraphics[scale = 0.22]{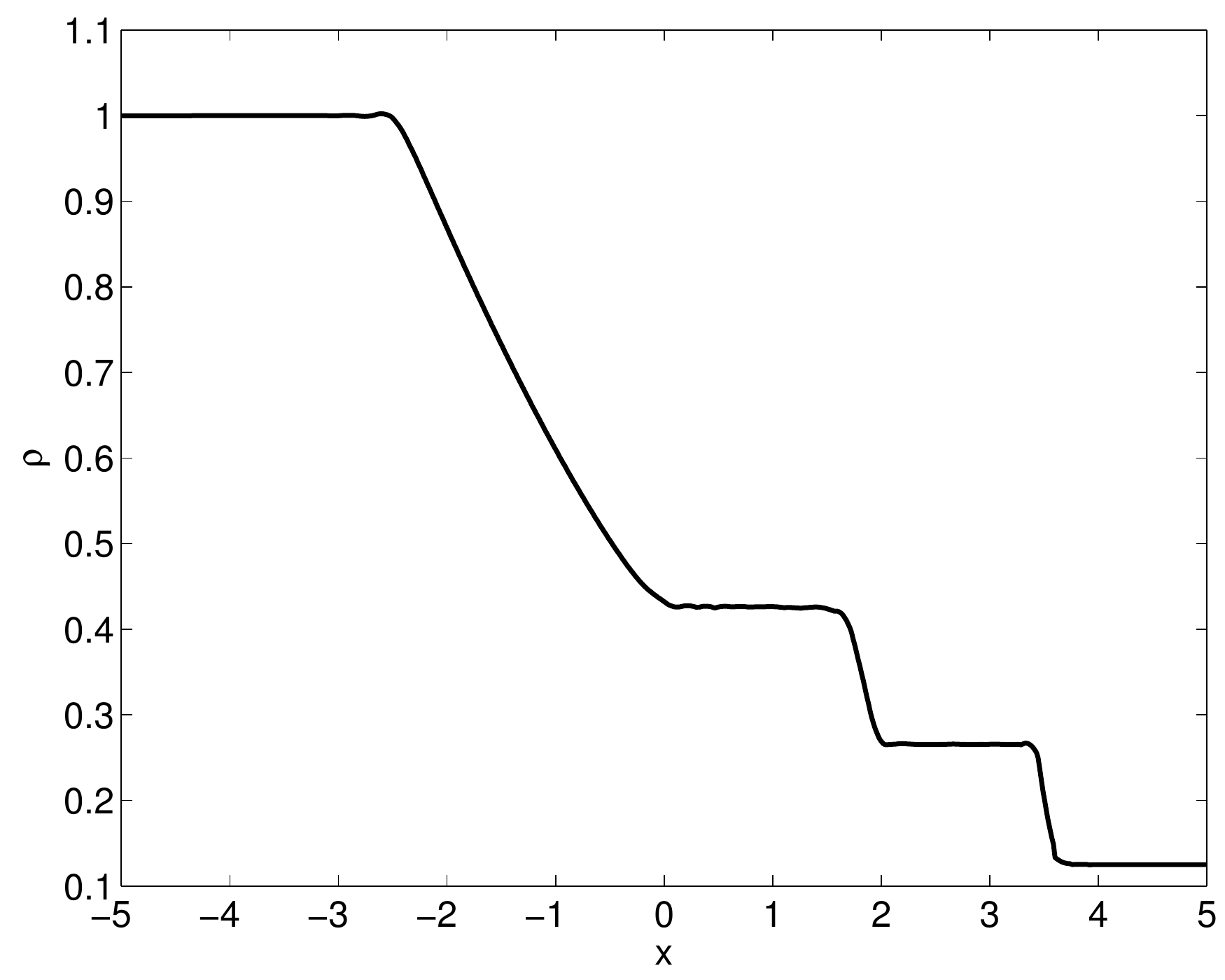}} 
 \subfigure[$C=0.5$]{\includegraphics[scale = 0.22]{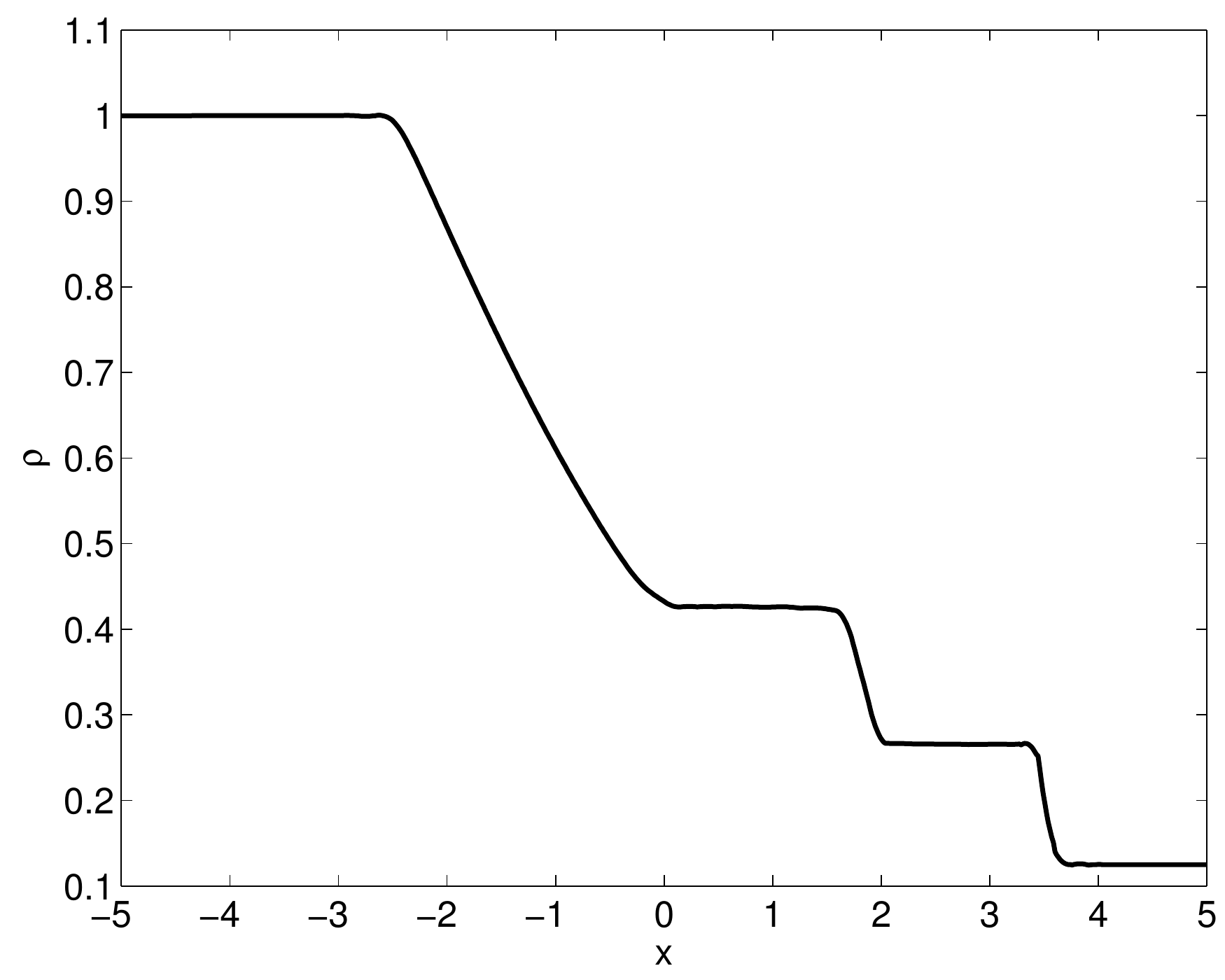}} 
 \subfigure[$C=0.1$]{\includegraphics[scale = 0.22]{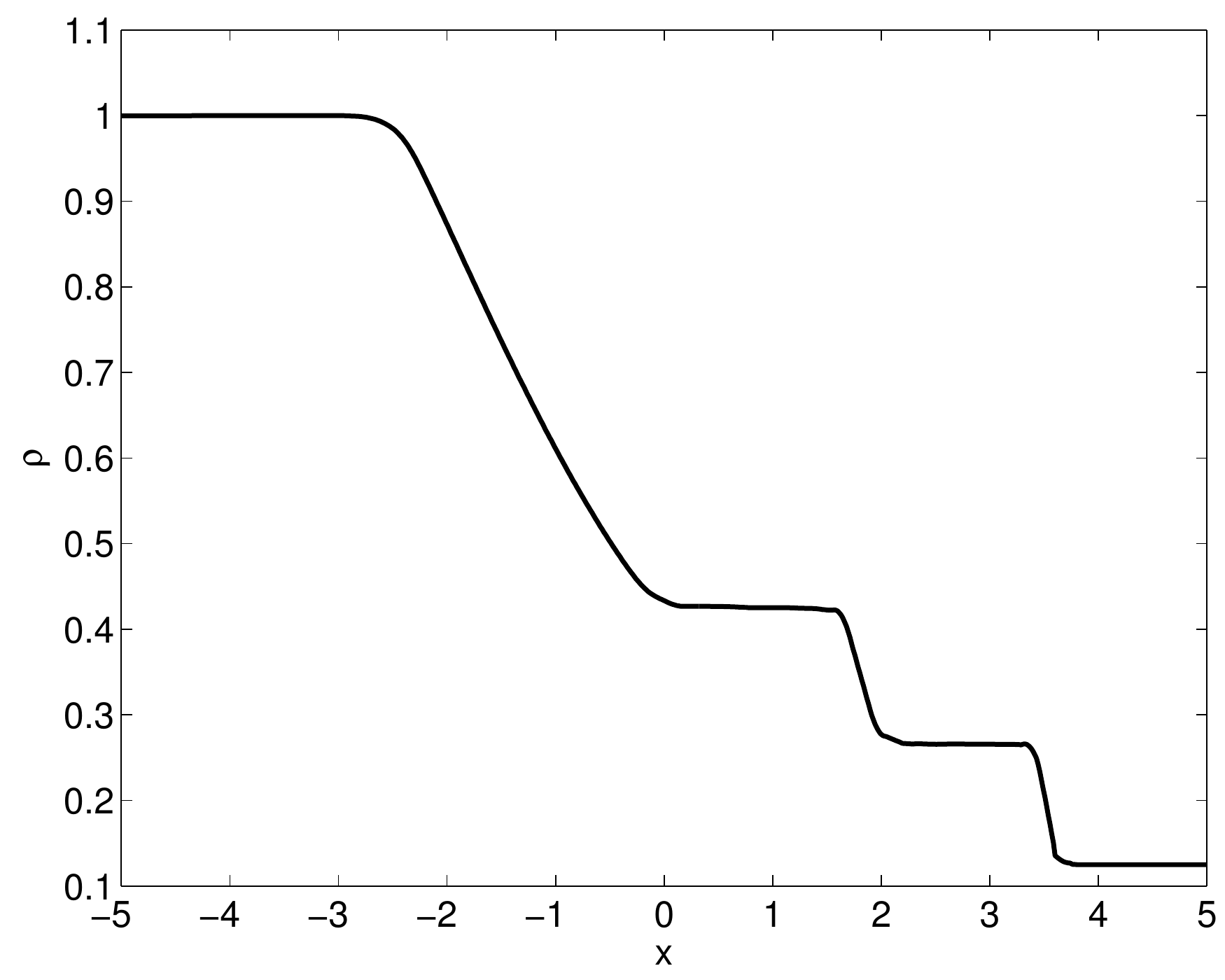}} \\
\vspace{-0.2cm}
\caption{Computed approximation at $T=2$, multiwavelet troubled-cell indicator (density), Sod, 64 elements. First row: $k=1$, second row: $k=2$.}\label{fig:SodCsol}
\end{figure}

The detected troubled cells using the multiwavelet indicator (on density) for different values of $C$ are shown in Figure \ref{fig:SodC}, both using linear and quadratic approximations. The corresponding approximate solutions at $T=2$ are given in Figure \ref{fig:SodCsol}. For $k=1$, it is clearly visible that $C=0.9$ only selects the strongest shock and part of the rarefaction wave. Therefore, we need to decrease the value of $C$, which makes the indicator more strict. For $C=0.1$, the shock, contact discontinuity and the end points of the rarefaction wave (where the derivative of the approximation is discontinuous) are detected. This means that our indicator is very accurate if the value of $C$ is chosen properly.  If not only density, but also entropy is used in our multiwavelet troubled-cell indicator, exactly the same elements are detected as troubled cells. This behavior can be seen in Figure \ref{fig:Soddensityentropy} for Sod's shock tube, and is generally true for each test problem that we investigated.  For $k = 2$, the multiwavelet indicator using density detects fewer elements than in the linear case, which means that $C$ should be smaller than in the piecewise linear case in order to select the same regions. 

\begin{figure}[h!]
 \centering
\subfigure[KXRCF, $k=1$]{\includegraphics[scale = 0.28]{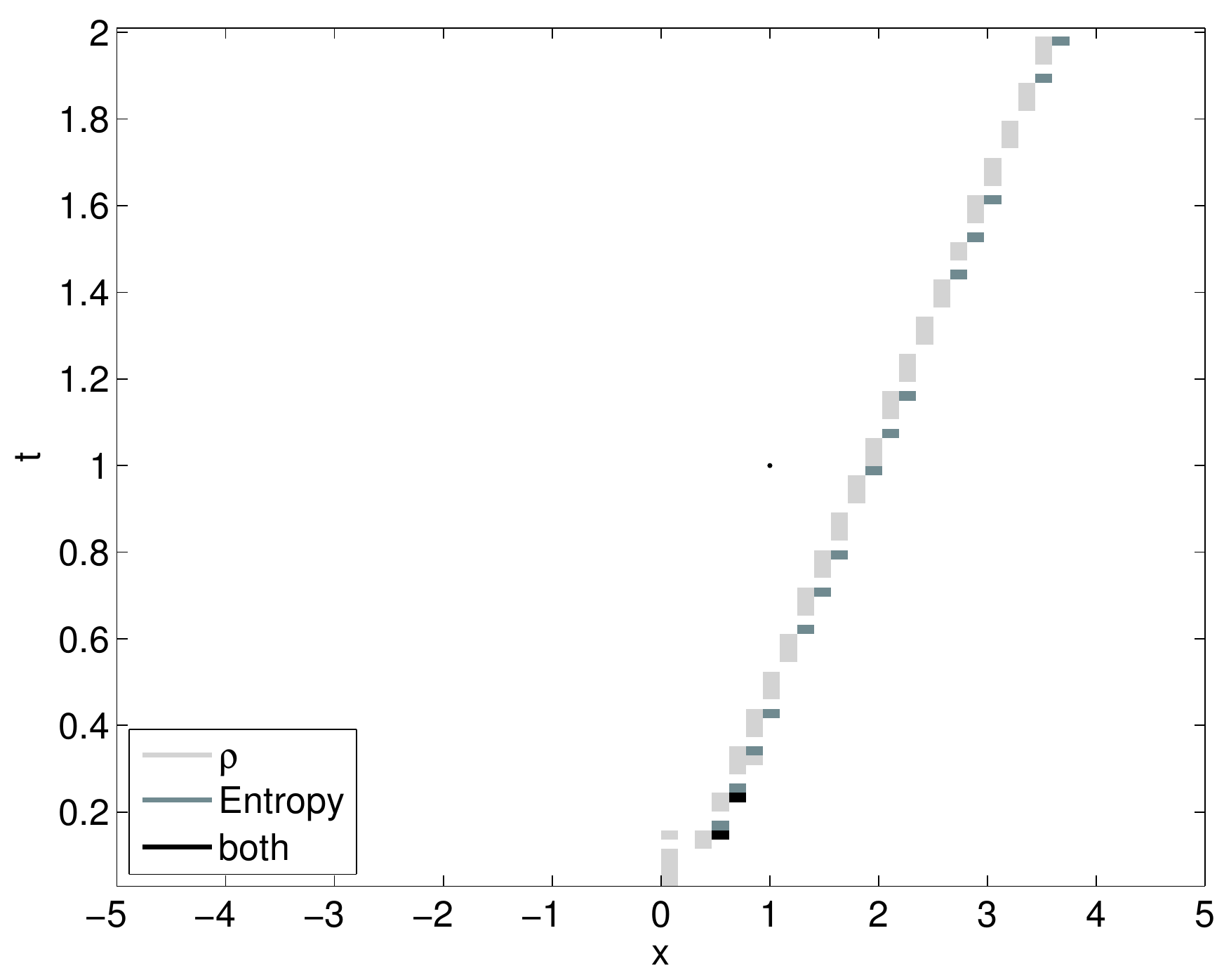}}
\subfigure[KXRCF, $k=2$]{\includegraphics[scale = 0.28]{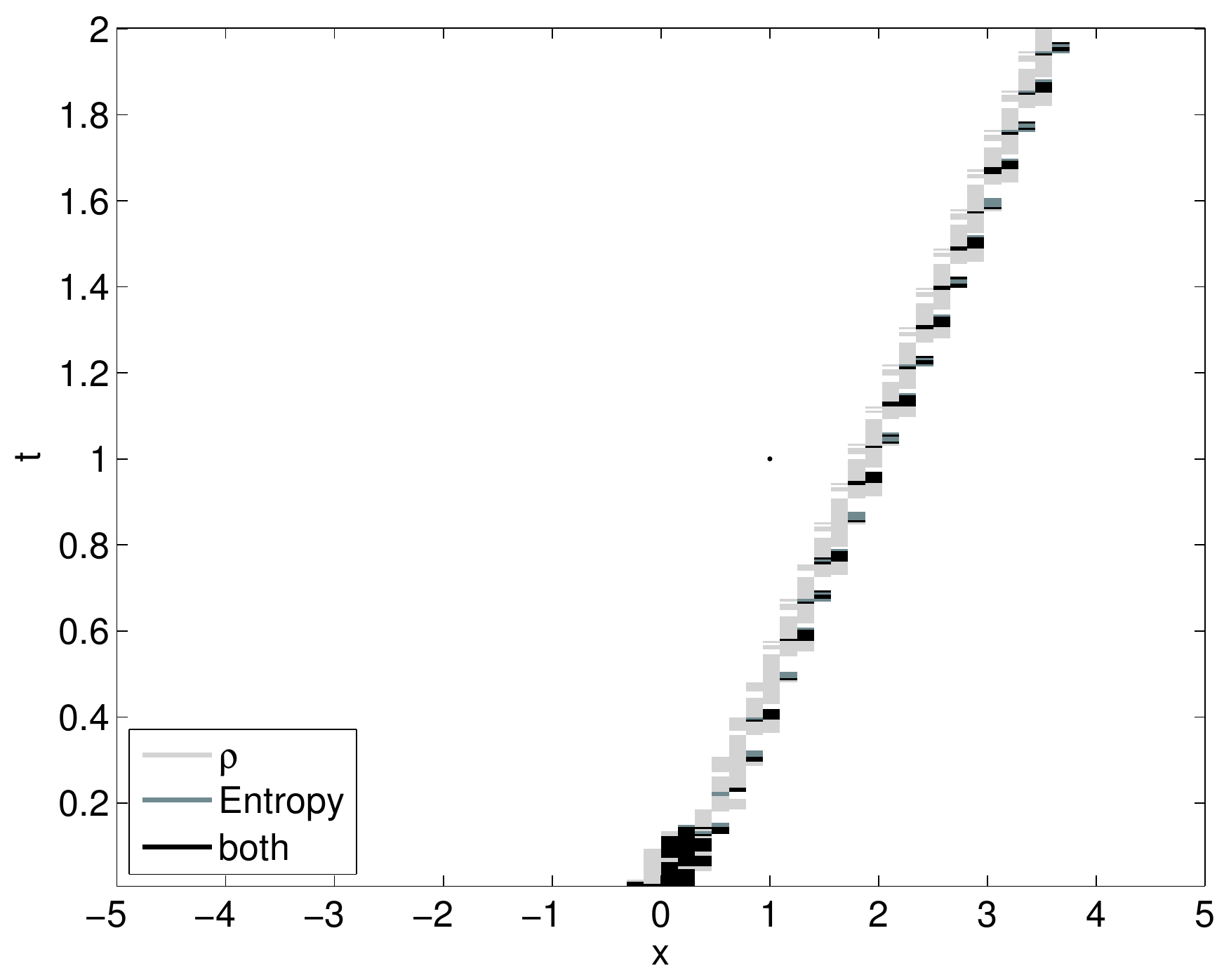}} \\
\vspace{-0.2cm}
\subfigure[Harten, $\alpha = 1.5$, $k=1$]{\includegraphics[scale = 0.28]{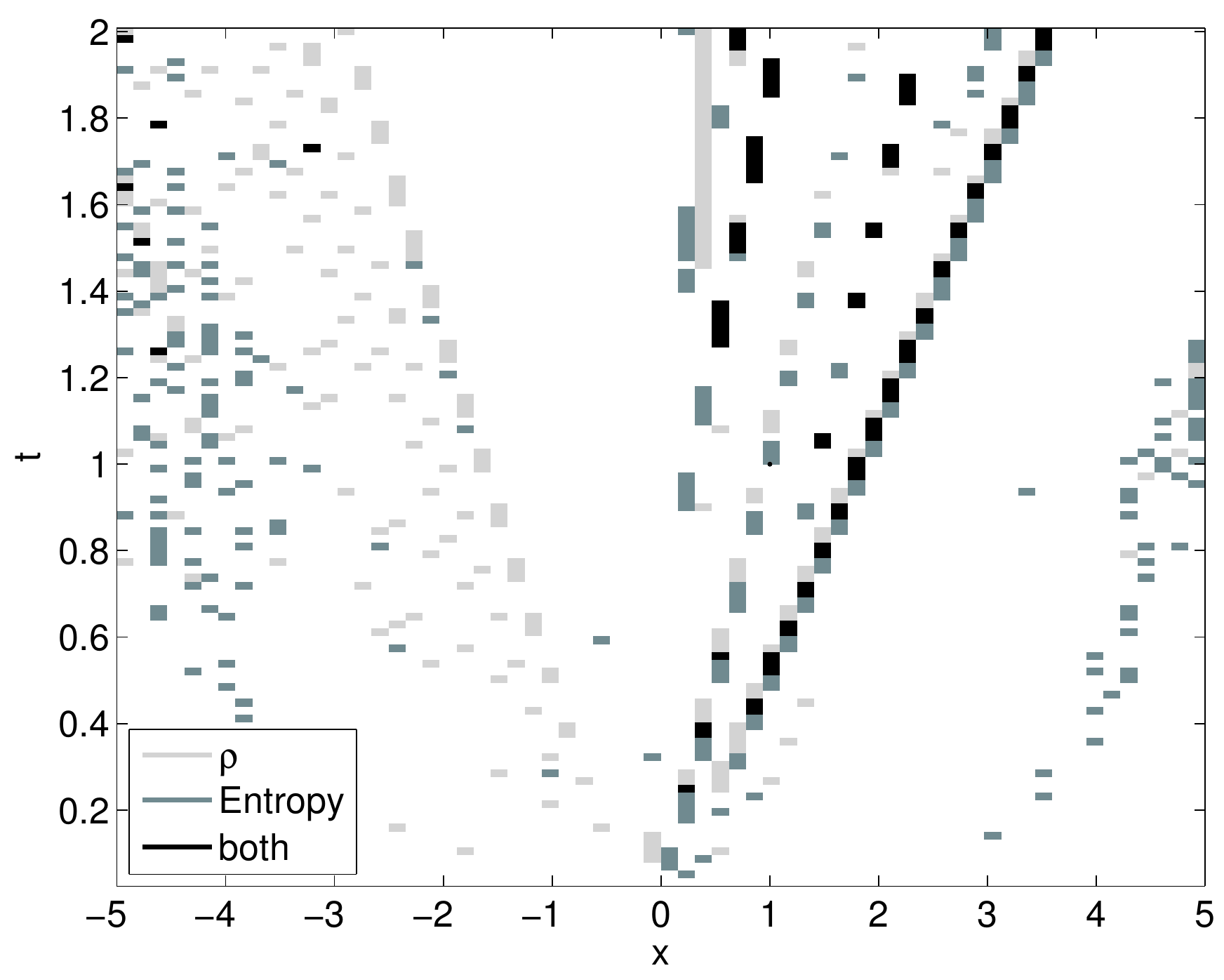}}
\subfigure[Harten, $\alpha = 1.5$, $k=2$]{\includegraphics[scale = 0.28]{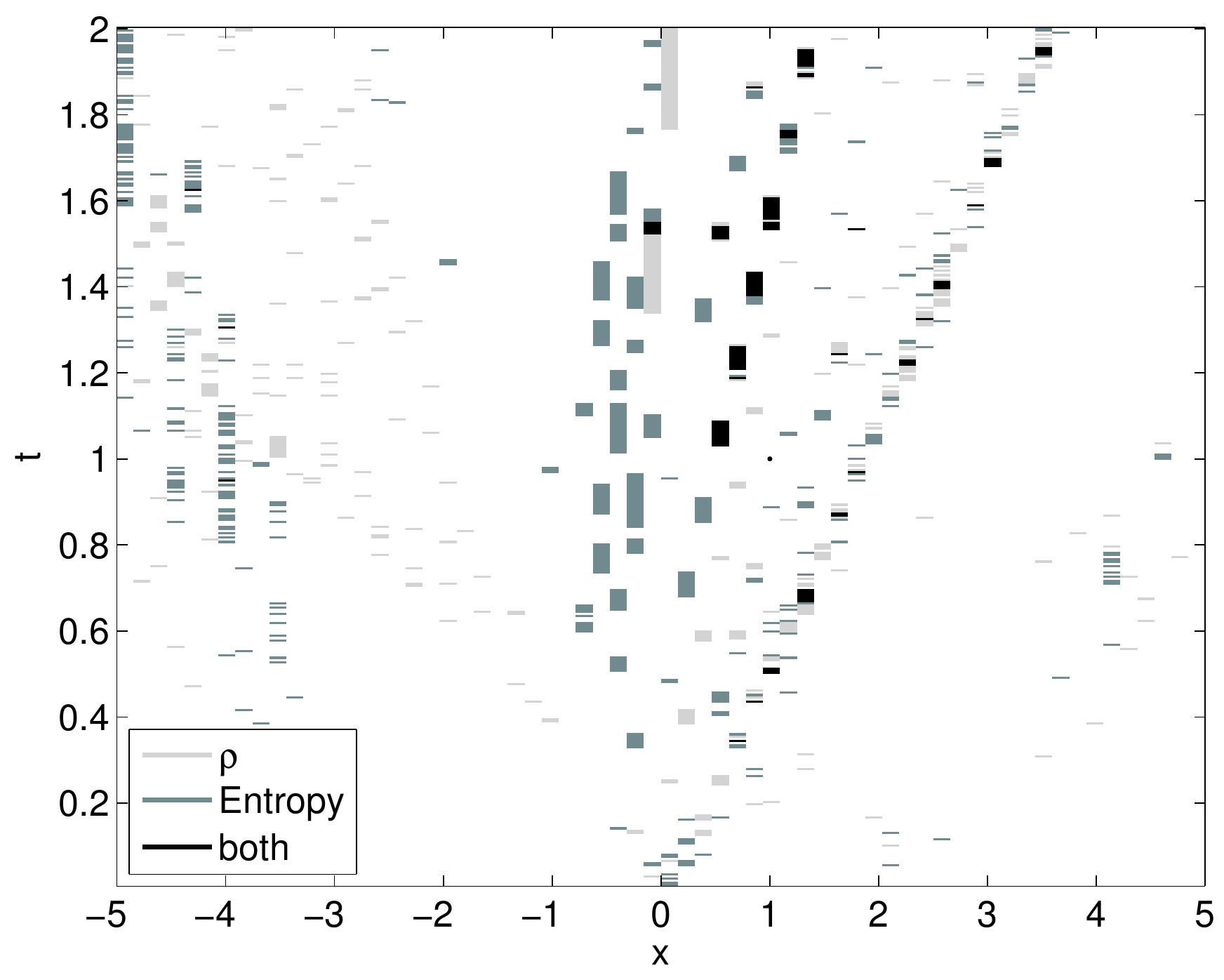}} \\
\vspace{-0.2cm}
\caption{The KXRCF or Harten's troubled-cell indicator on density and entropy, Sod, 64 elements.}\label{fig:SodKH}
\end{figure}

\begin{figure}[ht!]
 \centering
\subfigure[KXRCF, $k=1$]{\includegraphics[scale = 0.28]{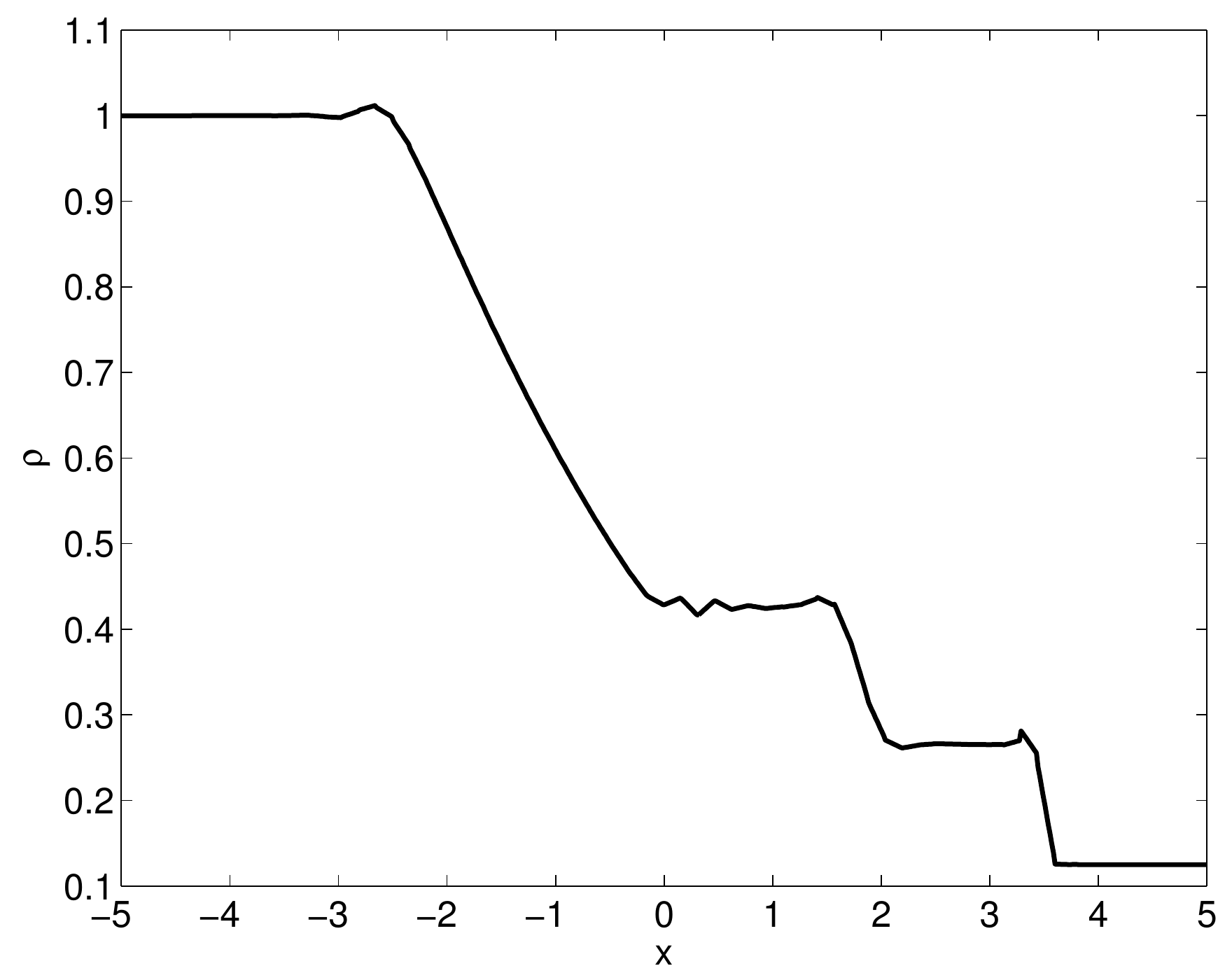}}
\subfigure[KXRCF, $k=2$]{\includegraphics[scale = 0.28]{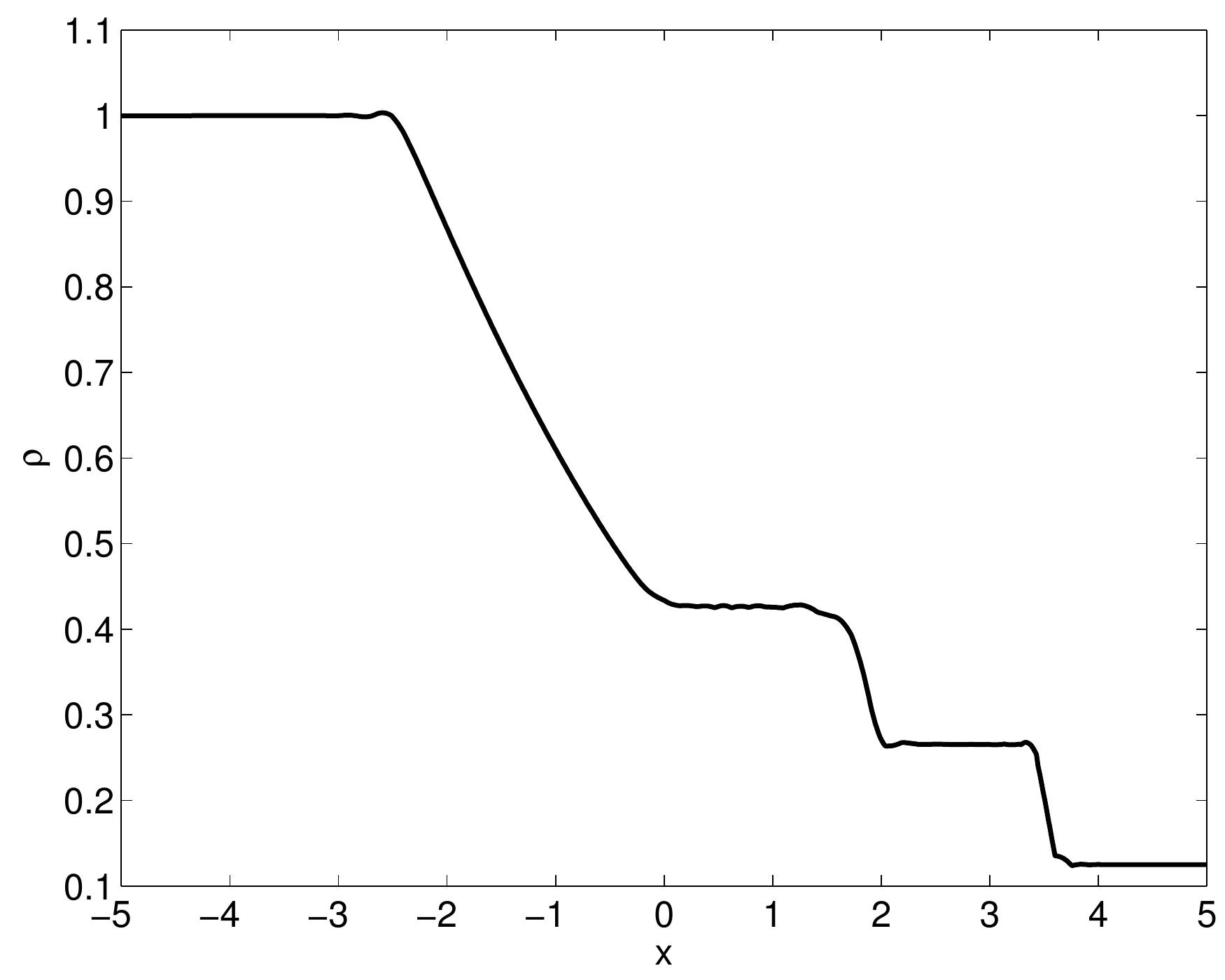}} \\
\subfigure[Harten, $\alpha = 1.5$, $k=1$]{\includegraphics[scale = 0.28]{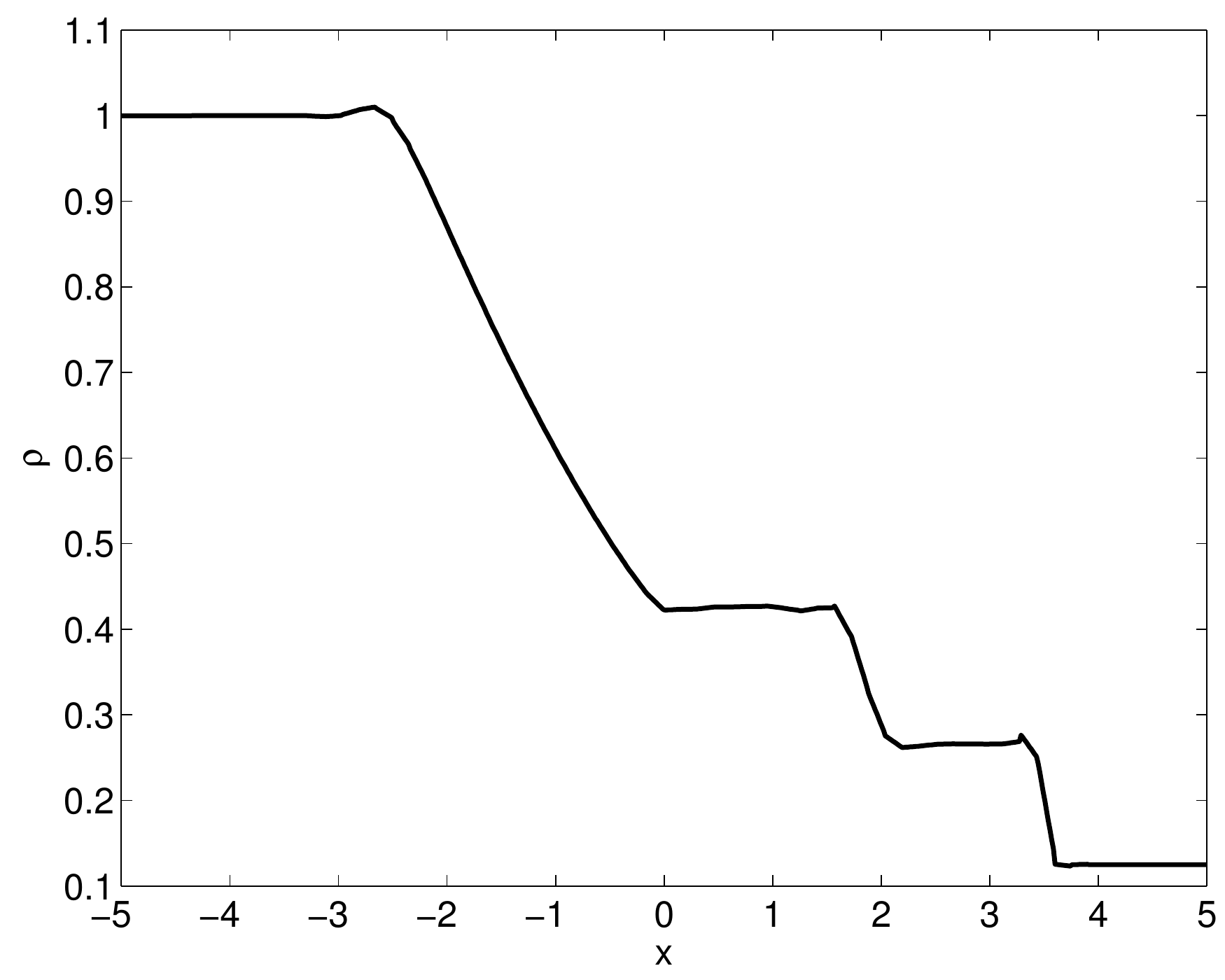}}
\subfigure[Harten, $\alpha = 1.5$, $k=2$]{\includegraphics[scale = 0.28]{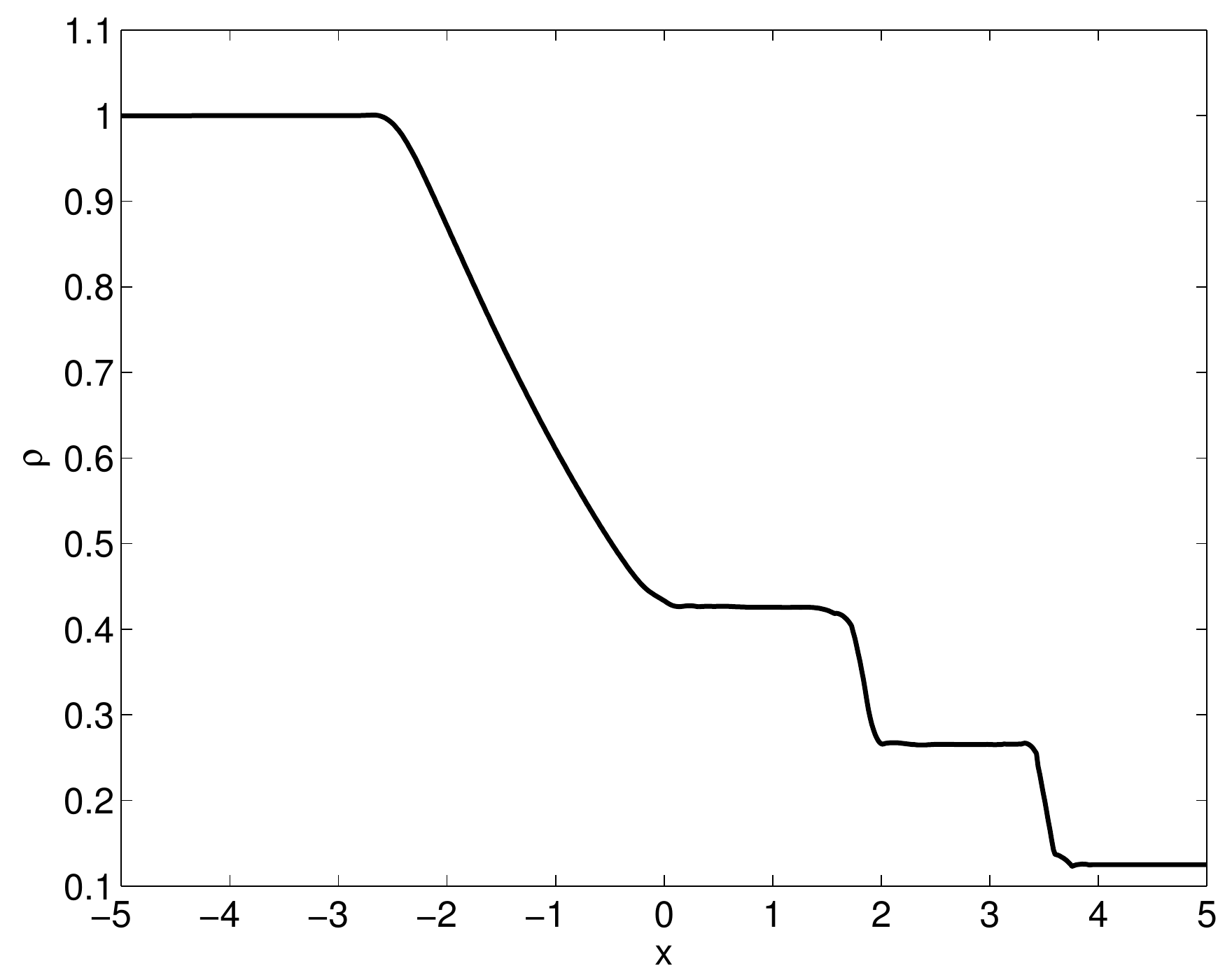}}
\caption{Computed approximation at $T=2$, KXRCF or Harten's troubled-cell indicator on density and entropy, Sod, 64 elements.}\label{fig:SodKHsol}
\end{figure}

The KXRCF and Harten ($\alpha = 1.5$) results using density and entropy as an indicator variable are visualized in Figures \ref{fig:SodKH} and \ref{fig:SodKHsol}. It is very surprising to notice that the KXRCF indicator is able to
detect the shock, but the contact discontinuity is not found (see Figure \ref{fig:SodKH}), such that the resulting approximation is very oscillatory. This was also noted in \cite{Kri04XRCF}.  Taking $k=2$ improves the solution, 
but still does not find the contact discontinuities. Using Harten's troubled-cell indicator, the detected elements are more scattered over the domain.
 
\newpage

\clearpage

\subsubsection{Lax's shock tube}
The second test problem that we consider is the shock tube problem of Lax \cite{Lax54}.  The initial conditions are given by 
\begin{subequations}
\begin{equation}
 \rho(x,0)  = \left\{ \begin{array}{ll} 0.445, & \mbox{for } x <0, \\ 0.5, & \mbox{for } x \geq 0,\end{array} \right.\quad
 p(x,0)  = \left\{ \begin{array}{ll} 3.528, & \mbox{for } x <0, \\ 0.571, & \mbox{for } x \geq 0, \end{array} \right.
\end{equation}
and
\begin{equation}
 u(x,0)  = \left\{ \begin{array}{ll} 0.698, & \mbox{for } x <0, \\ 0, & \mbox{for } x \geq 0, \end{array} \right.
\end{equation}
 and constant initial state boundary conditions are used.
\end{subequations}

\begin{figure}[ht!]
\centering
 \subfigure[$C=0.9$]{\includegraphics[scale = 0.22]{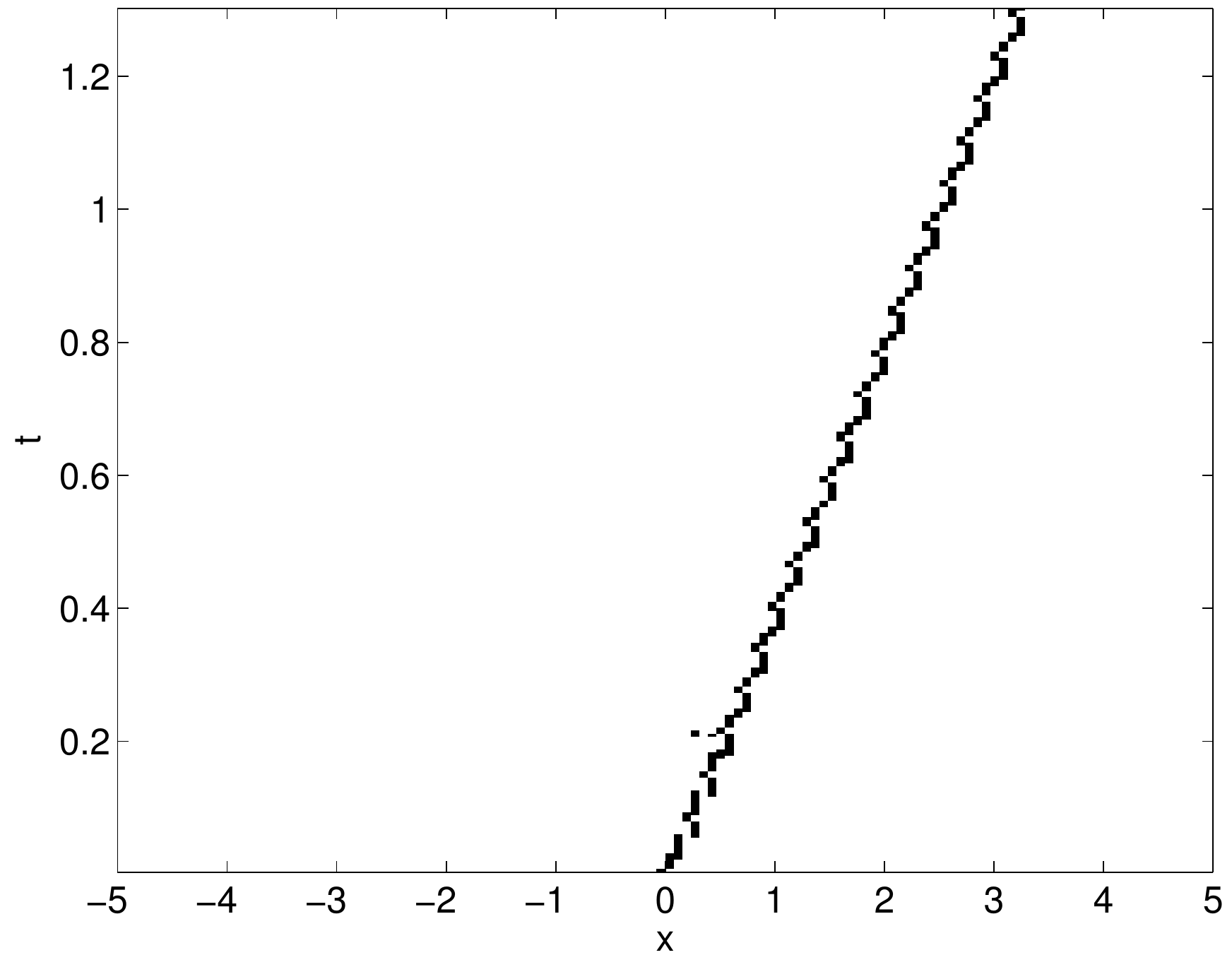}}
 \subfigure[$C=0.5$]{\includegraphics[scale = 0.22]{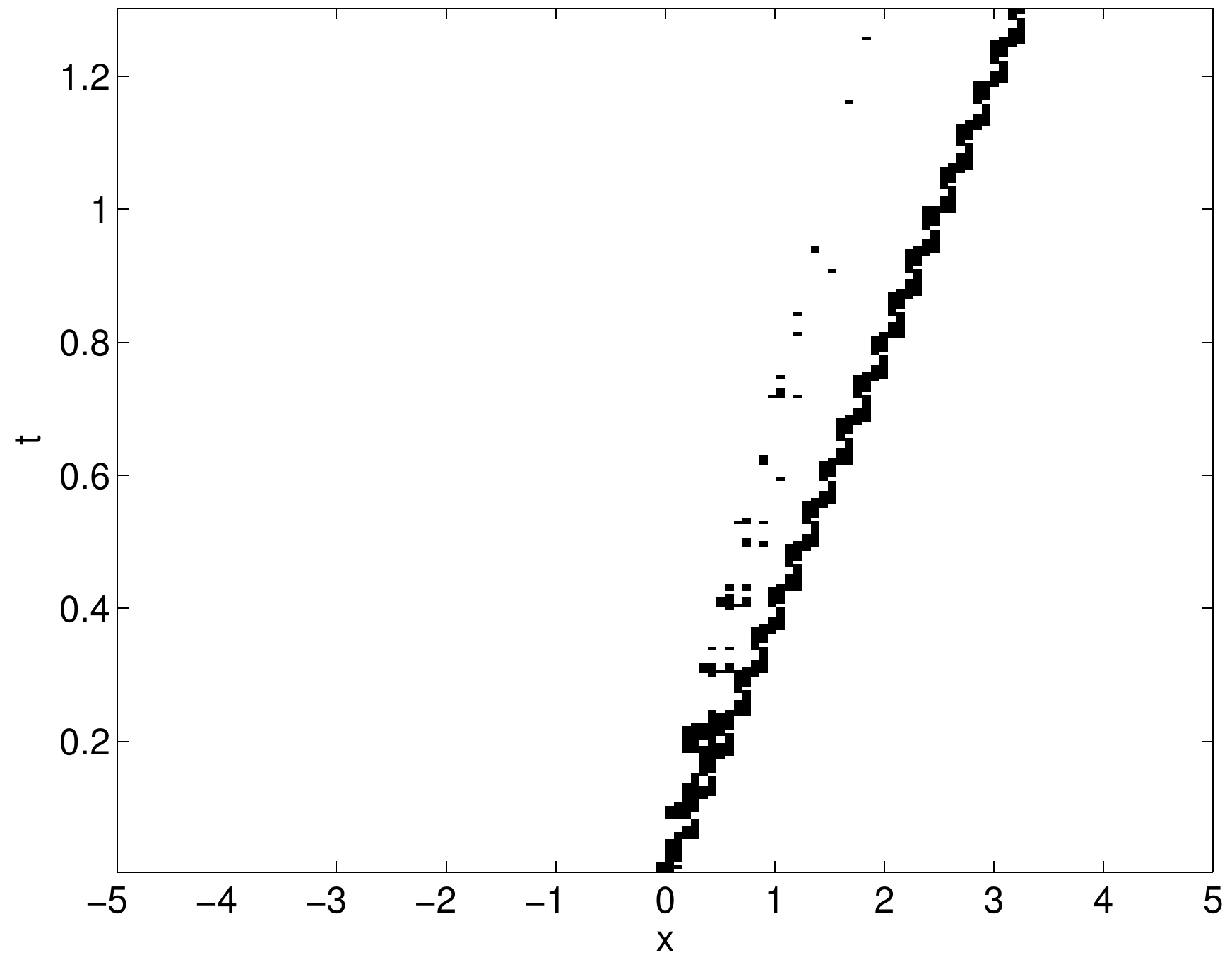}}
 \subfigure[$C=0.1$]{\includegraphics[scale = 0.22]{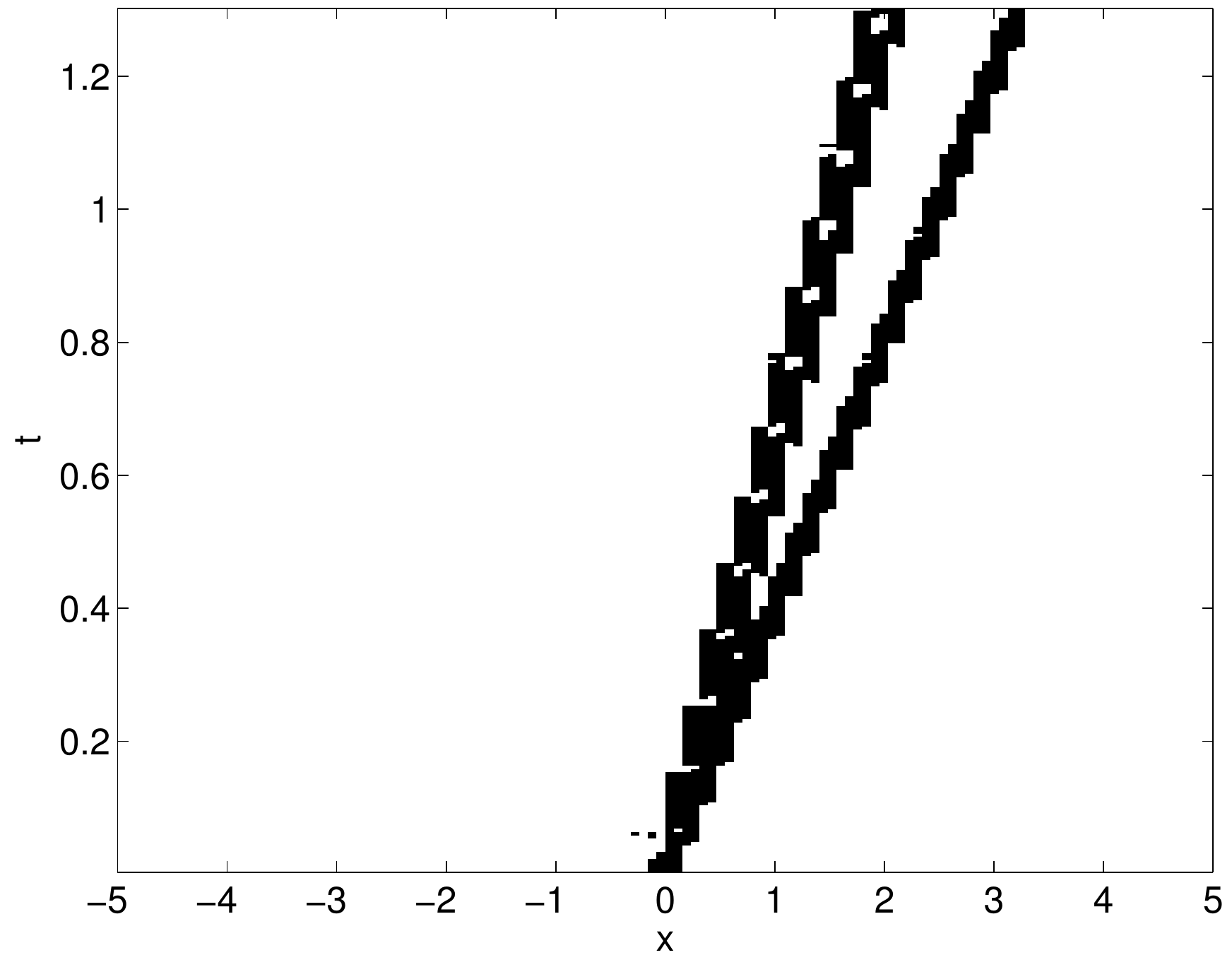}} \\
 \subfigure[$C=0.9$]{\includegraphics[scale = 0.22]{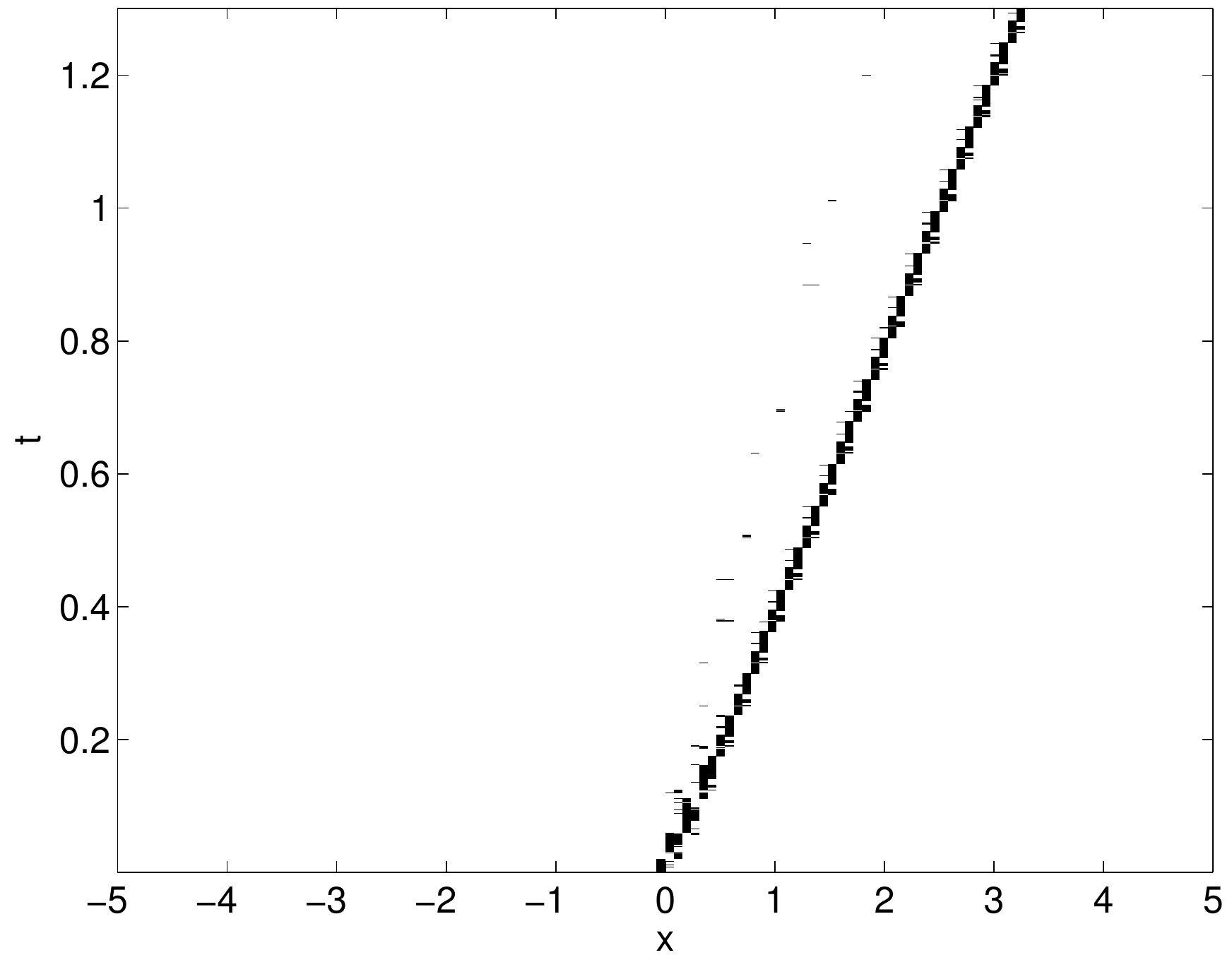}} 
 \subfigure[$C=0.5$]{\includegraphics[scale = 0.22]{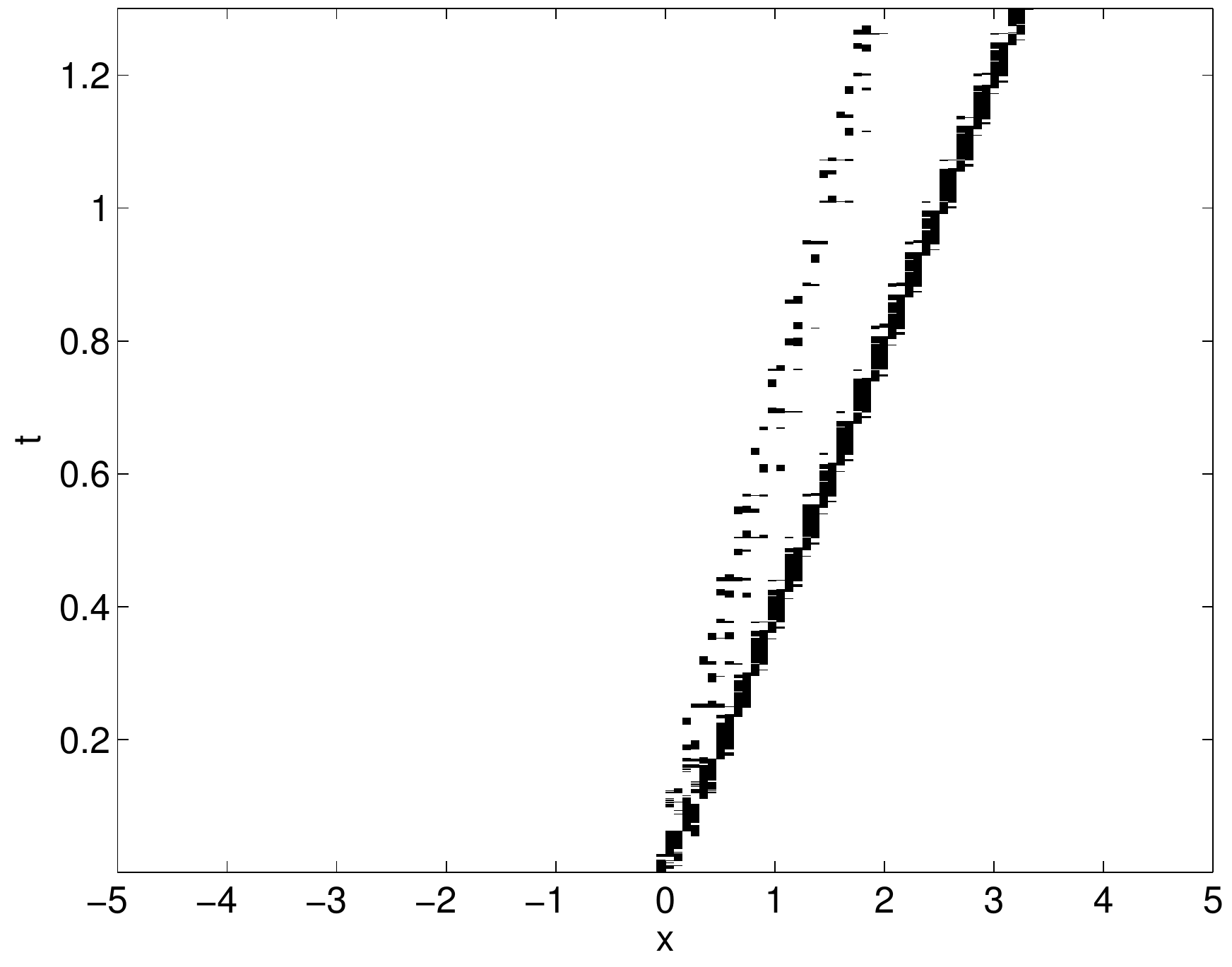}} 
 \subfigure[$C=0.1$]{\includegraphics[scale = 0.22]{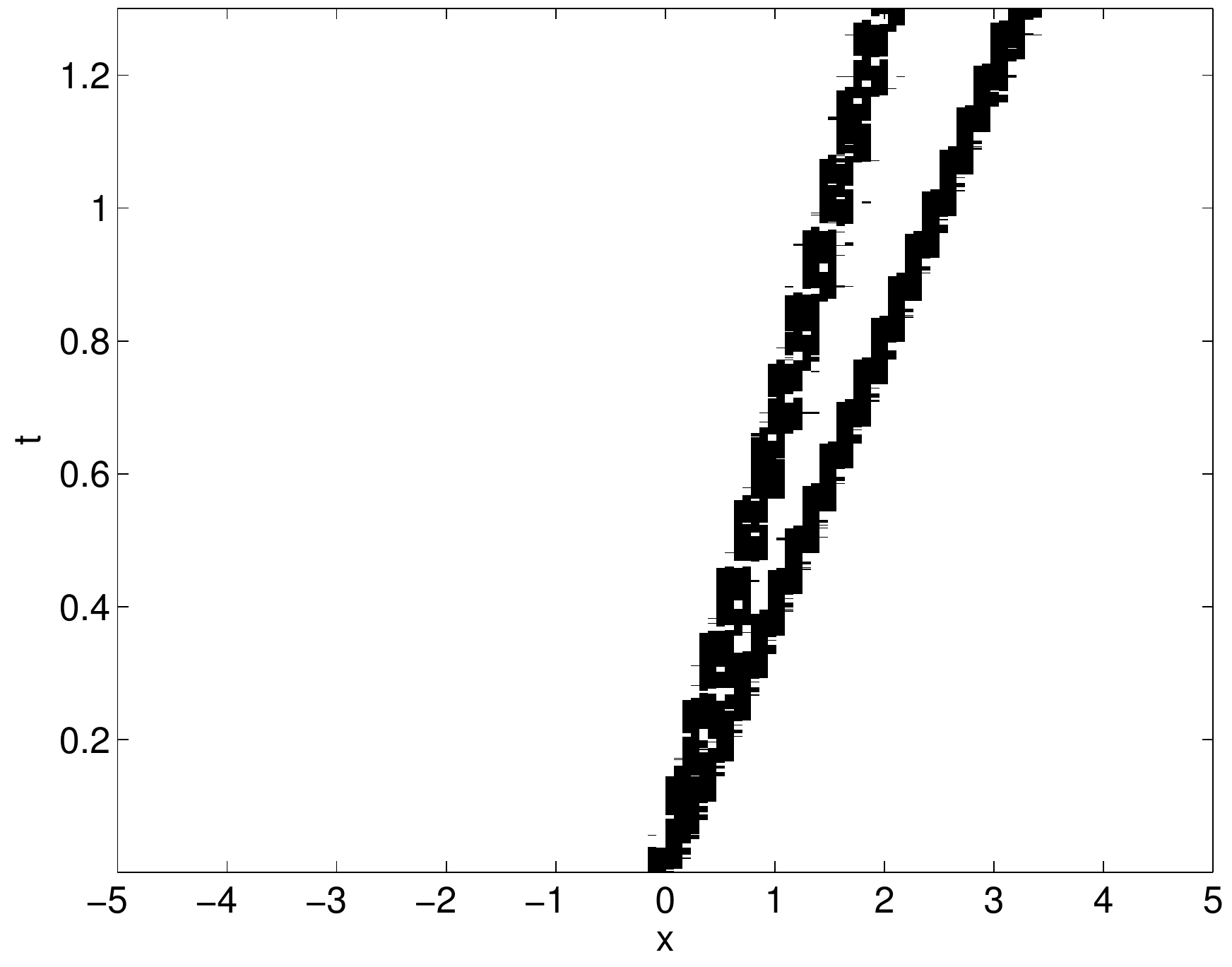}} \\
\vspace{-0.2cm}
\caption{Time history plot of detected troubled cells, multiwavelet shock detector on density, Lax, 128 elements. First row: $k=1$, second row: $k=2$.}\label{fig:LaxC}
\end{figure}

\begin{figure}[ht!]
 \centering
\subfigure[KXRCF, $k=1$]{\includegraphics[scale = 0.29]{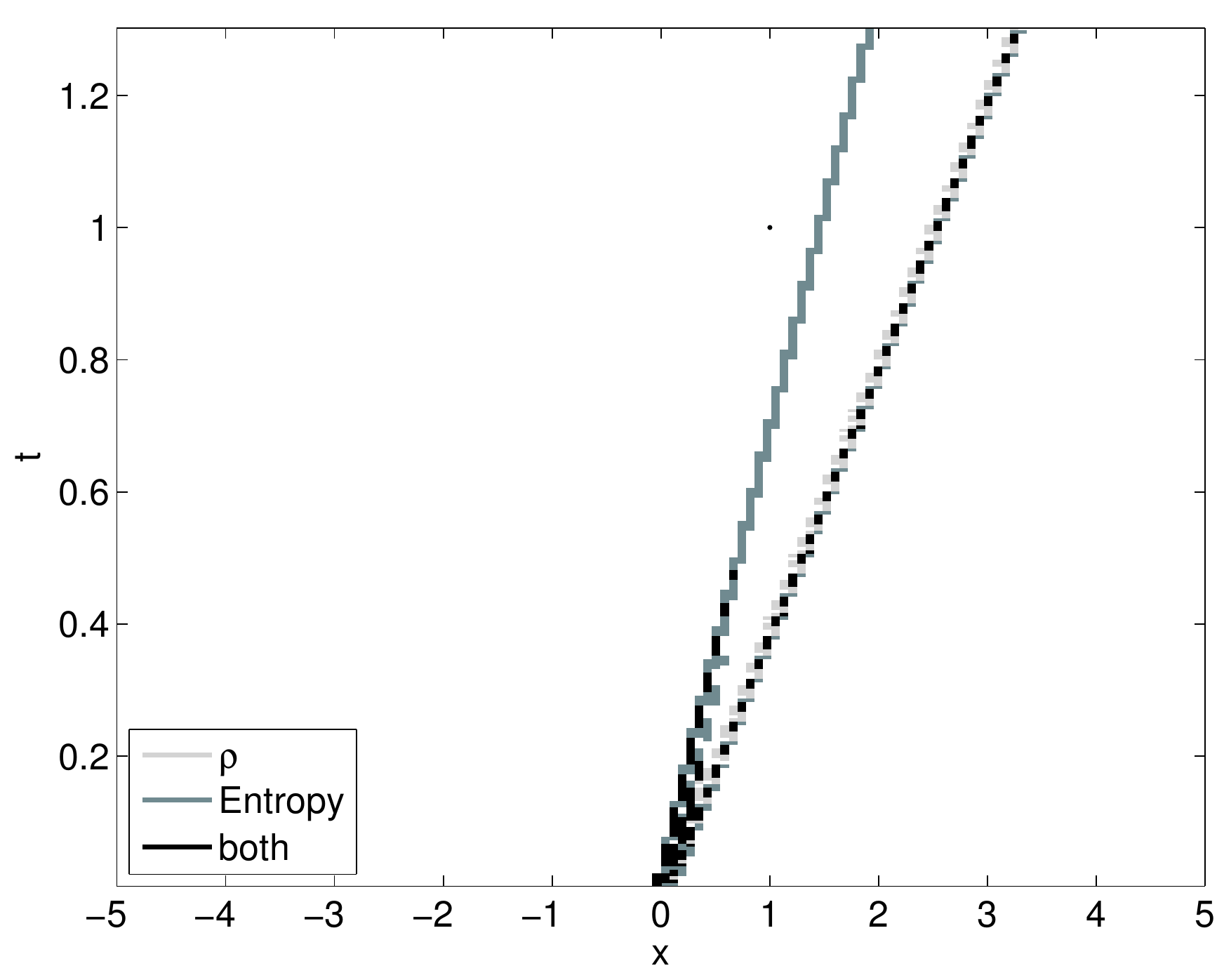}}
\subfigure[Harten, $k=1, \alpha = 1.5$]{\includegraphics[scale = 0.29]{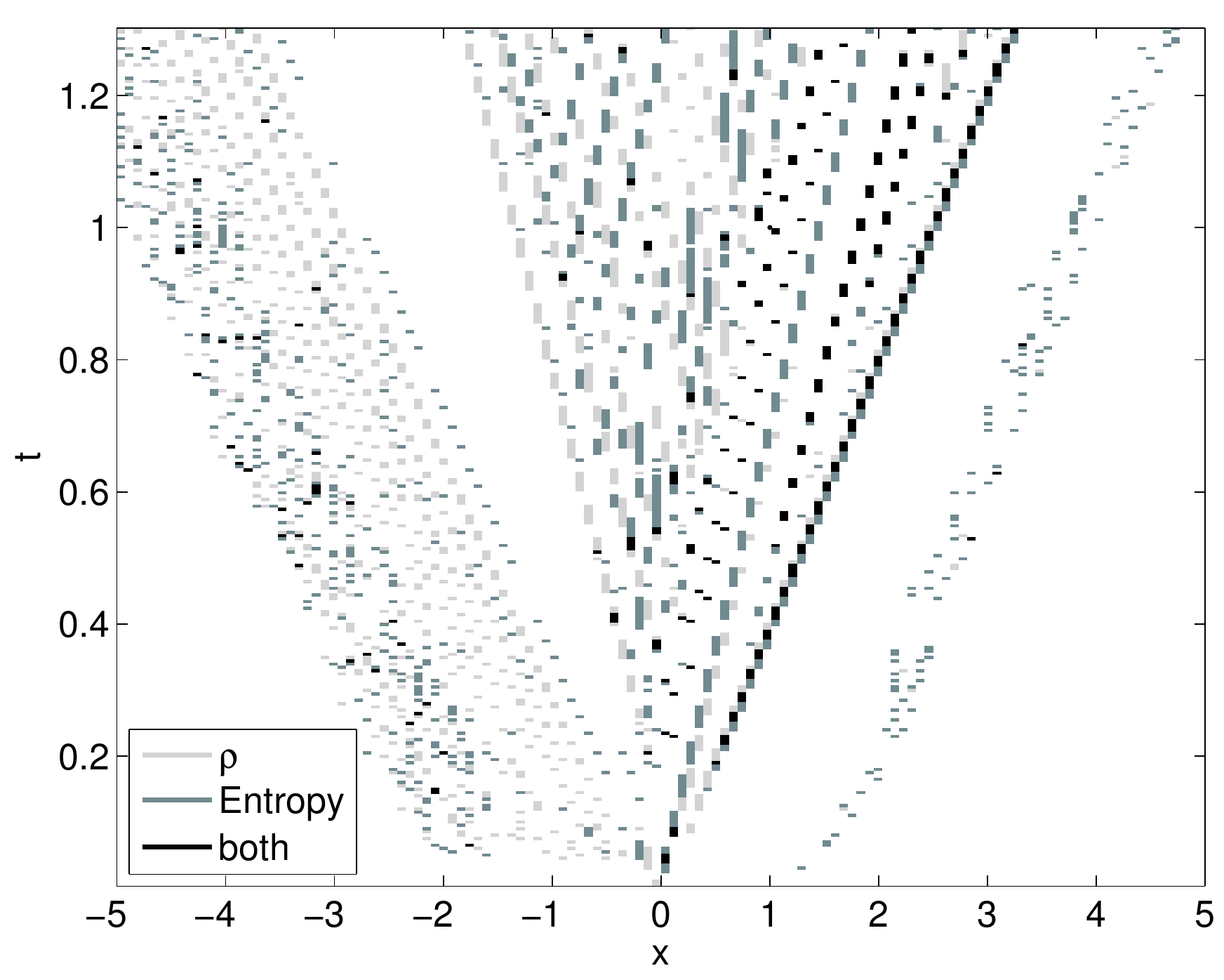}} \\
\vspace{-0.3cm}
\subfigure[KXRCF, $k=2$]{\includegraphics[scale = 0.29]{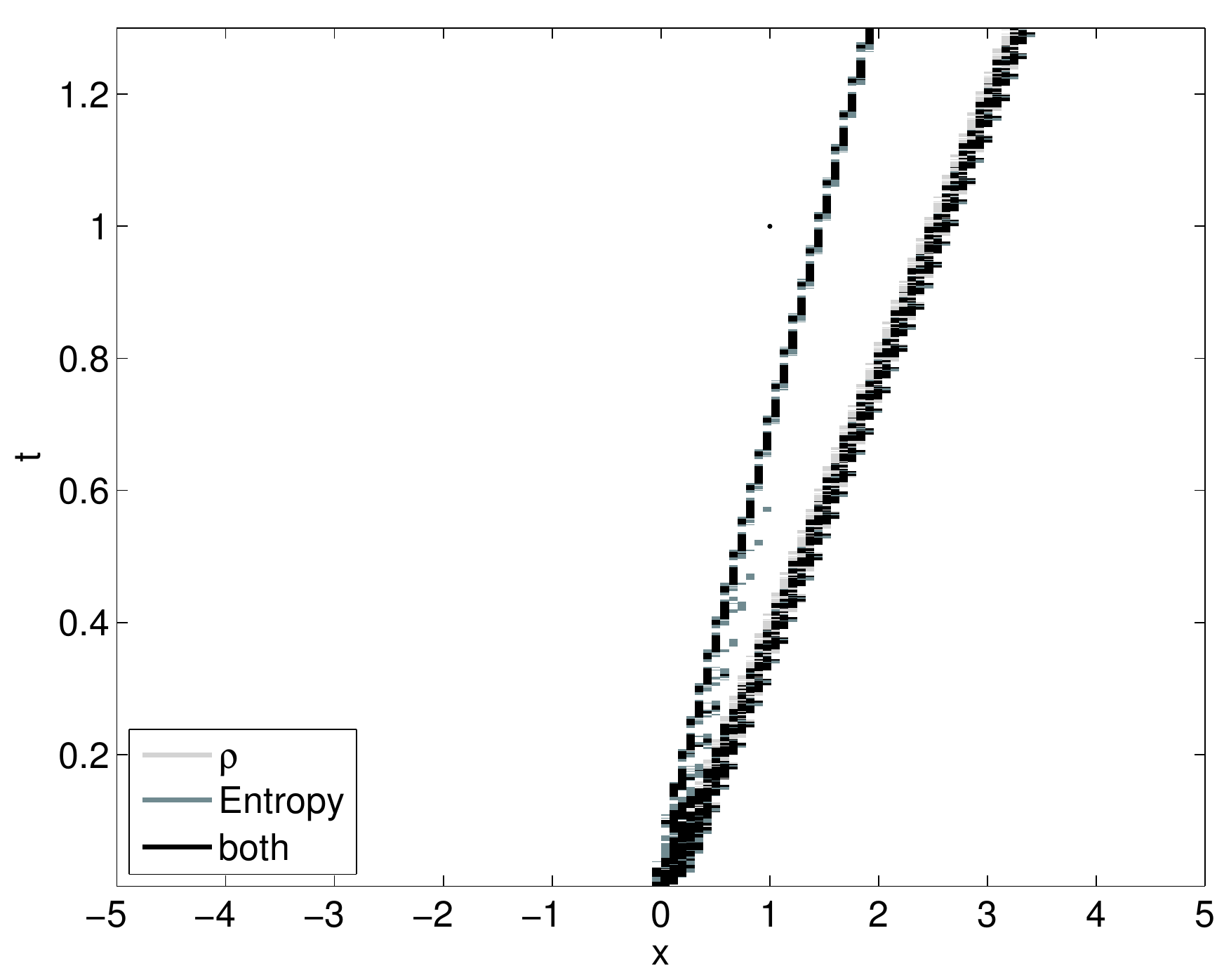}}
\subfigure[Harten, $k=2, \alpha = 1.5$]{\includegraphics[scale = 0.29]{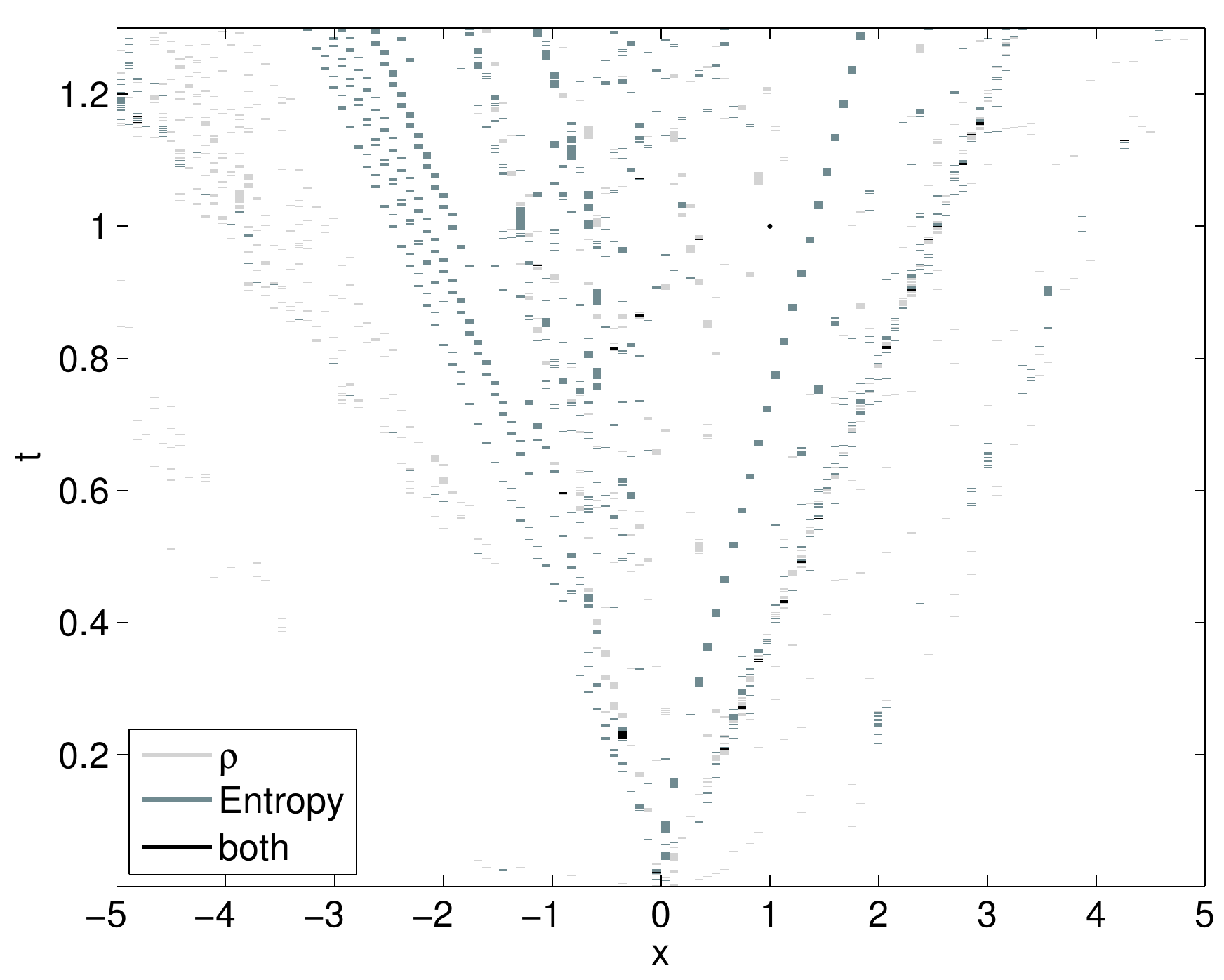}} \\
\vspace{-0.3cm}
\caption{The KXRCF or Harten's troubled-cell indicator on density and entropy, Lax, 128 elements.}\label{fig:LaxKH}
\end{figure}

\begin{figure}[ht!]
\centering
 \subfigure[$C=0.9$]{\includegraphics[scale = 0.22]{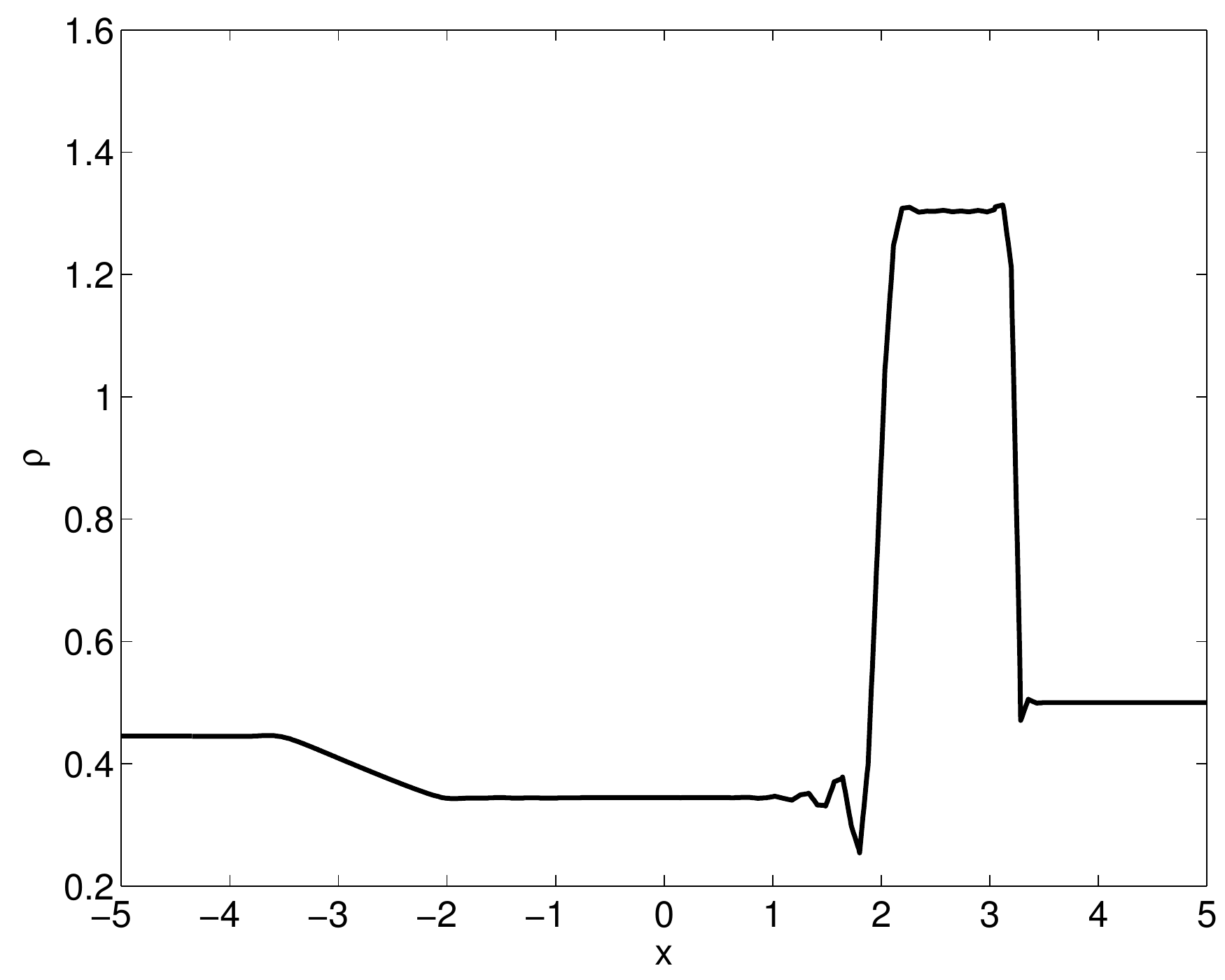}}
 \subfigure[$C=0.5$]{\includegraphics[scale = 0.22]{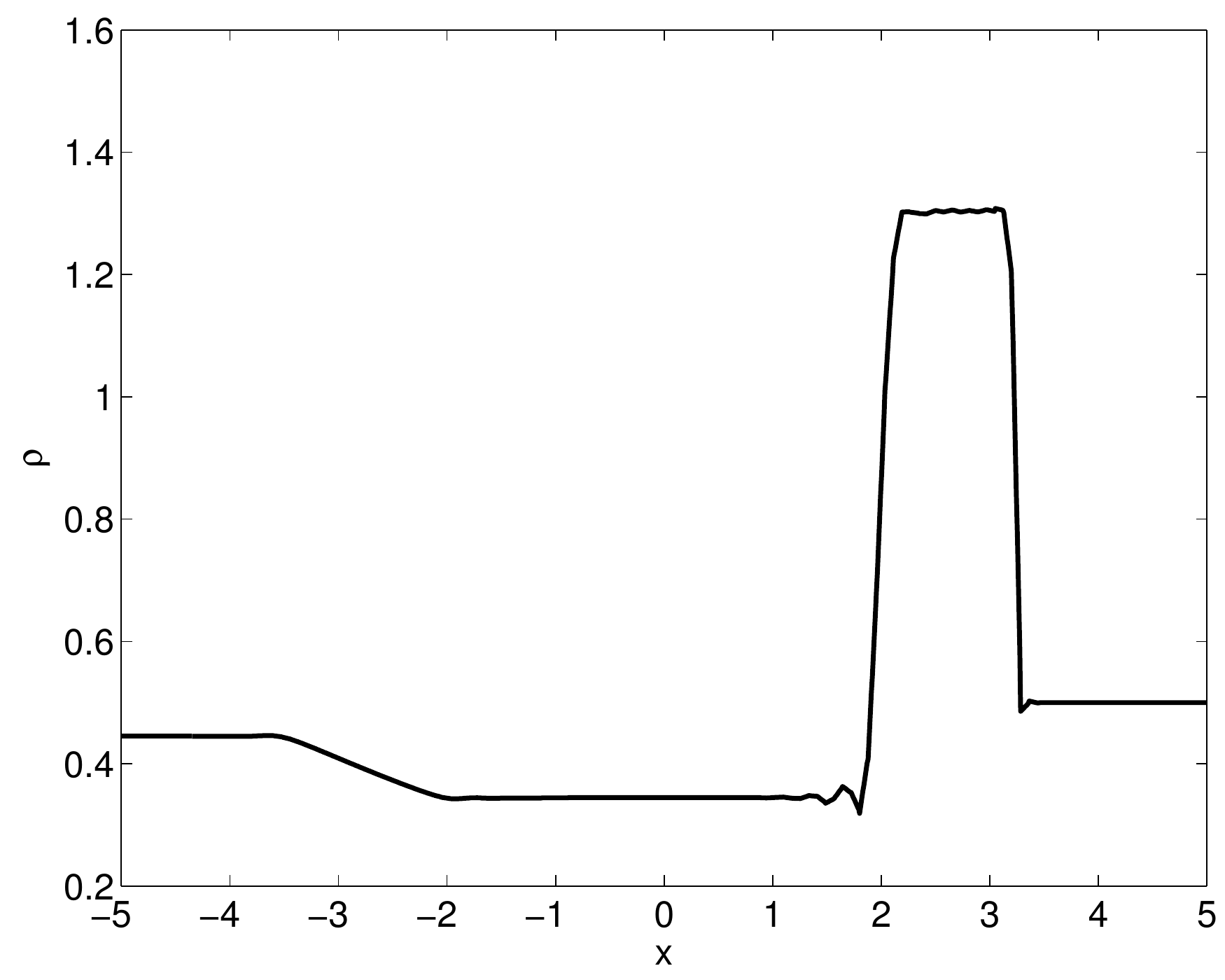}}
 \subfigure[$C=0.1$]{\includegraphics[scale = 0.22]{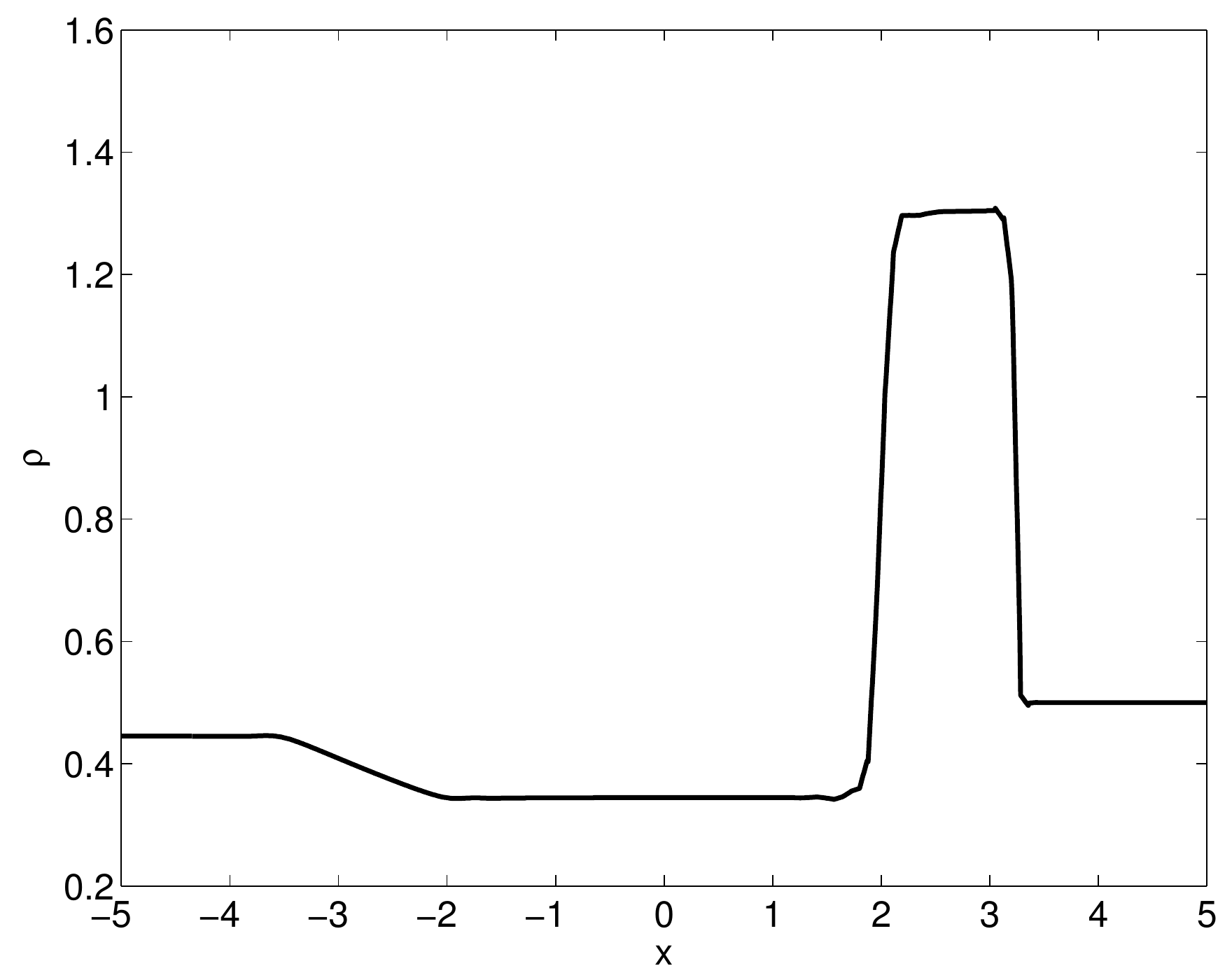}} \\
\vspace{-0.3cm}
 \subfigure[$C=0.9$]{\includegraphics[scale = 0.22]{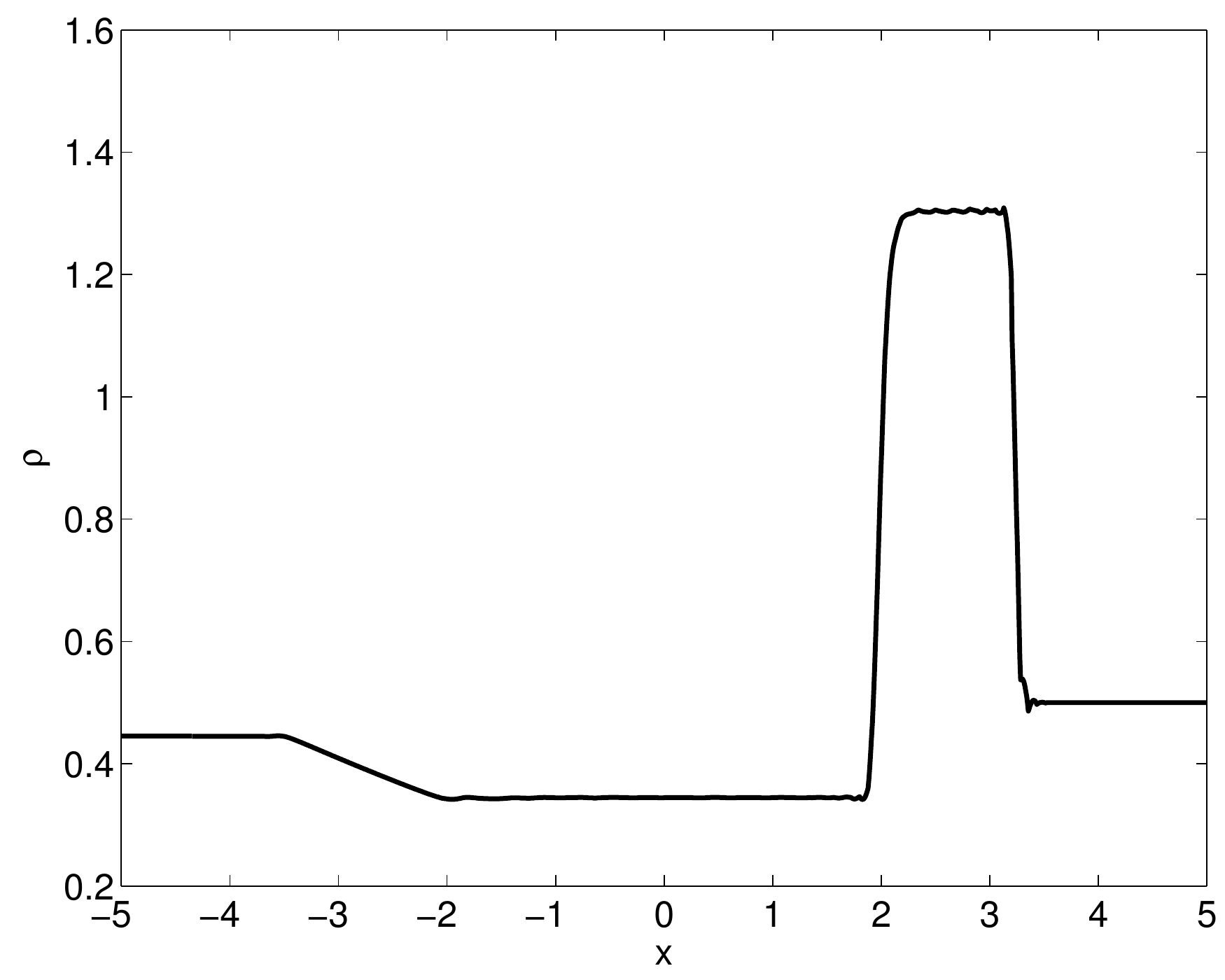}} 
 \subfigure[$C=0.5$]{\includegraphics[scale = 0.22]{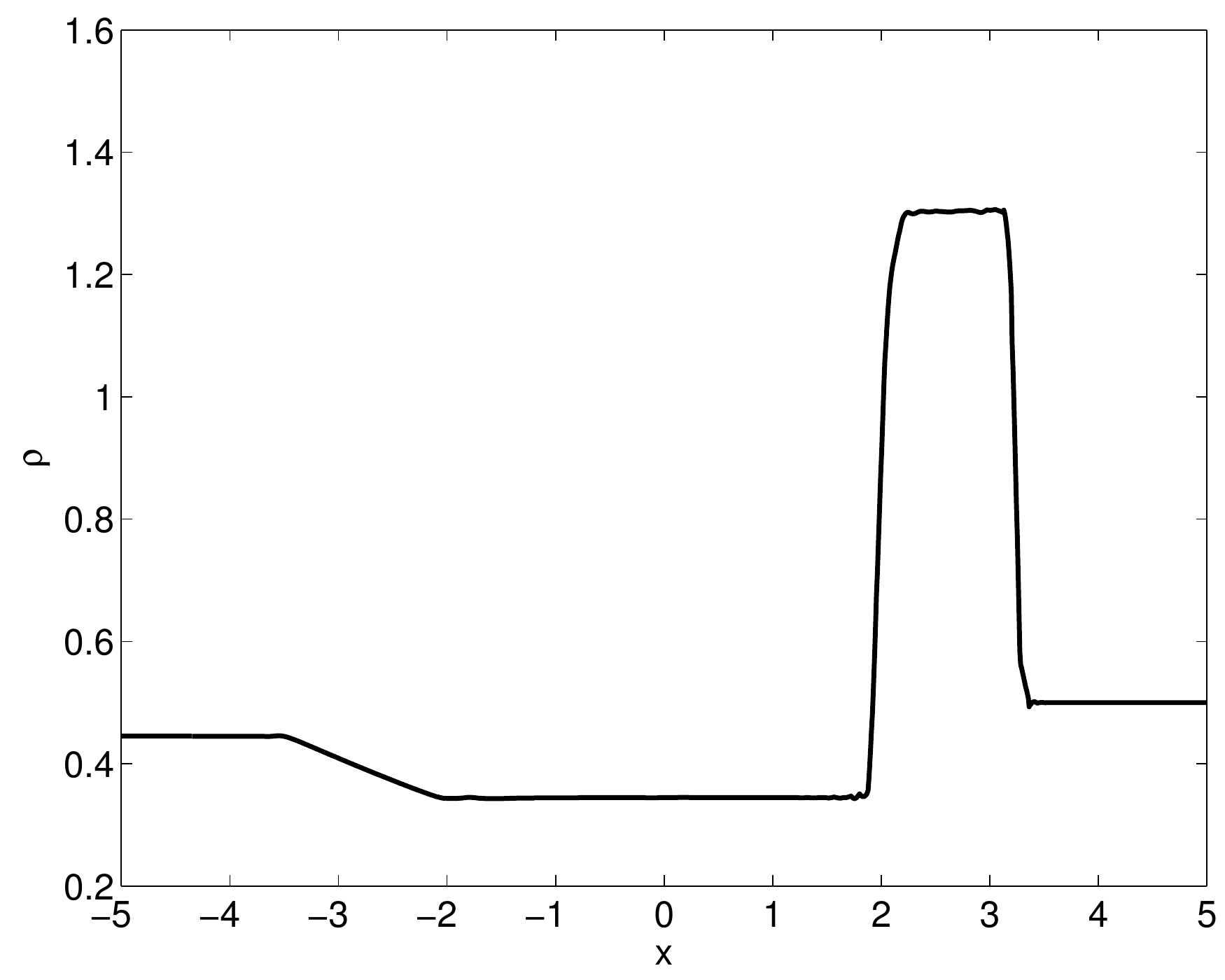}} 
 \subfigure[$C=0.1$]{\includegraphics[scale = 0.22]{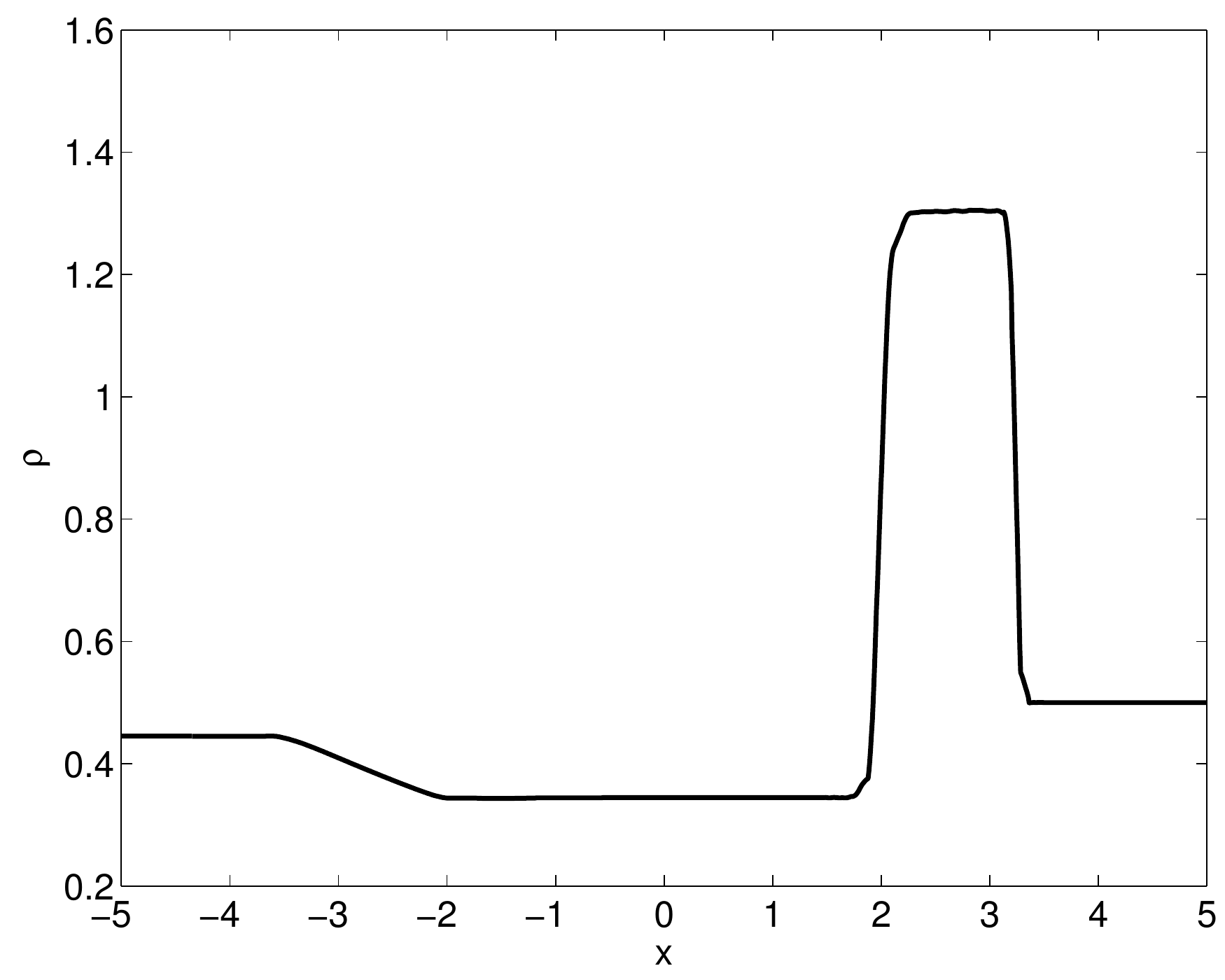}} \\
\vspace{-0.3cm}
\caption{Approximation at $T=1.3$, multiwavelet solution detector on density, Lax, 128 elements. First row: $k=1$, second row: $k=2$.} \label{fig:LaxCsol}
\end{figure}

\begin{figure}[ht!]
 \centering
\subfigure[KXRCF, $k=1$]{\includegraphics[scale = 0.29]{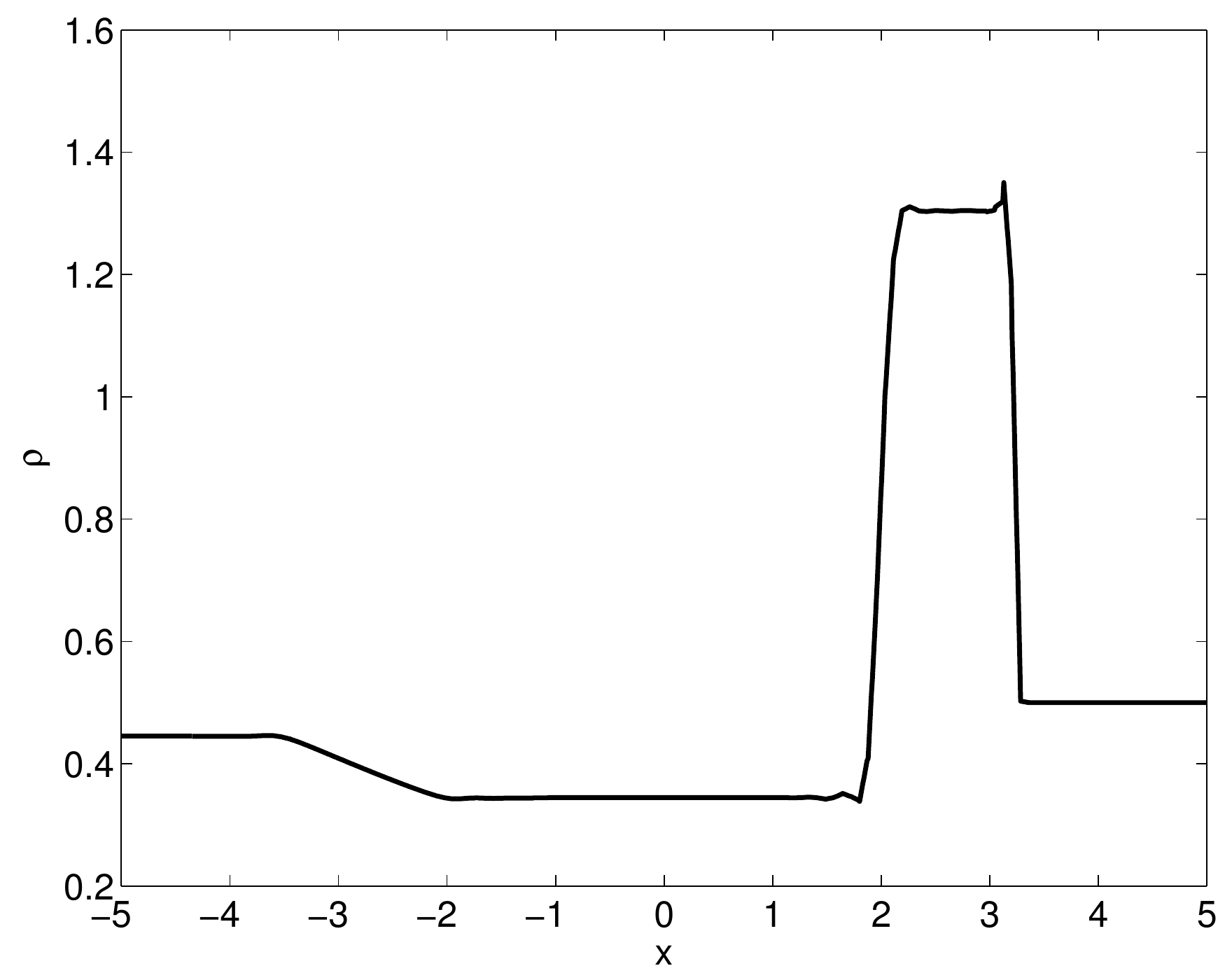}}
\subfigure[Harten, $k=1, \alpha = 1.5$]{\includegraphics[scale = 0.29]{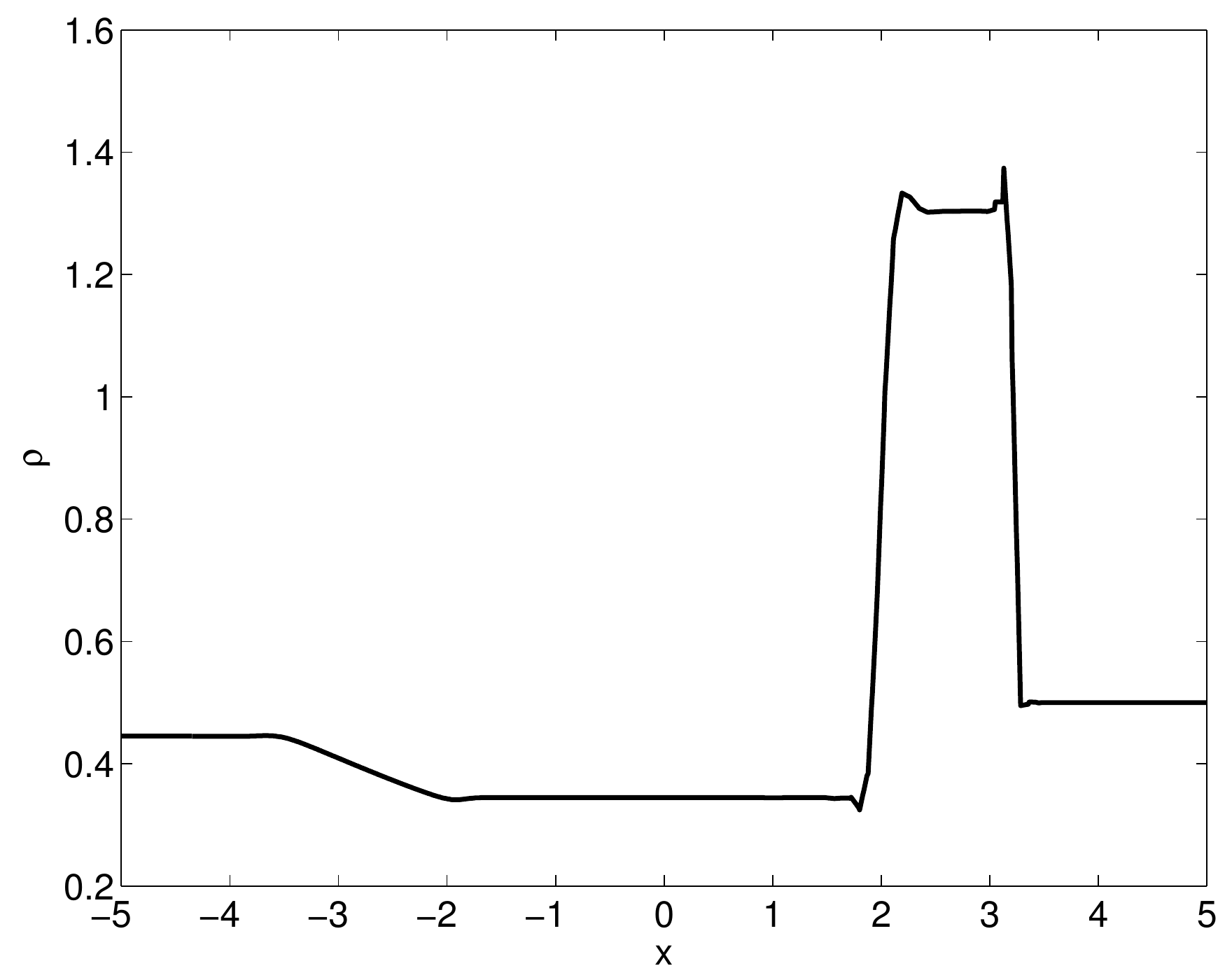}} \\
\vspace{-0.2cm}
\subfigure[KXRCF, $k=2$]{\includegraphics[scale = 0.29]{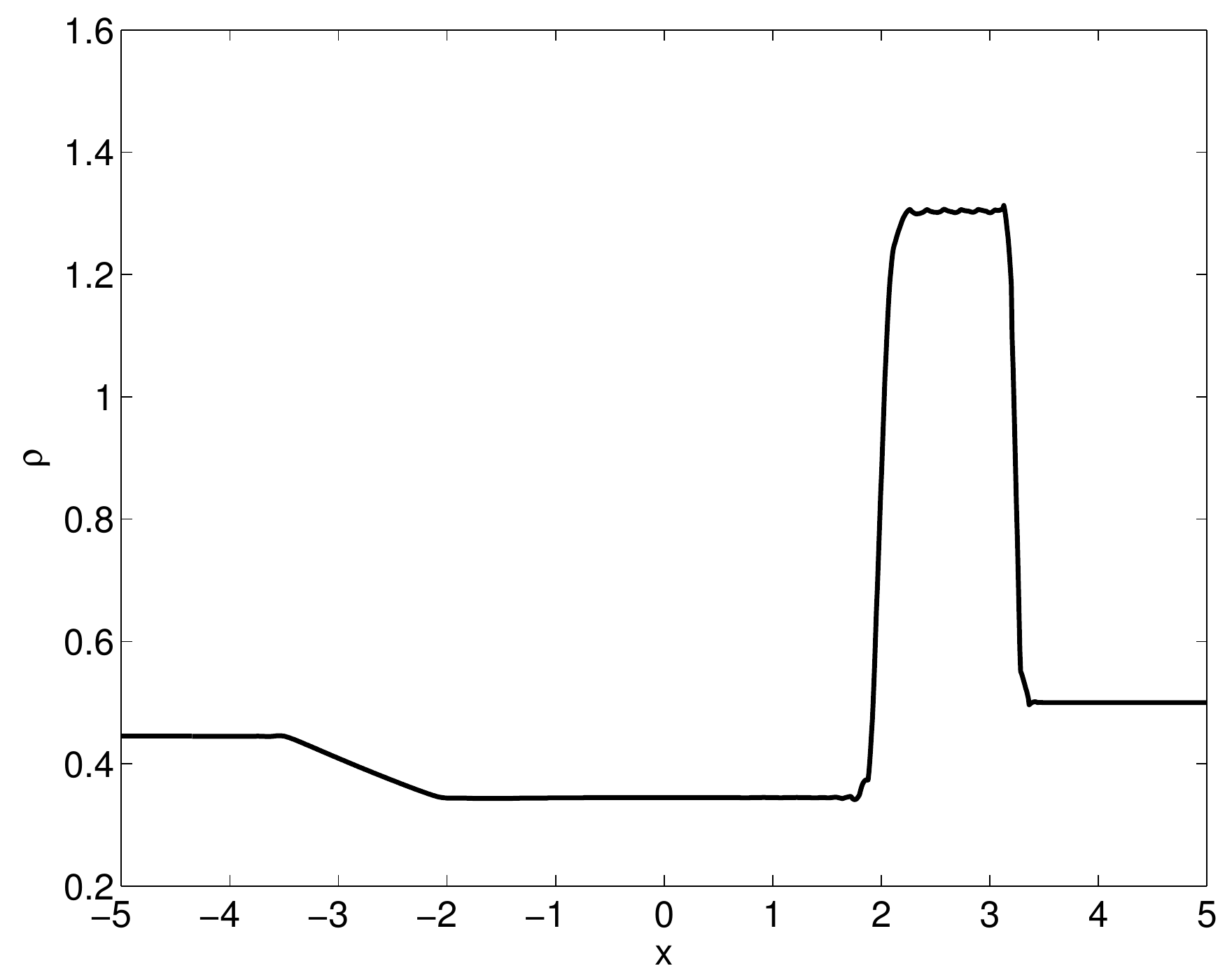}}
\subfigure[Harten, $k=2, \alpha = 1.5$]{\includegraphics[scale = 0.29]{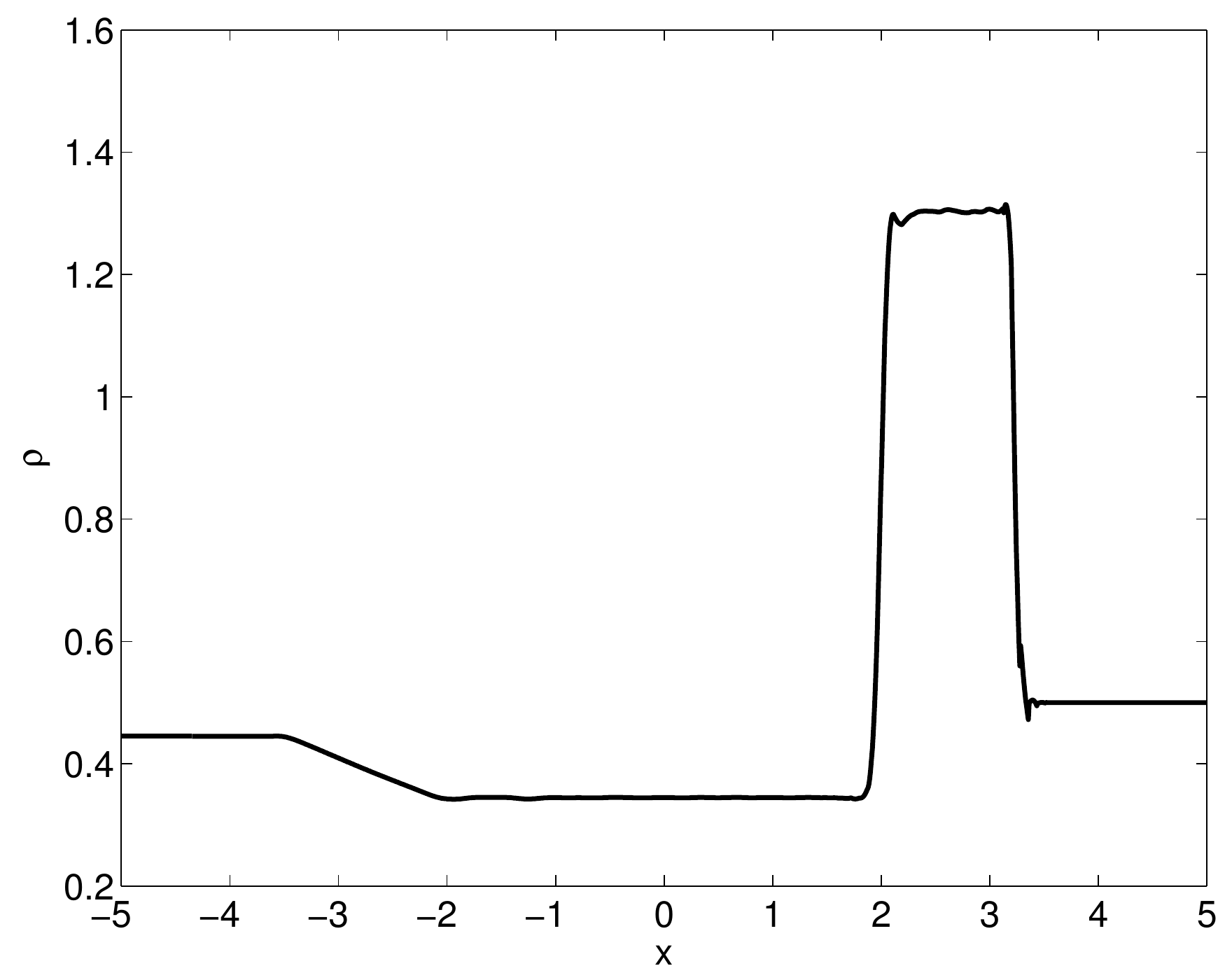}} \\
\vspace{-0.2cm}
\caption{Approximation at $T=1.3$, KXRCF or Harten's troubled-cell indicator on density and entropy, Lax, 128 elements.}\label{fig:LaxKHsol}
\end{figure}

The results using the multiwavelet indicator on density can be seen in Figures \ref{fig:LaxC} and \ref{fig:LaxCsol}, the KXRCF and Harten results (using density and entropy) are visualized in Figures \ref{fig:LaxKH} and \ref{fig:LaxKHsol} (final time $T=1.3$).   Note that the multiwavelet indicator does not detect the rarefaction wave for the given values of $C$ as this wave is more smooth than in Sod's shock tube. The value $C=0.9$ is too big, and oscillations are present in the solution. $C=0.1$ gives much better results.  The KXRCF indicator also does not detect the rarefaction wave. It is clearly visible that in the linear case, entropy detects the contact discontinuity in the solution. The approximations, however, remain oscillatory.  The detected elements using Harten's subcell resolution are scattered. Although Harten's method detects all regions with interesting features, it seems that this detector does not select enough neighboring elements to remove the oscillations. This can, however, be influenced by the choice of $\alpha$ as well as the choice of the limiter.

\newpage

\newpage

\subsubsection{Blast wave problem}
The third initial condition that is considered simulates the interaction of two blast waves  \cite{Woo84C}.  This is given by
 
\[
  \rho (x,0)  = 1, \quad  u(x,0) = 0, \quad   p(x,0) = \left\{ \begin{array}{cl} 1000, & 0 \leq x < 0.1, \\ 0.01, & 0.1 \leq x < 0.9; \\ 100, & 0.9 \leq x \leq 1. \end{array} \right.
\]

Here, the boundary conditions of Shu et al. \cite{Shu89O} are used and the detected troubled cells in time are compared using different troubled-cell indicators.  The combination of Harten's troubled-cell indicator and the moment limiter is unstable for this example.  This possibility was also noticed in \cite{Zhu09Q}. Therefore, the multiwavelet approach will be tested against the KXRCF indicator only.

\begin{figure}[ht!]
\centering
 \subfigure[$C=0.25$]{\includegraphics[scale = 0.21]{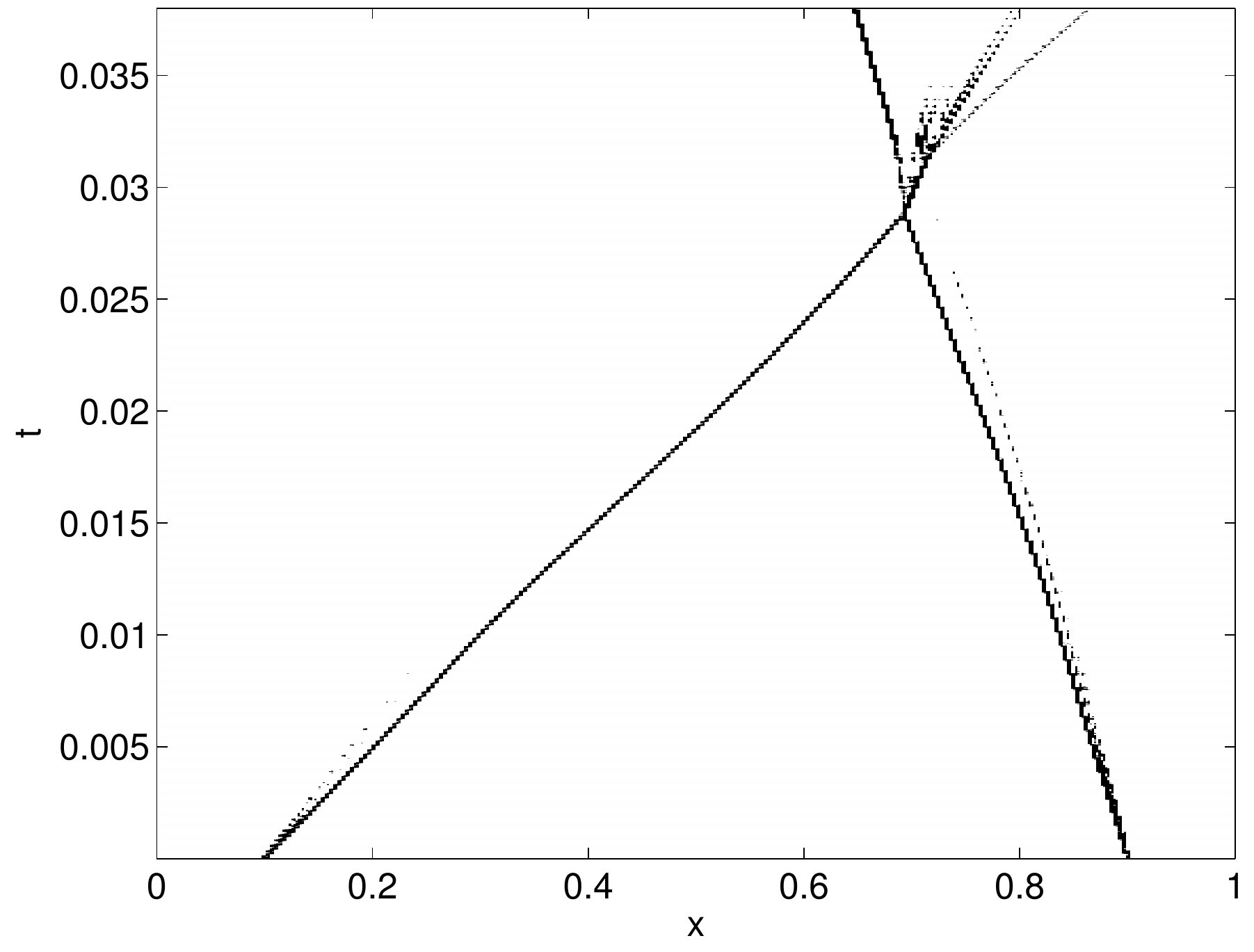}}
 \subfigure[$C=0.1$]{\includegraphics[scale = 0.21]{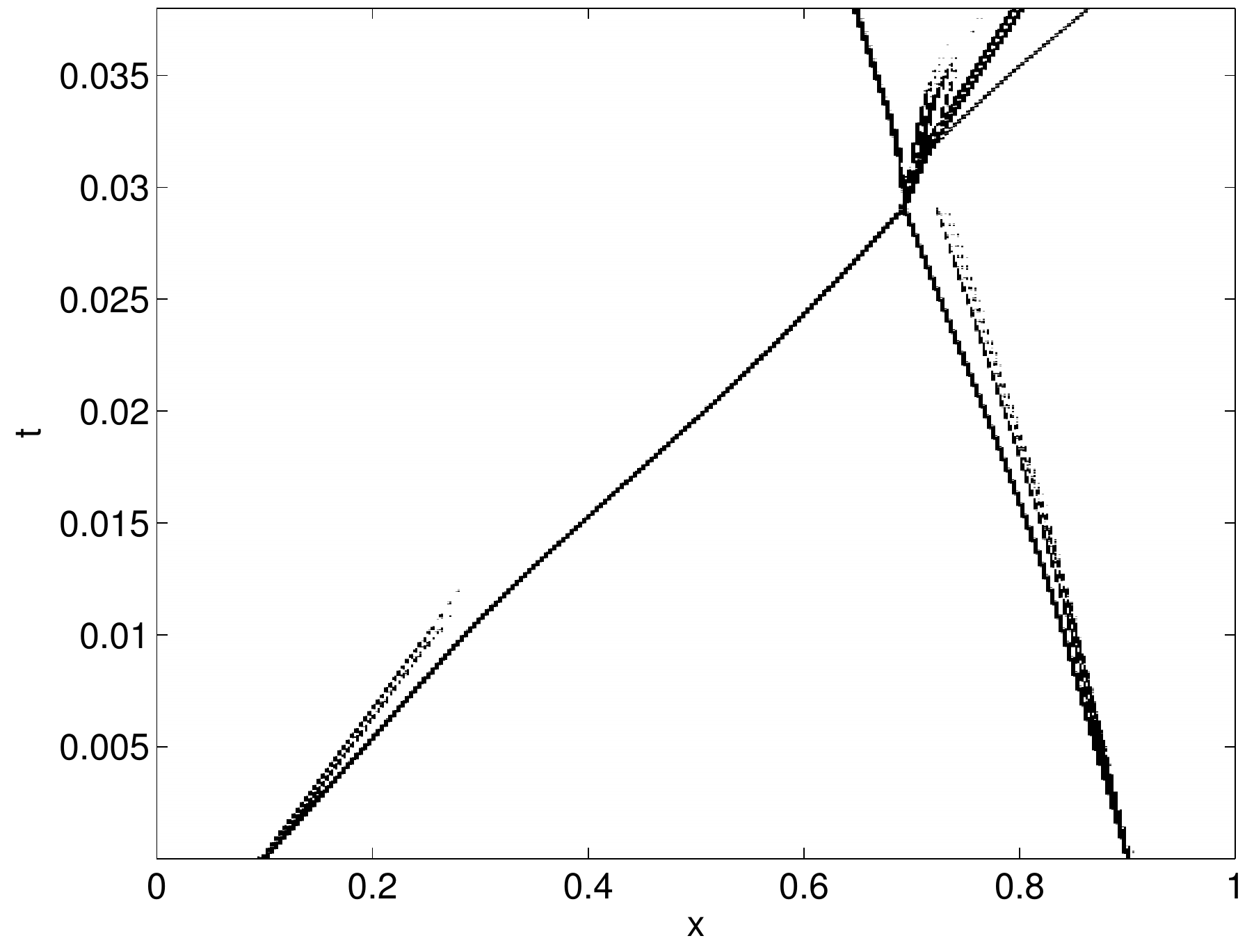}} 
 \subfigure[$C=0.05$]{\includegraphics[scale = 0.21]{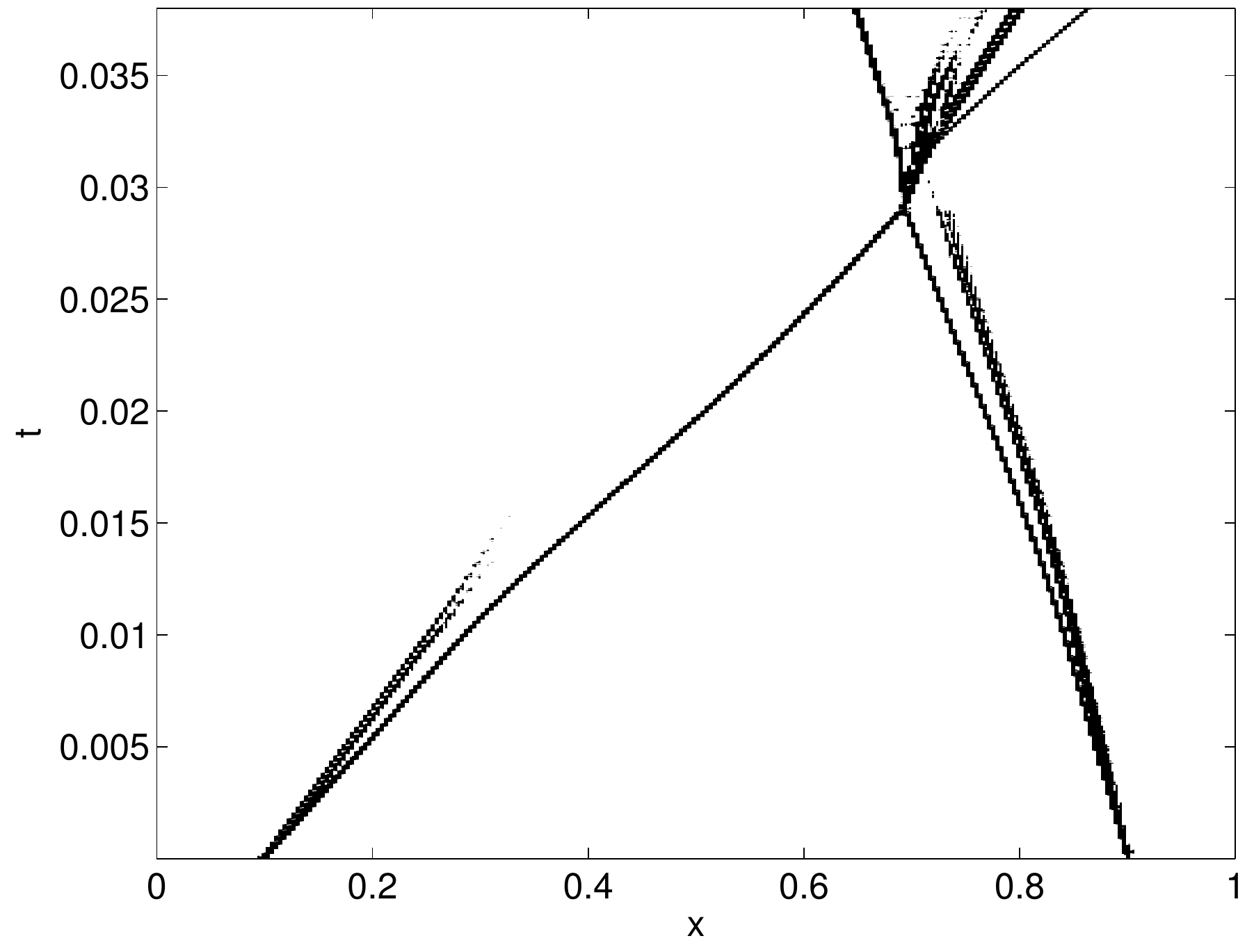}} \\
 \subfigure[$C=0.1$]{\includegraphics[scale = 0.21]{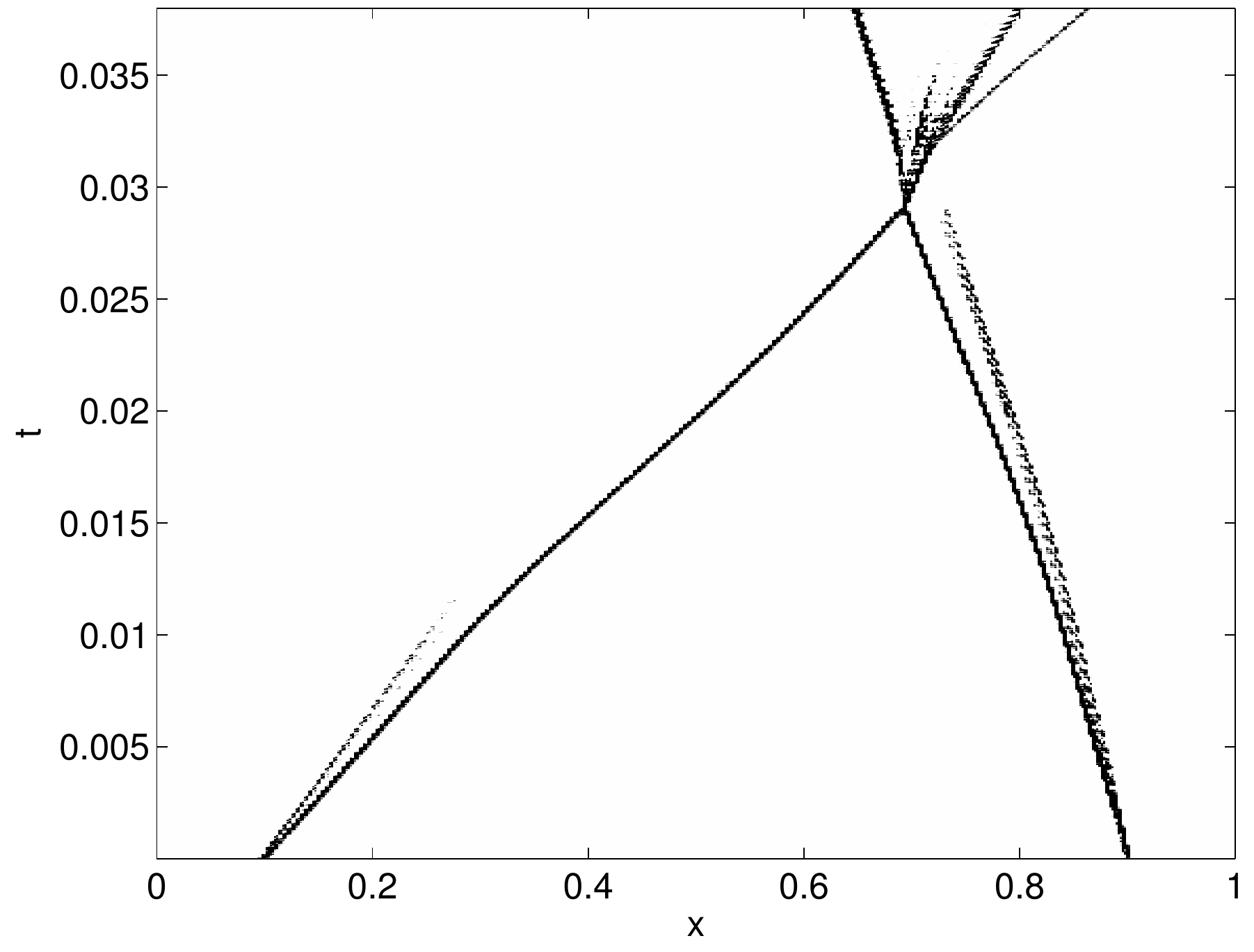}} 
 \subfigure[$C=0.05$]{\includegraphics[scale = 0.21]{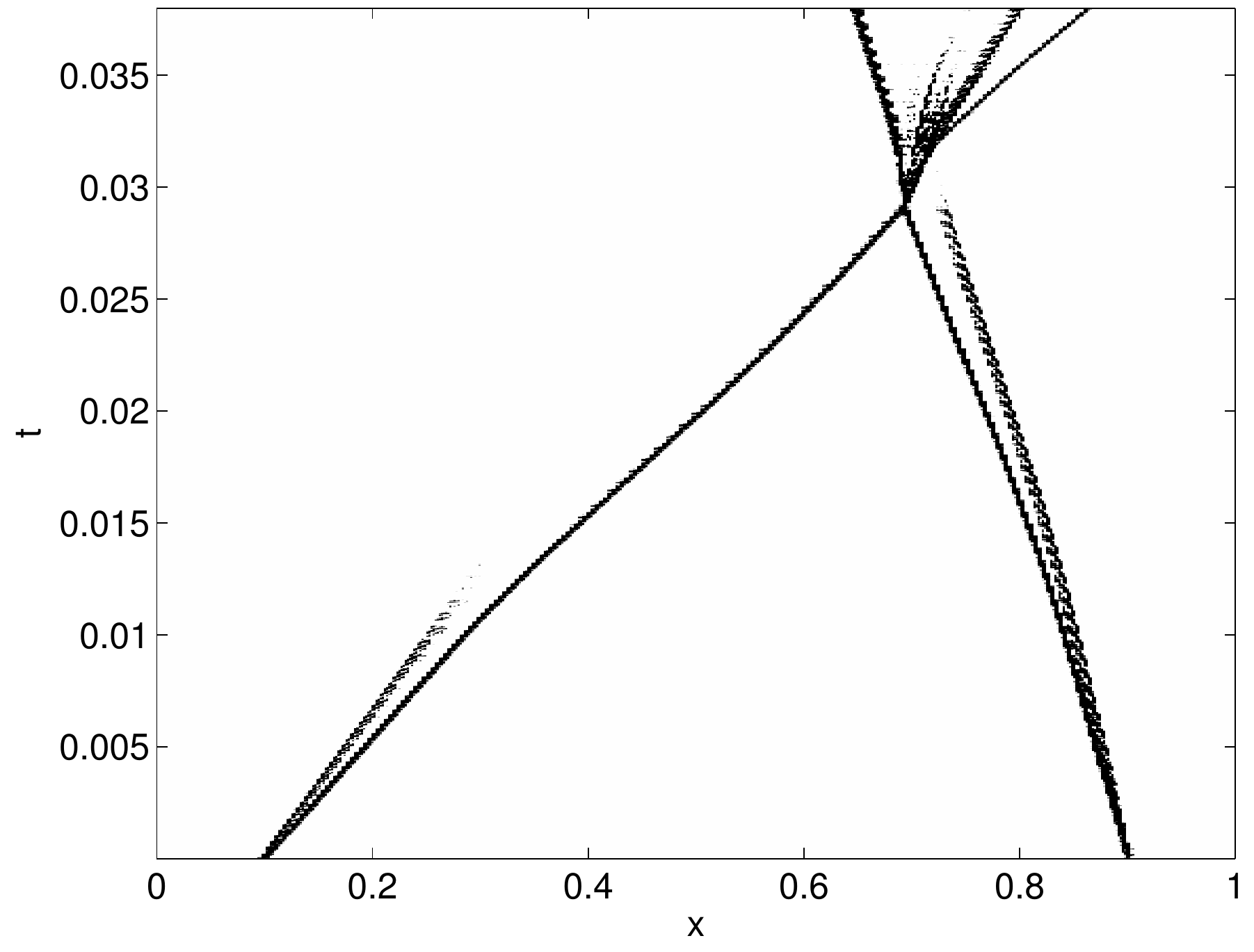}} 
 \subfigure[$C=0.01$]{\includegraphics[scale = 0.21]{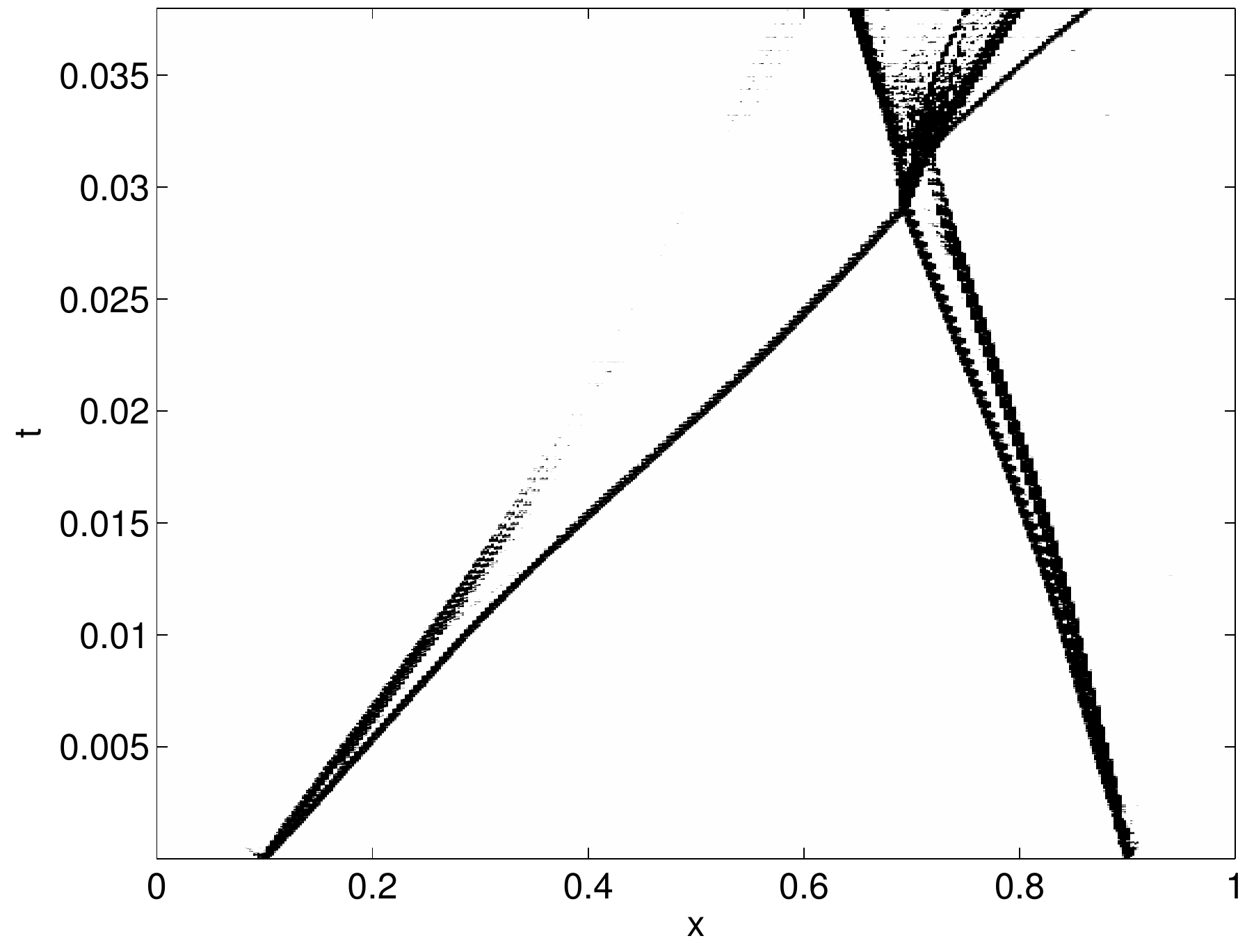}} \\
\caption{Time history plot of detected troubled cells, multiwavelet shock detector on density, Blast, 512 elements. First row: $k=1$, second row: $k=2$.}\label{fig:BlastCn9}
\end{figure}

\begin{figure}[ht!]
 \centering
\subfigure[KXRCF, $k=1$]{\includegraphics[scale = 0.27]{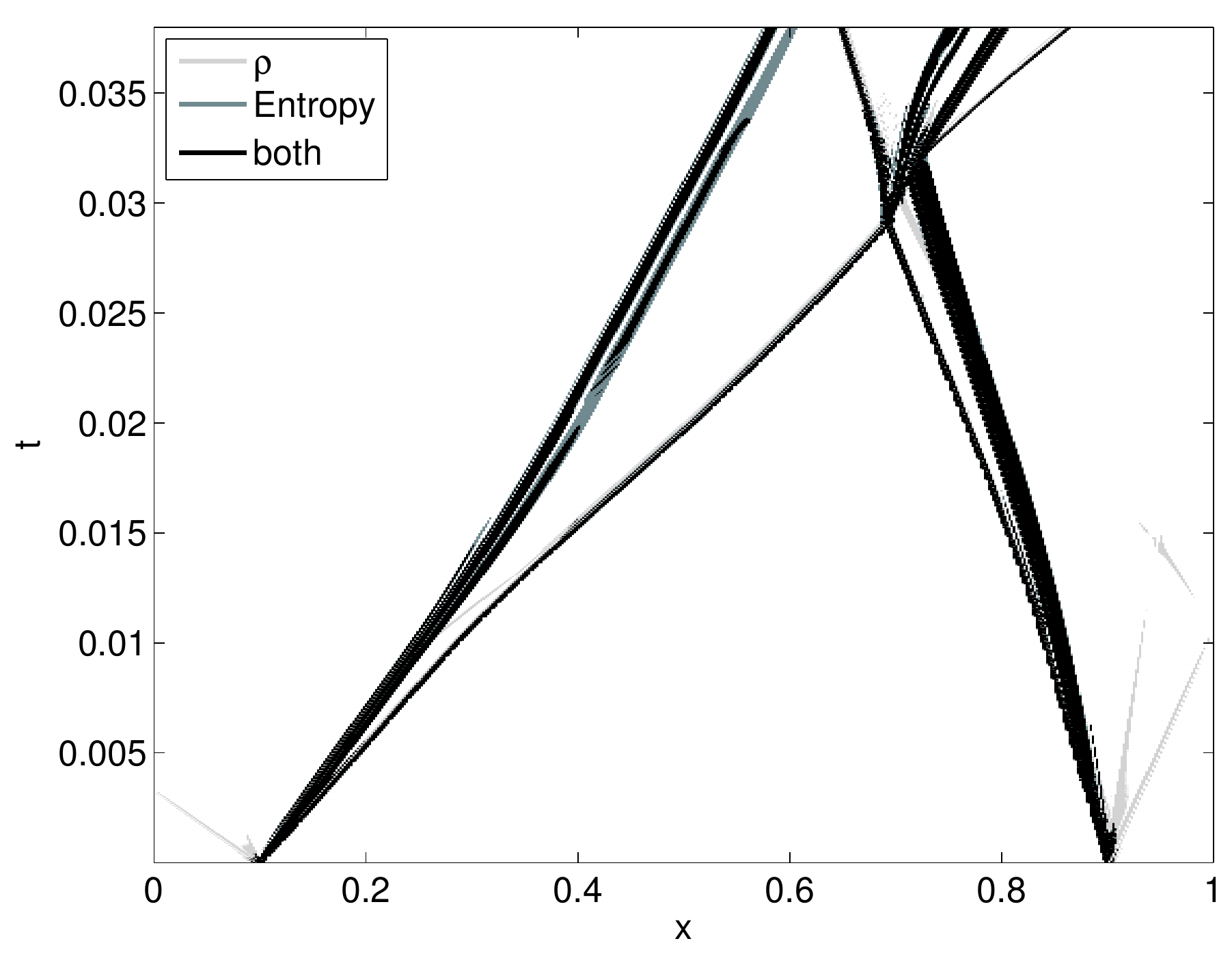}}
\subfigure[KXRCF, $k=2$]{\includegraphics[scale = 0.27]{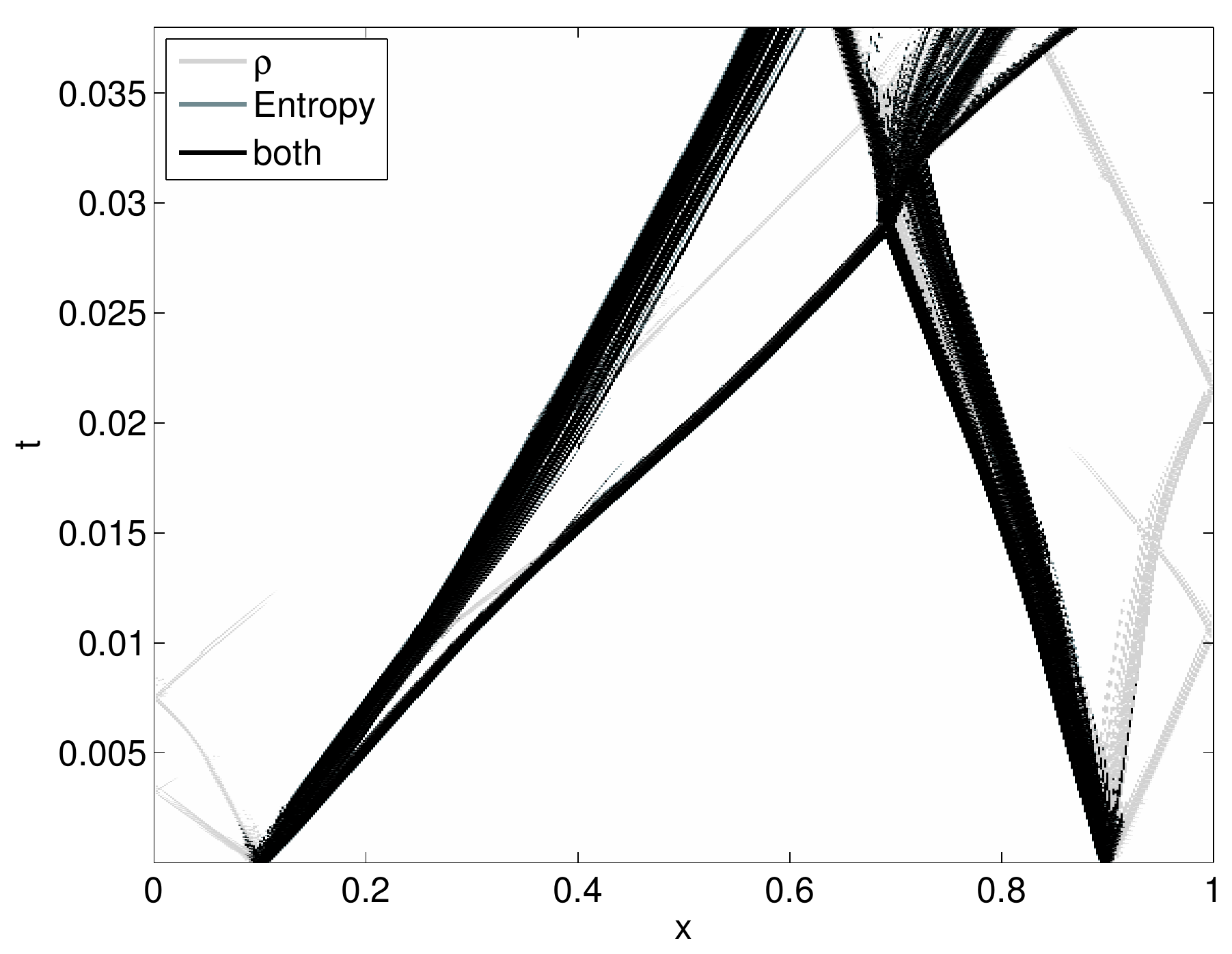}}
\caption{The KXRCF troubled-cell indicator (density and entropy), Blast, 512 elements.}\label{fig:BlastKH}
\end{figure}

\begin{figure}[ht!]
\centering
 \subfigure[$C=0.25$]{\includegraphics[scale = 0.22]{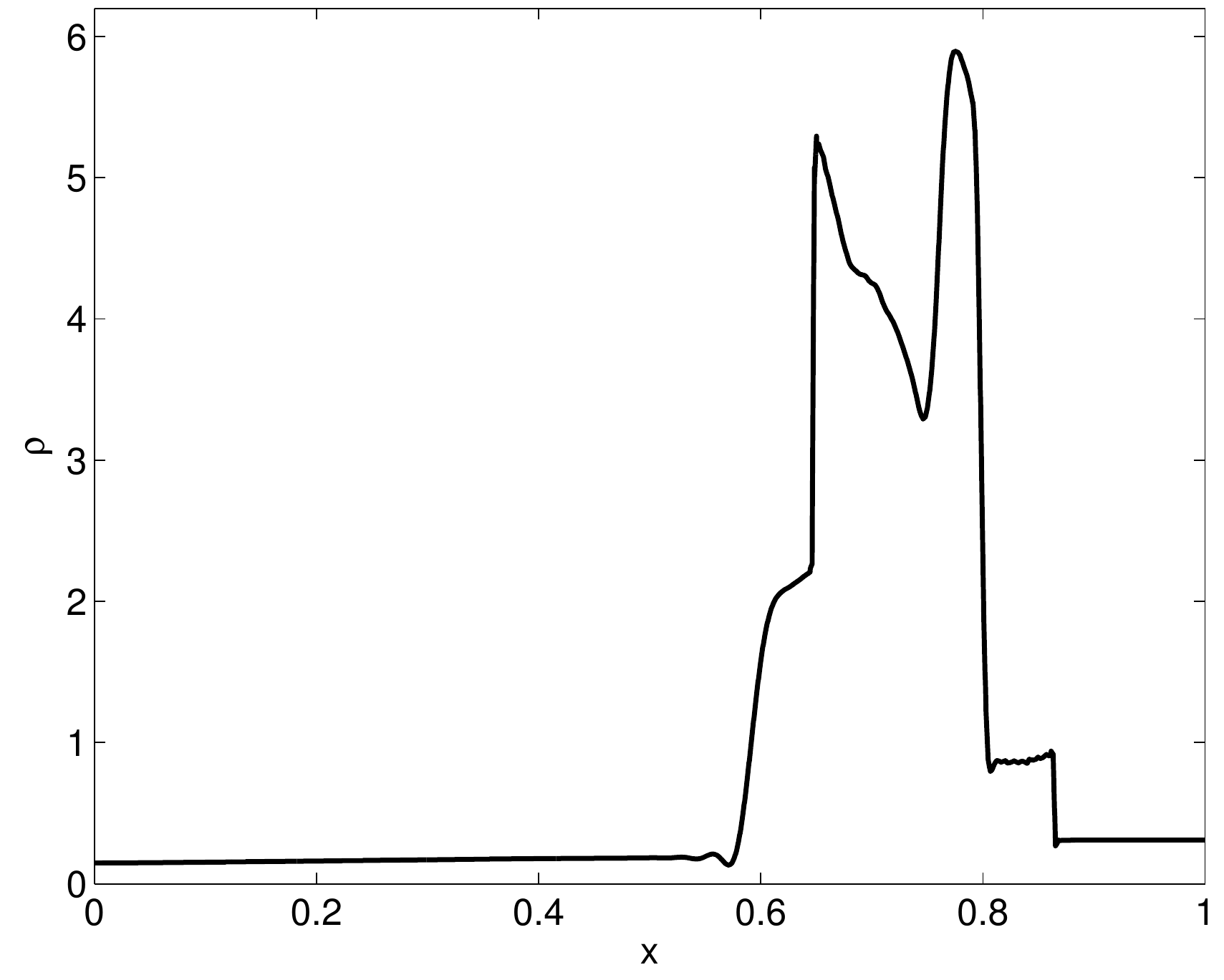}}
 \subfigure[$C=0.1$]{\includegraphics[scale = 0.22]{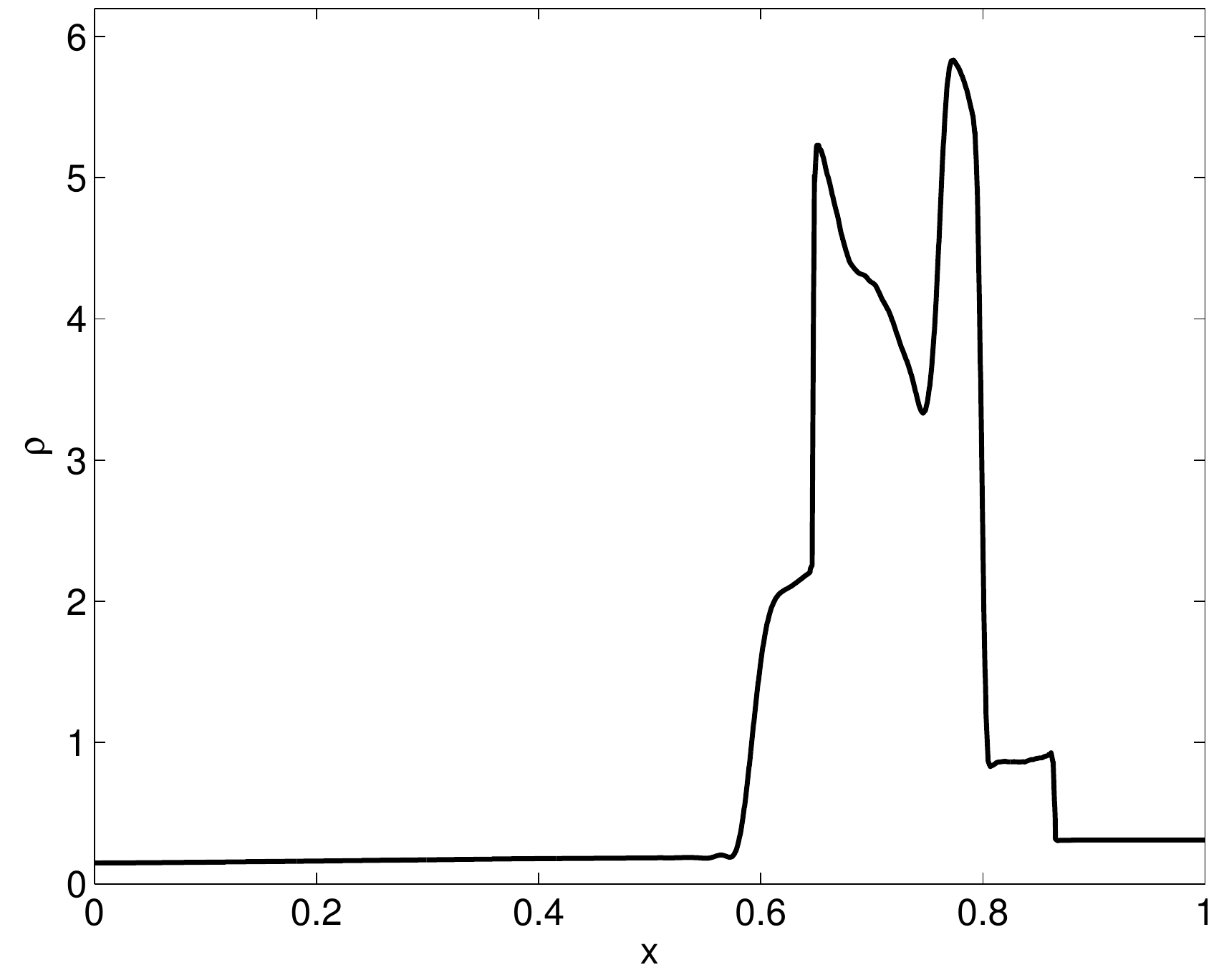}} 
 \subfigure[$C=0.05$]{\includegraphics[scale = 0.22]{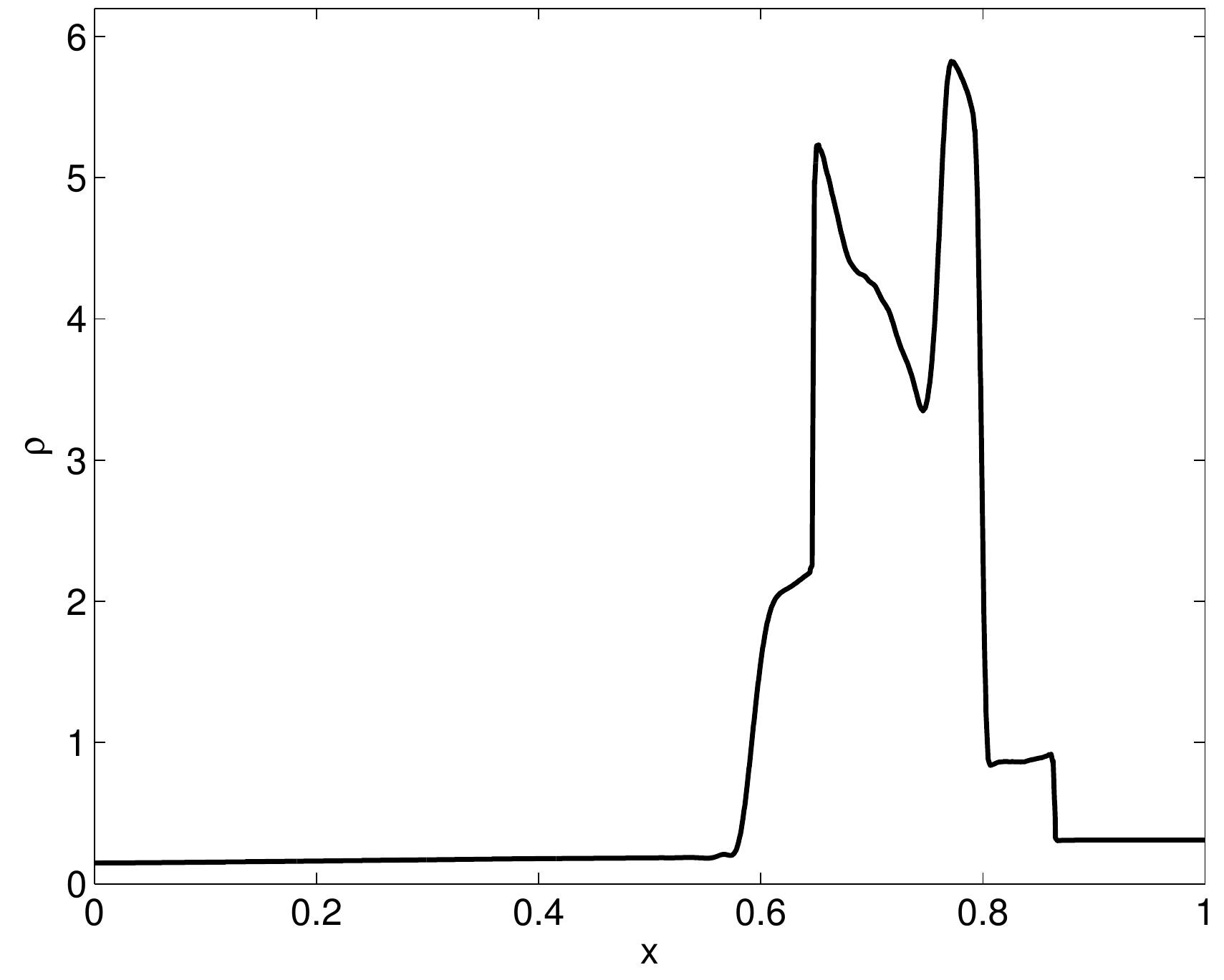}} \\
 \subfigure[$C=0.1$]{\includegraphics[scale = 0.22]{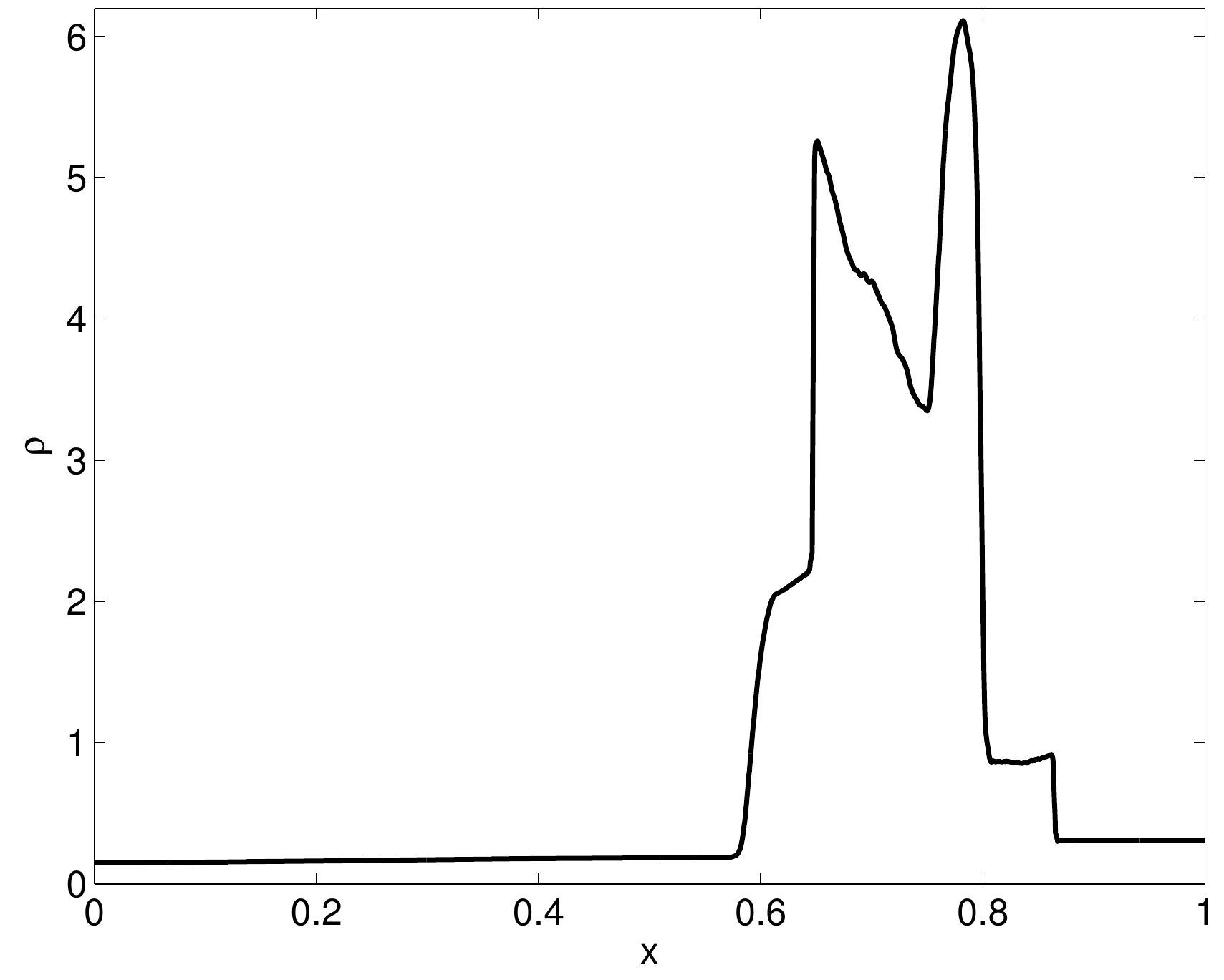}} 
 \subfigure[$C=0.05$]{\includegraphics[scale = 0.22]{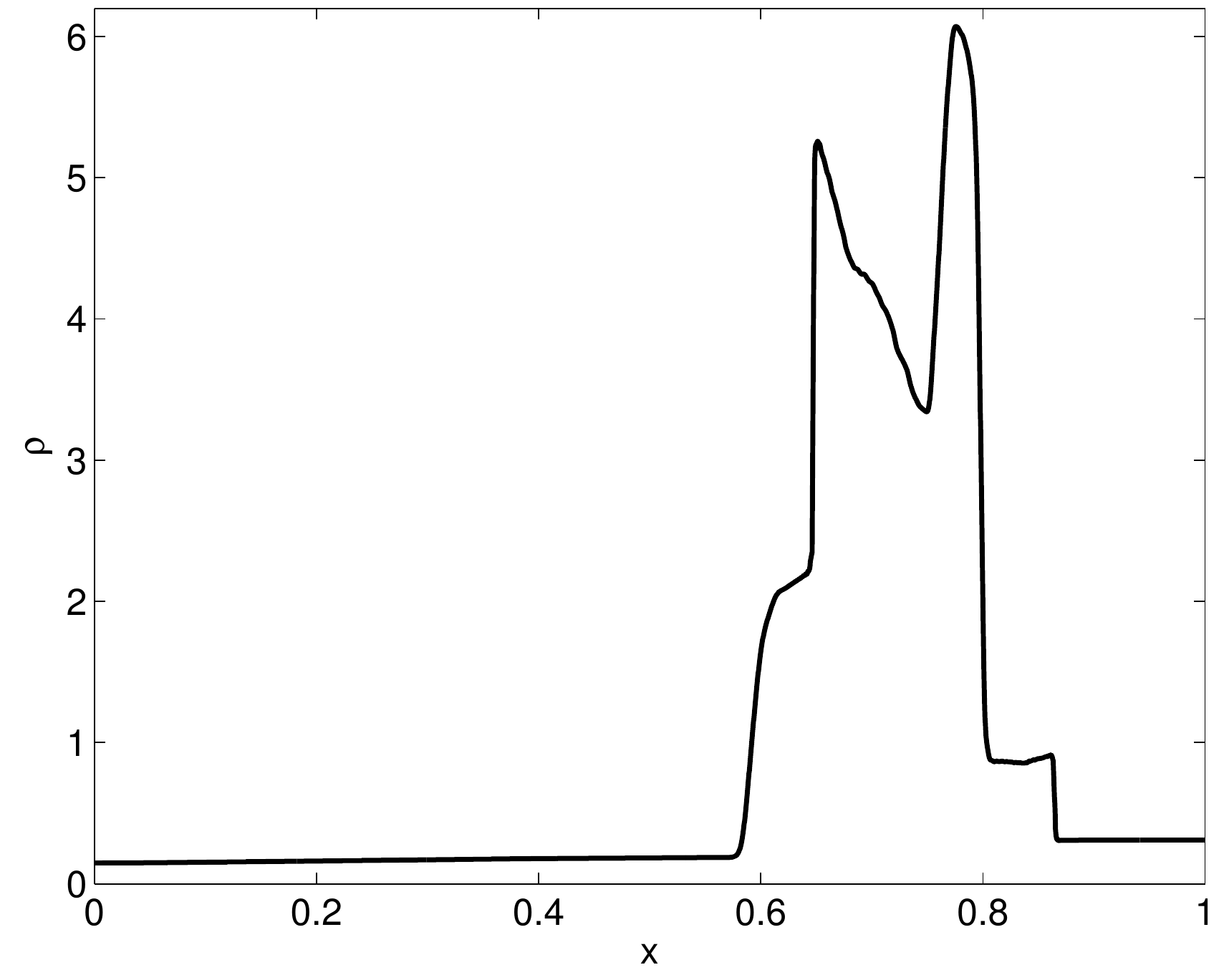}} 
 \subfigure[$C=0.01$]{\includegraphics[scale = 0.22]{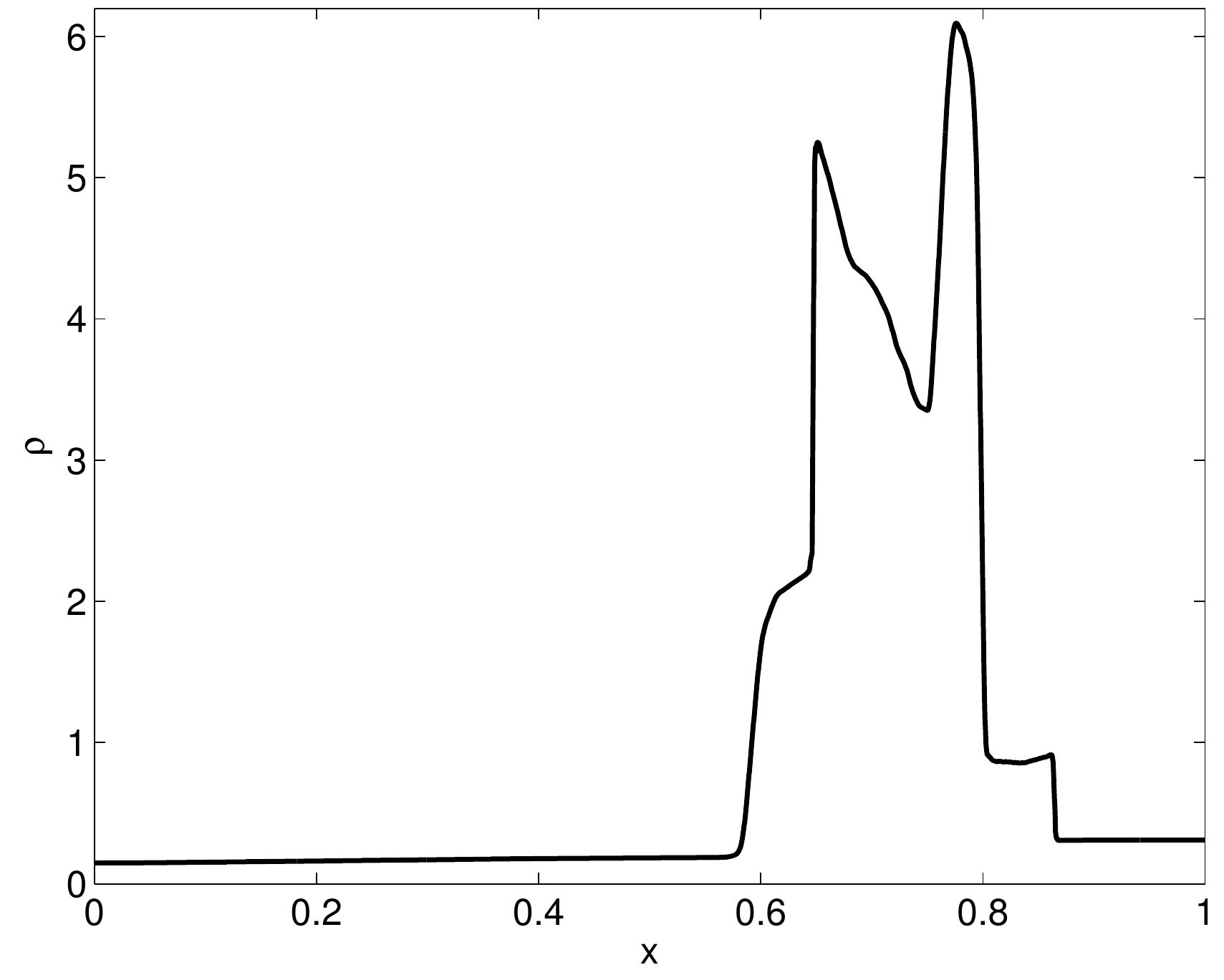}} \\
\caption{Approximation at $T=0.038$, multiwavelet detector on density, Blast, 512 elements. First row: $k=1$, second row: $k=2$.}\label{fig:BlastCn9sol}
\end{figure}

\begin{figure}[ht!]
 \centering
\subfigure[KXRCF, $k=1$]{\includegraphics[scale = 0.27]{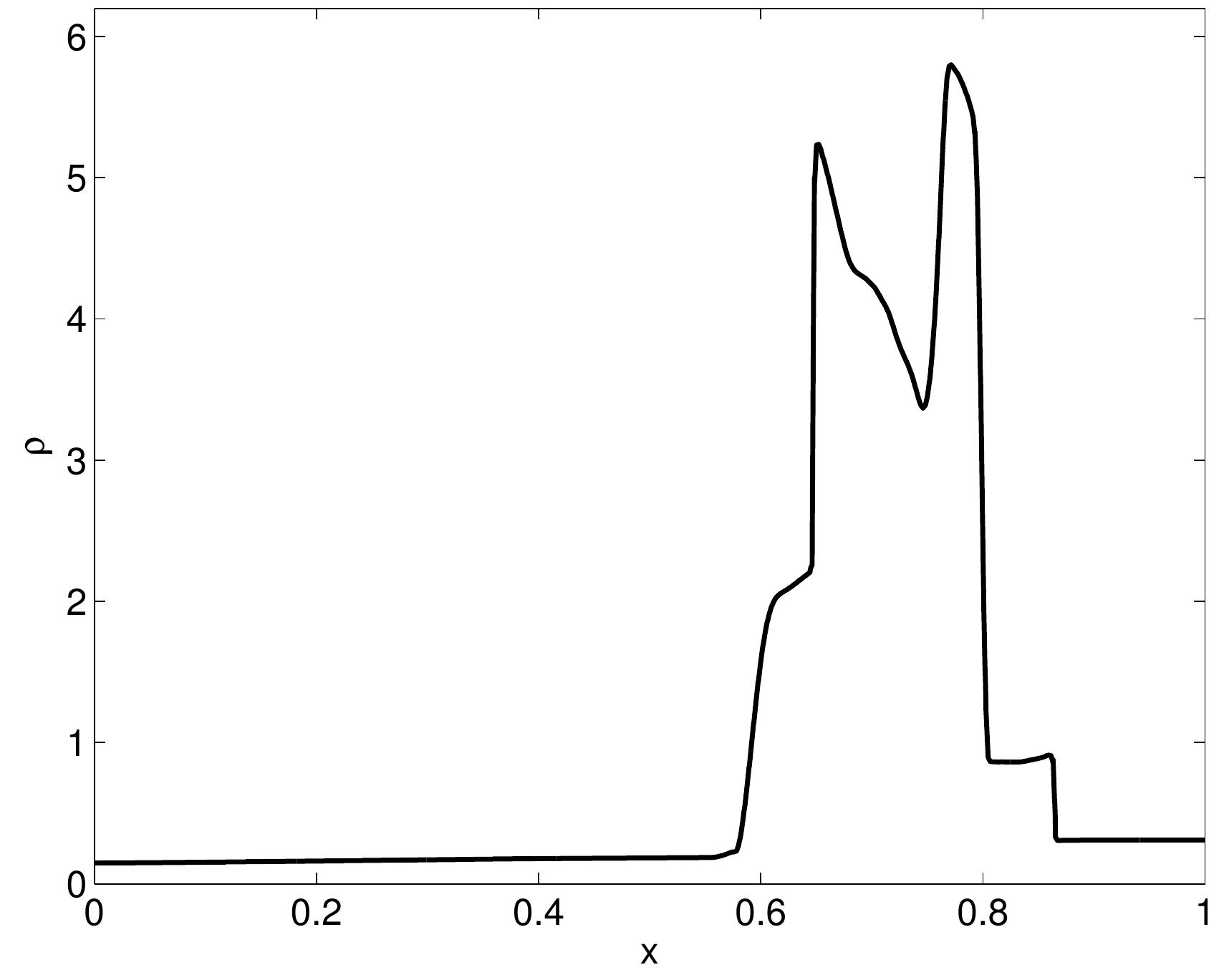}}
\subfigure[KXRCF, $k=2$]{\includegraphics[scale = 0.27]{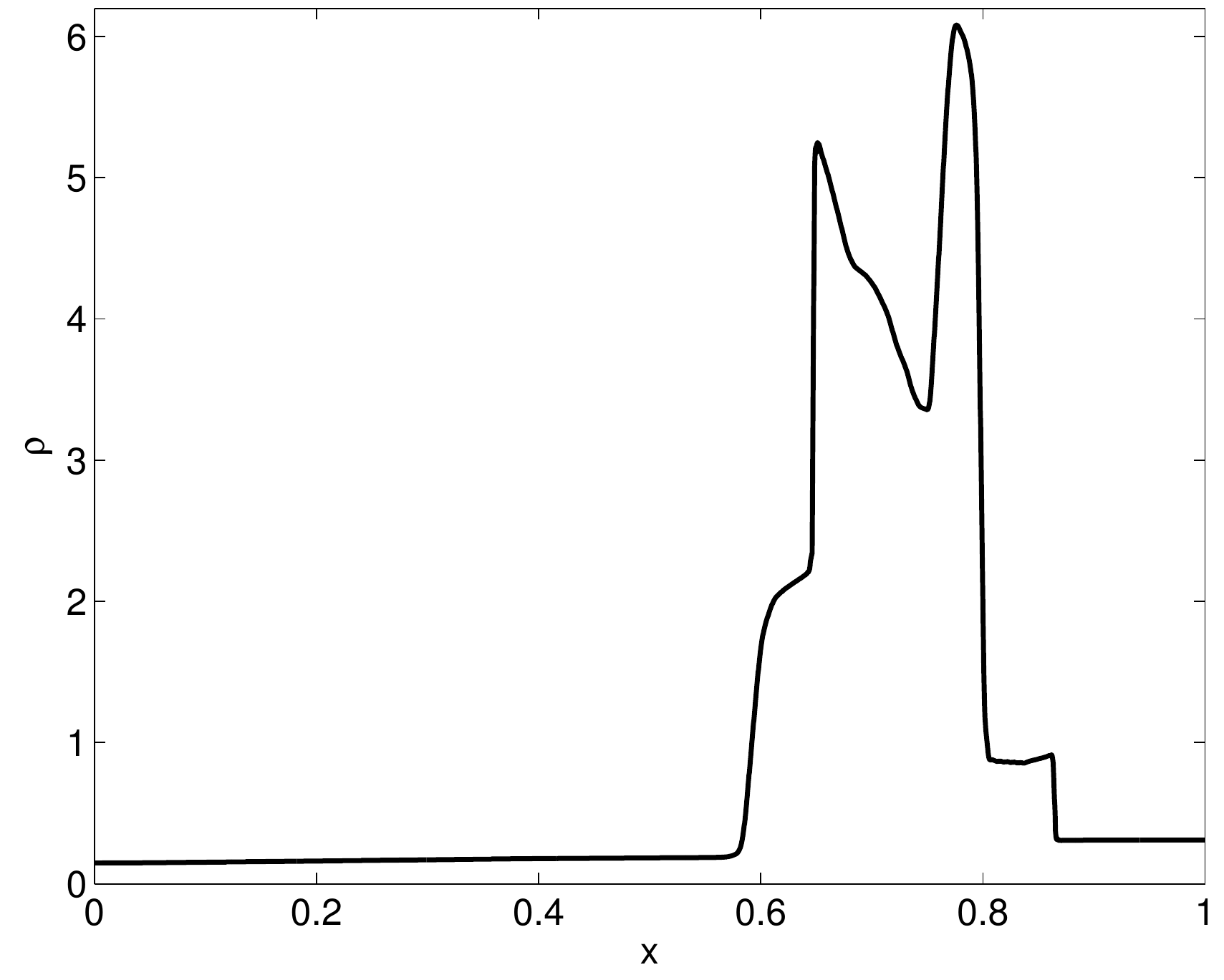}}
\caption{Approximation at $T=0.038$, KXRCF troubled-cell indicator (density and entropy), Blast, 512 elements.}\label{fig:BlastKHsol}
\end{figure}

Time history plots of detected troubled cells using the multiwavelet troubled-cell indicator with 512 elements and $k=1$ or $k=2$ can be seen in Figure \ref{fig:BlastCn9}, with corresponding approximation at $T=0.038$ in Figure \ref{fig:BlastCn9sol}. Note that although this is an extremely nonlinear problem, only a few elements elements should be limited in order to get nonoscillatory results. Our parameter $C$ is a useful tool to prevent limiting too many elements. The KXRCF indicator, however, selects more elements, as shown in Figures \ref{fig:BlastKH}, and \ref{fig:BlastKHsol}.  Both the multiwavelet and the KXRCF indicator detect regions that are visible in the exact shock solution, which was given by Woodward et al. \cite{Woo84C}. We speculate that the multiwavelet indicator will detect the same regions as KXRCF if a smaller $C$ is chosen. Note that this is the only example where different values of $C$ are used for the linear and quadratic case. This is due to the fact that  the interacting blast waves form an extremely nonlinear problem, thereby requiring a very accurate choice of $C$.

\subsubsection{Shock density wave interaction problem}
\label{sec:sine}

\begin{figure}[ht!]
\centering
 \includegraphics[scale = 0.3]{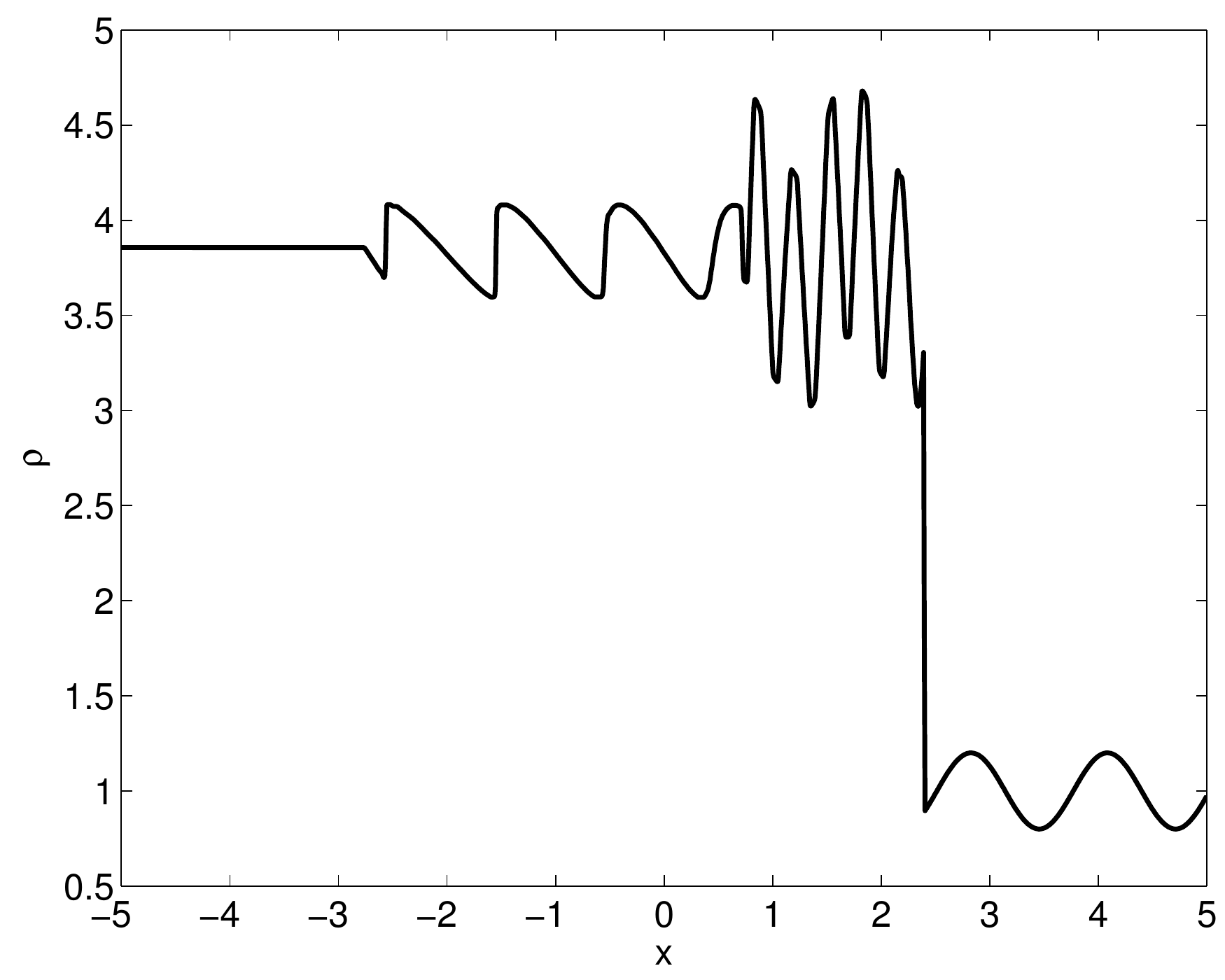}
\caption{'Exact' solution of the Shu-Osher problem at $T=1.8$.} \label{fig:exactsine}
\end{figure}

The final set of initial conditions that we consider for the one-dimensional Euler equations is given by Shu et al. \cite{Shu89O}:
\begin{subequations}
\begin{equation}
 \rho(x,0) = \left\{ \begin{array}{ll}3.857143, & x < -4,\\ 1+0.2 \sin(5x), & x \geq -4, \end{array} \right.
\end{equation}
\begin{equation}
  u(x,0) = \left\{ \begin{array}{ll}2.629369, & x < -4, \\ 0, & x \geq -4,\end{array} \right. P(x,0) = \left\{ \begin{array}{ll}10.33333, & x < -4, \\ 1, & x \geq -4, \end{array} \right.
\end{equation}
\end{subequations}
together with constant boundary conditions.  The exact solution at $T=1.8$ is approximated using a fine mesh, and is shown in Figure \ref{fig:exactsine}. Here, we see that the discontinuity in the initial condition is still apparent, and some shocks are formed in the left part of the solution. 

\begin{figure}[ht!]
\centering
 \subfigure[$k=1$, limited elements]{\includegraphics[scale = 0.25]{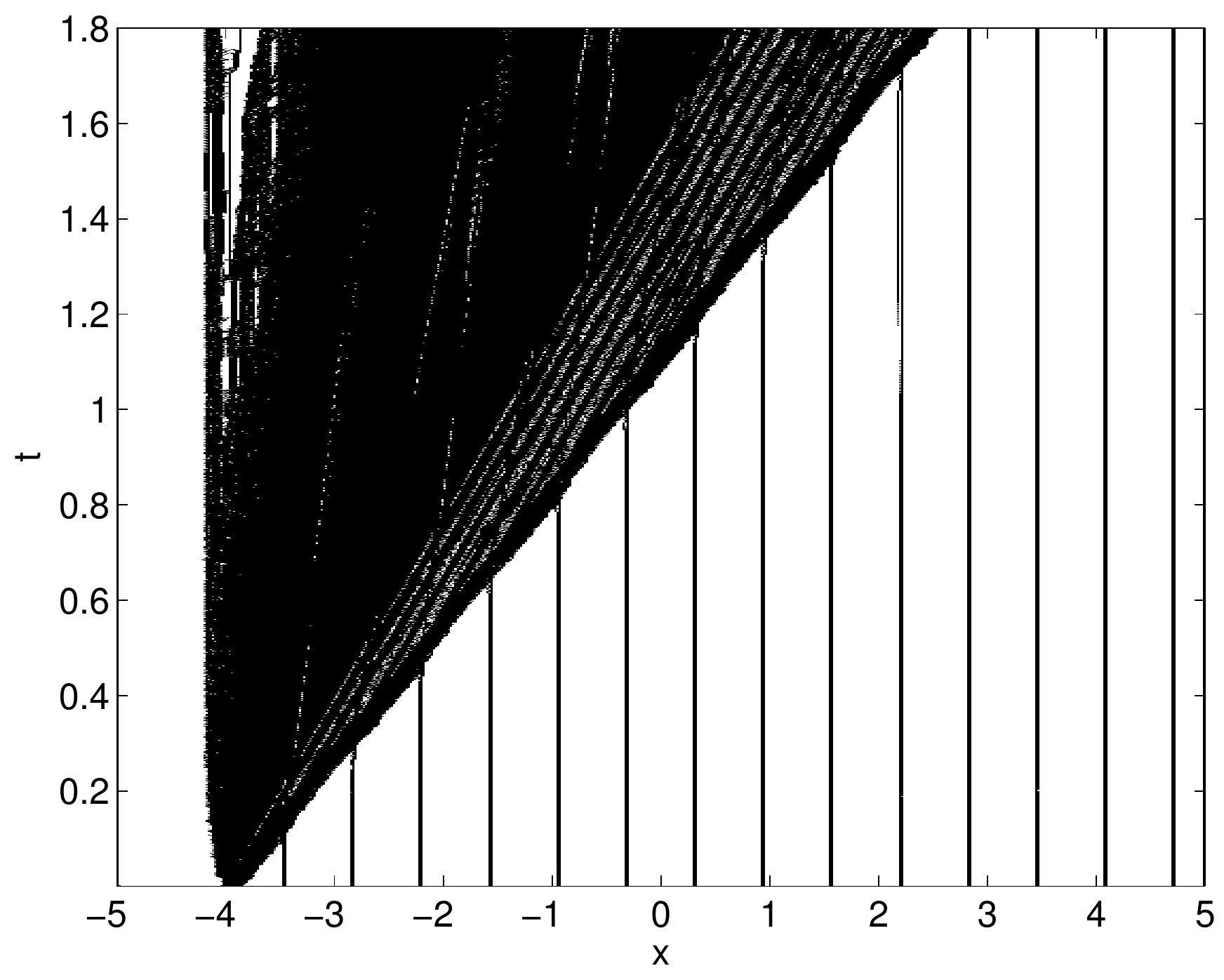}}
 \subfigure[$k=2$, limited elements]{\includegraphics[scale = 0.25]{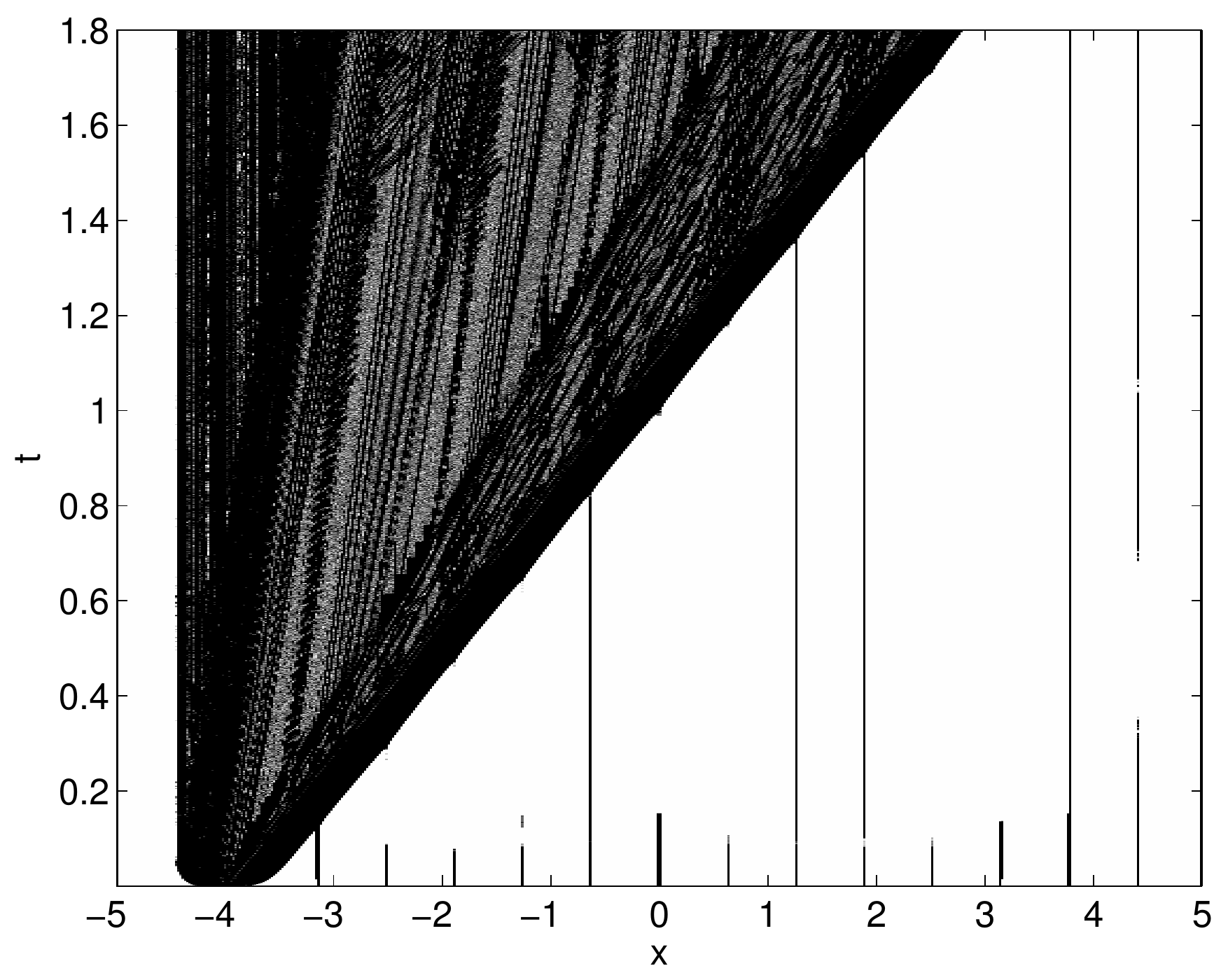}} \\
 \subfigure[$k=1$, solution]{\includegraphics[scale = 0.25]{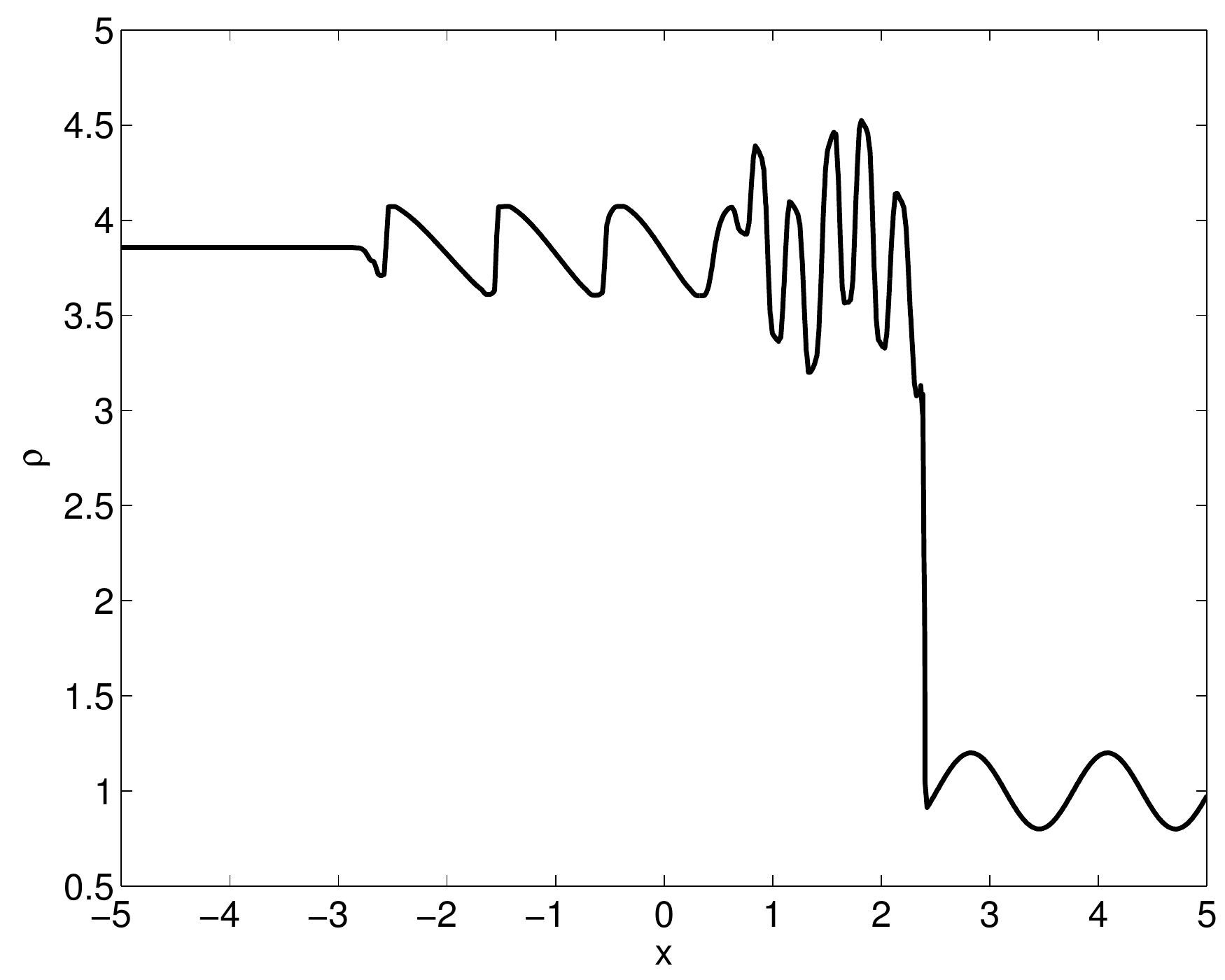}}
 \subfigure[$k=2$, solution]{\includegraphics[scale = 0.25]{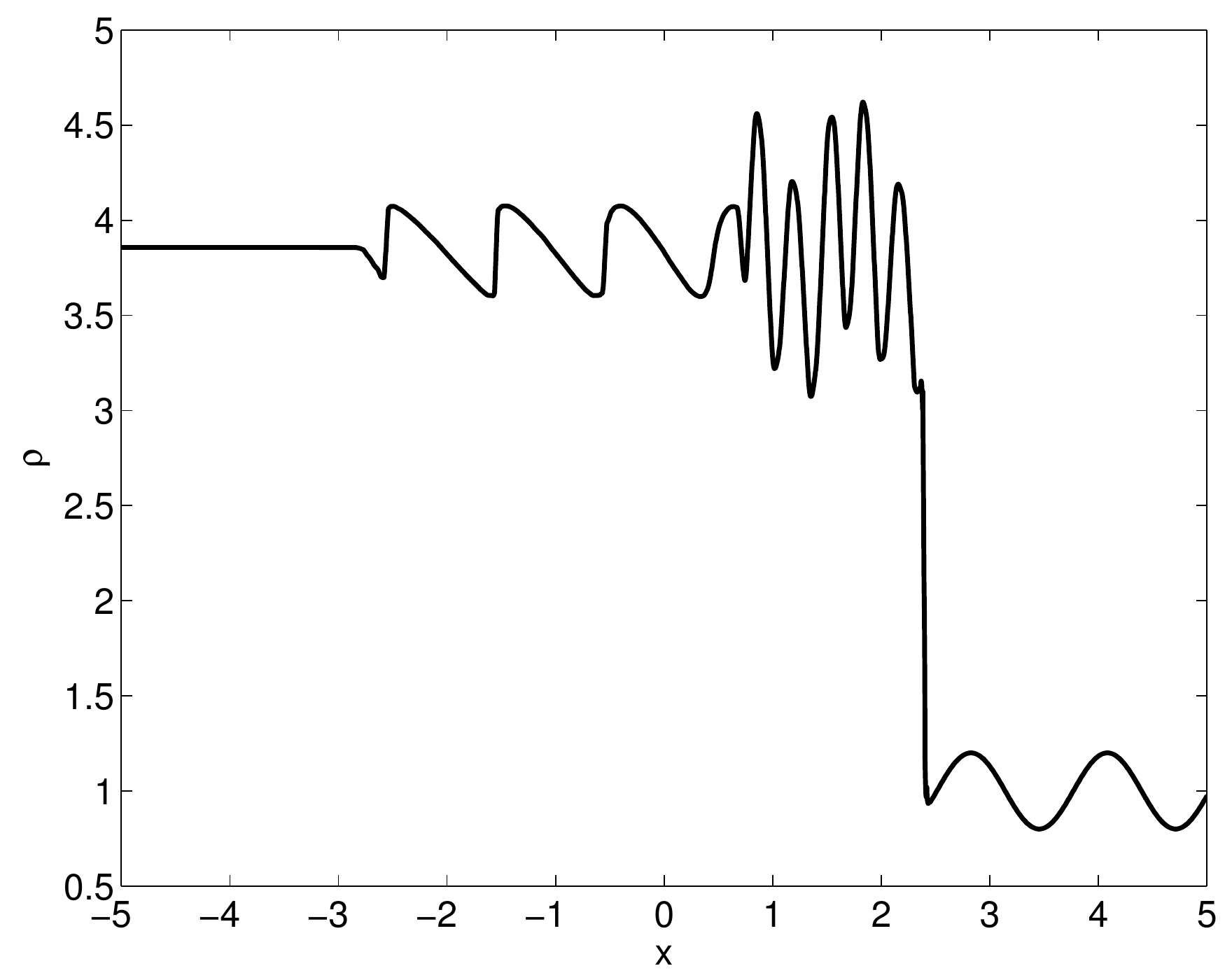}} \\
\caption{Solution at $T=1.8$ and time history plot of limited cells, using the unmodified moment limiter, 512 elements, Shu-Osher problem.}\label{fig:Sinemoment}
\end{figure}

\begin{figure}[ht!]
\centering
 \subfigure[$C=0.5$]{\includegraphics[scale = 0.22]{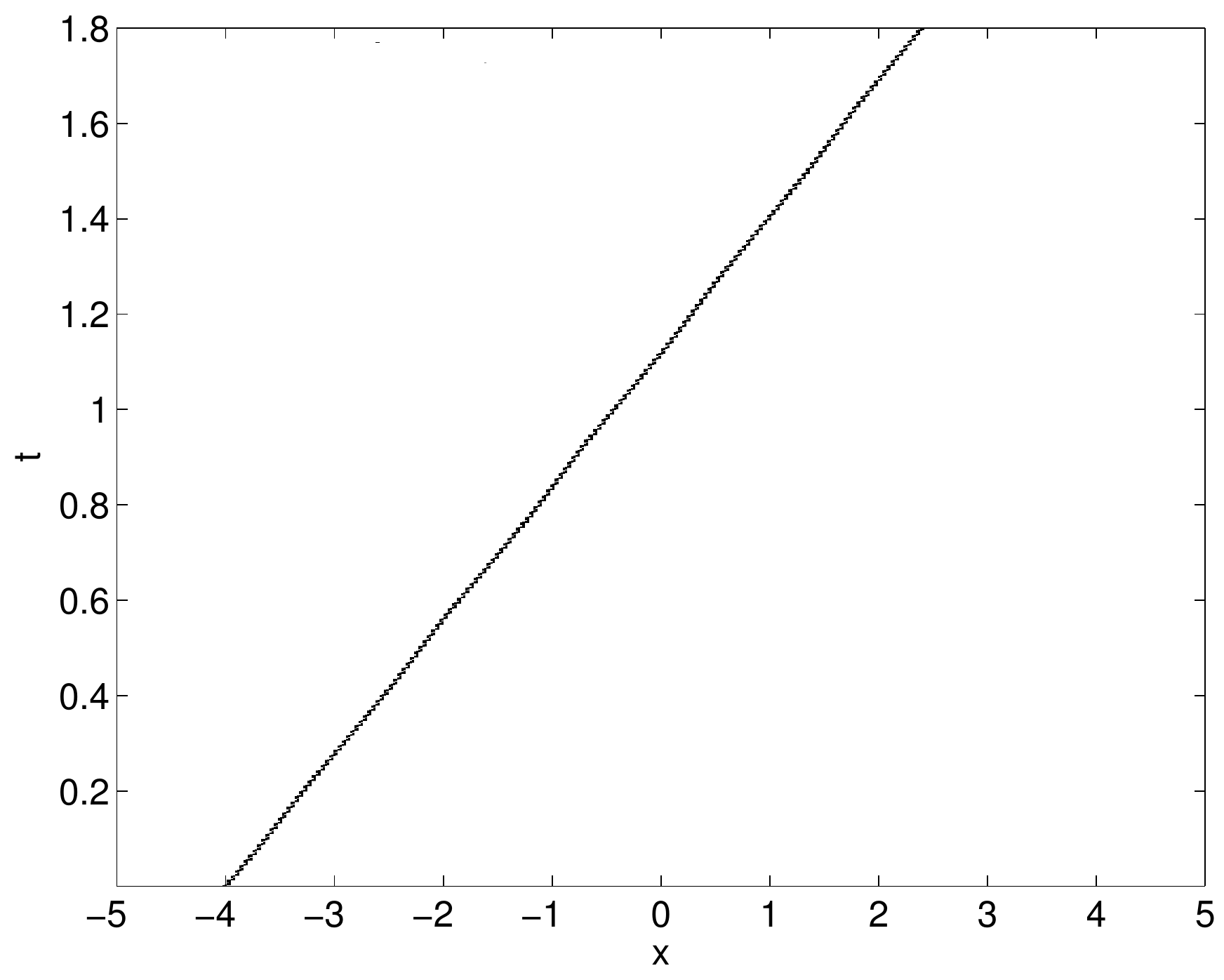}}
 \subfigure[$C=0.1$]{\includegraphics[scale = 0.22]{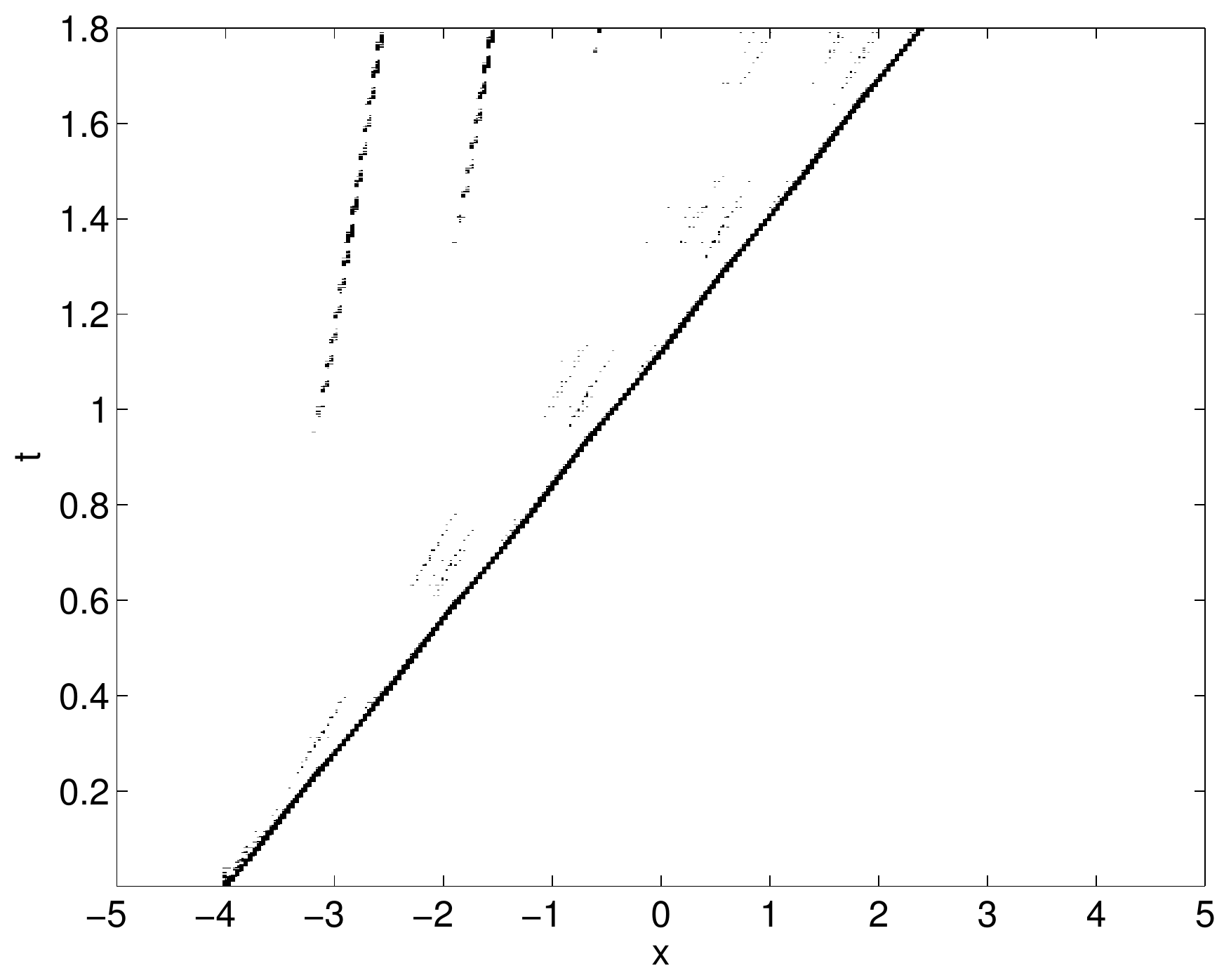}}
 \subfigure[$C=0.05$]{\includegraphics[scale = 0.22]{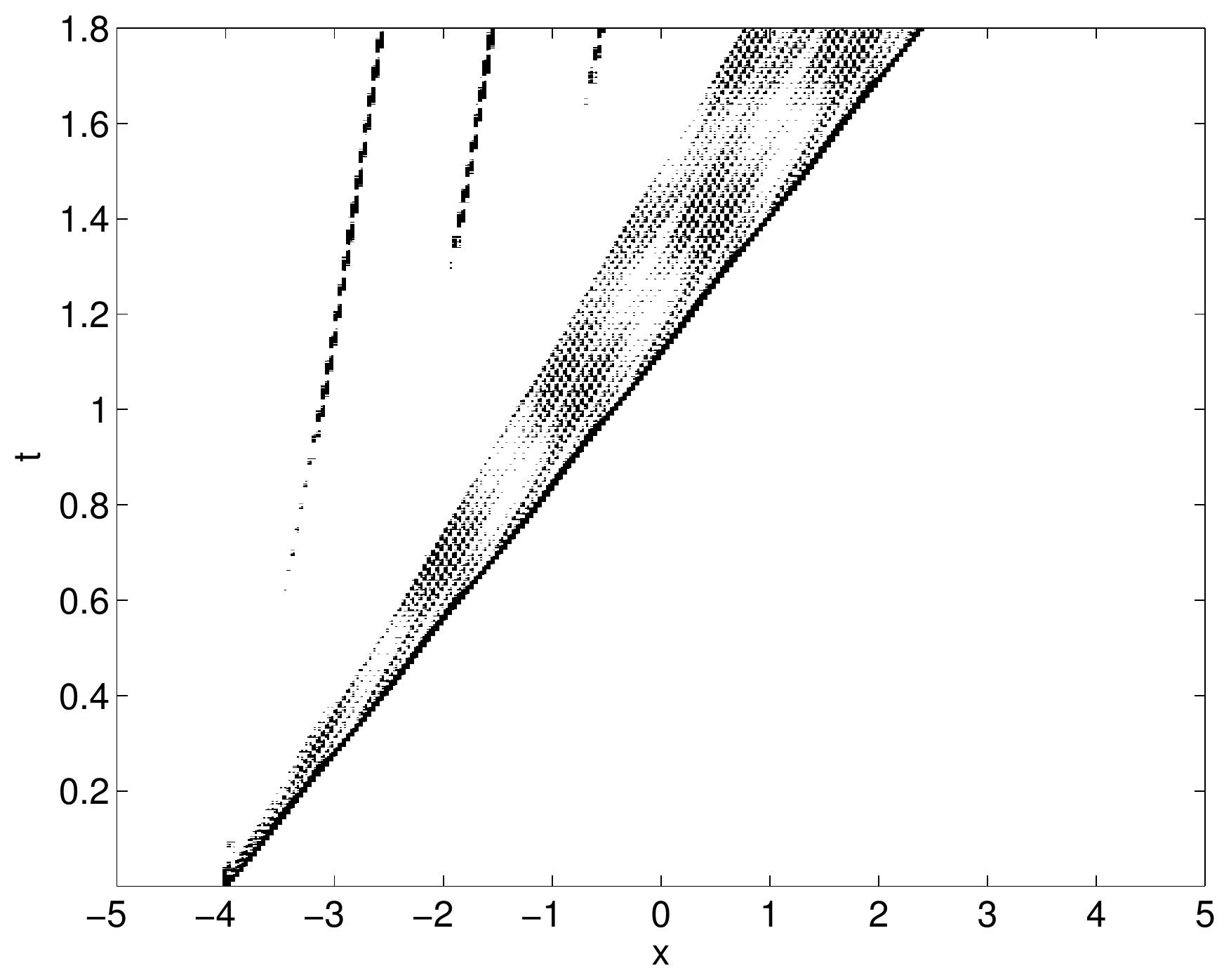}} \\
 \subfigure[$C=0.5$]{\includegraphics[scale = 0.22]{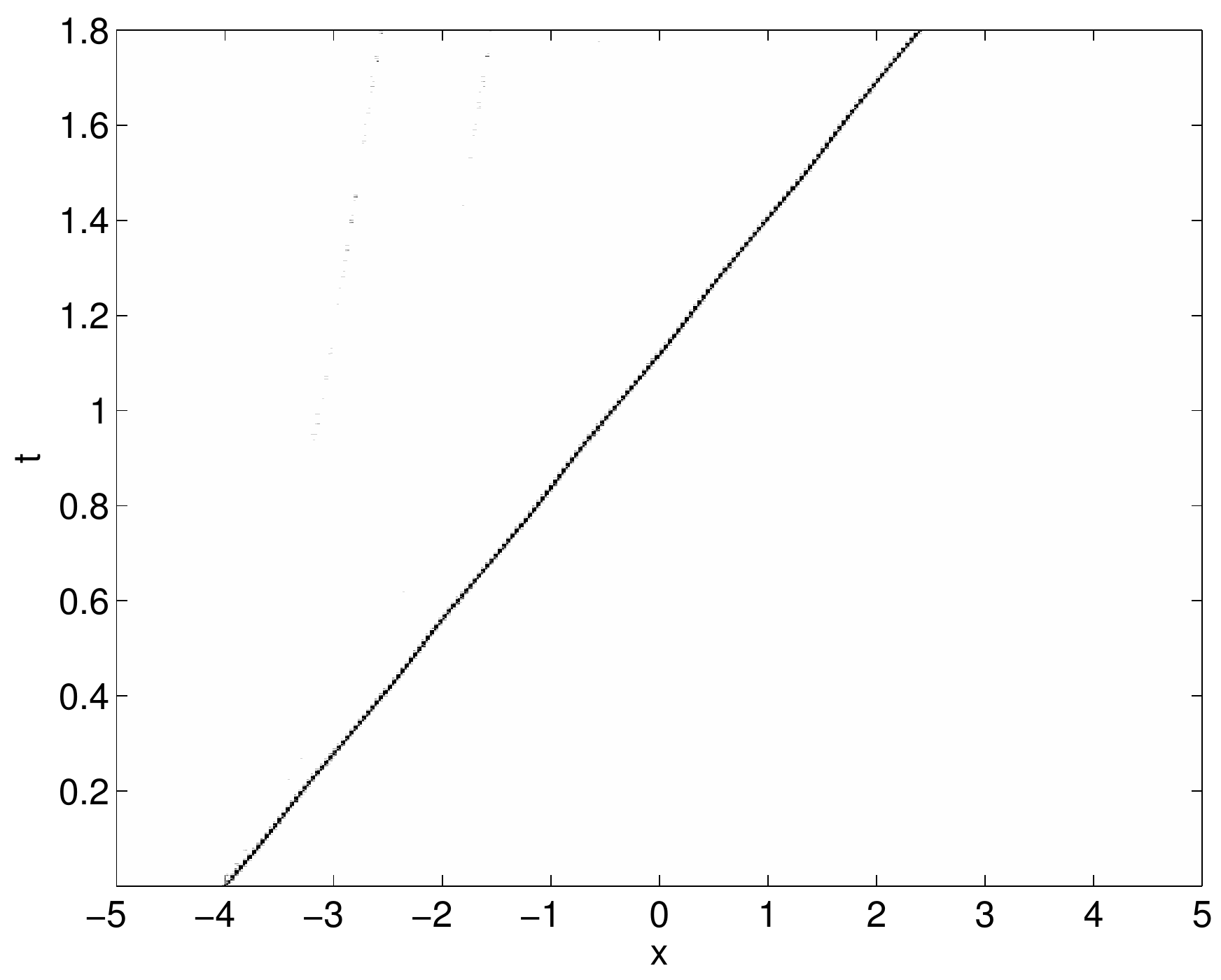}} 
 \subfigure[$C=0.1$]{\includegraphics[scale = 0.22]{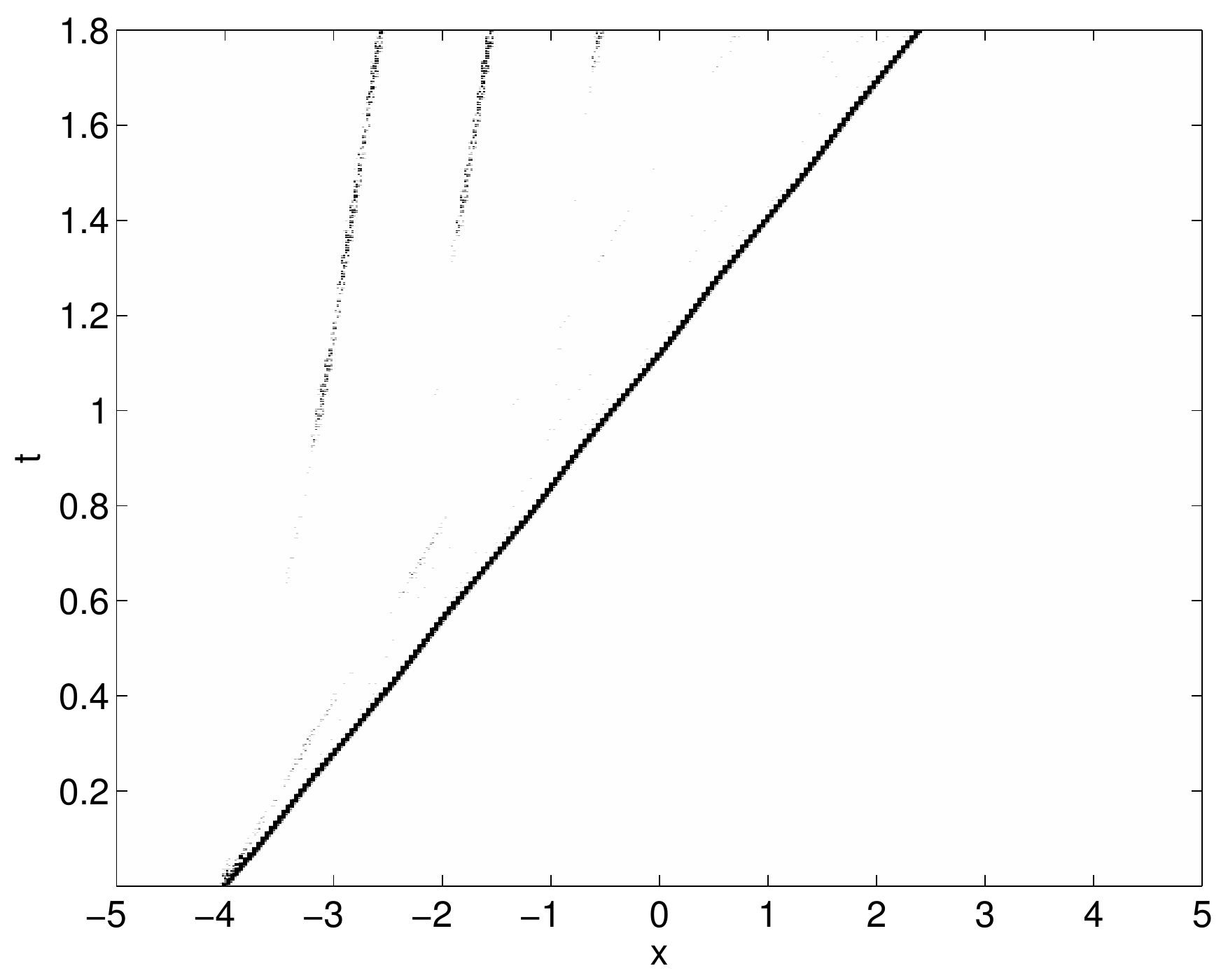}} 
 \subfigure[$C=0.05$]{\includegraphics[scale = 0.22]{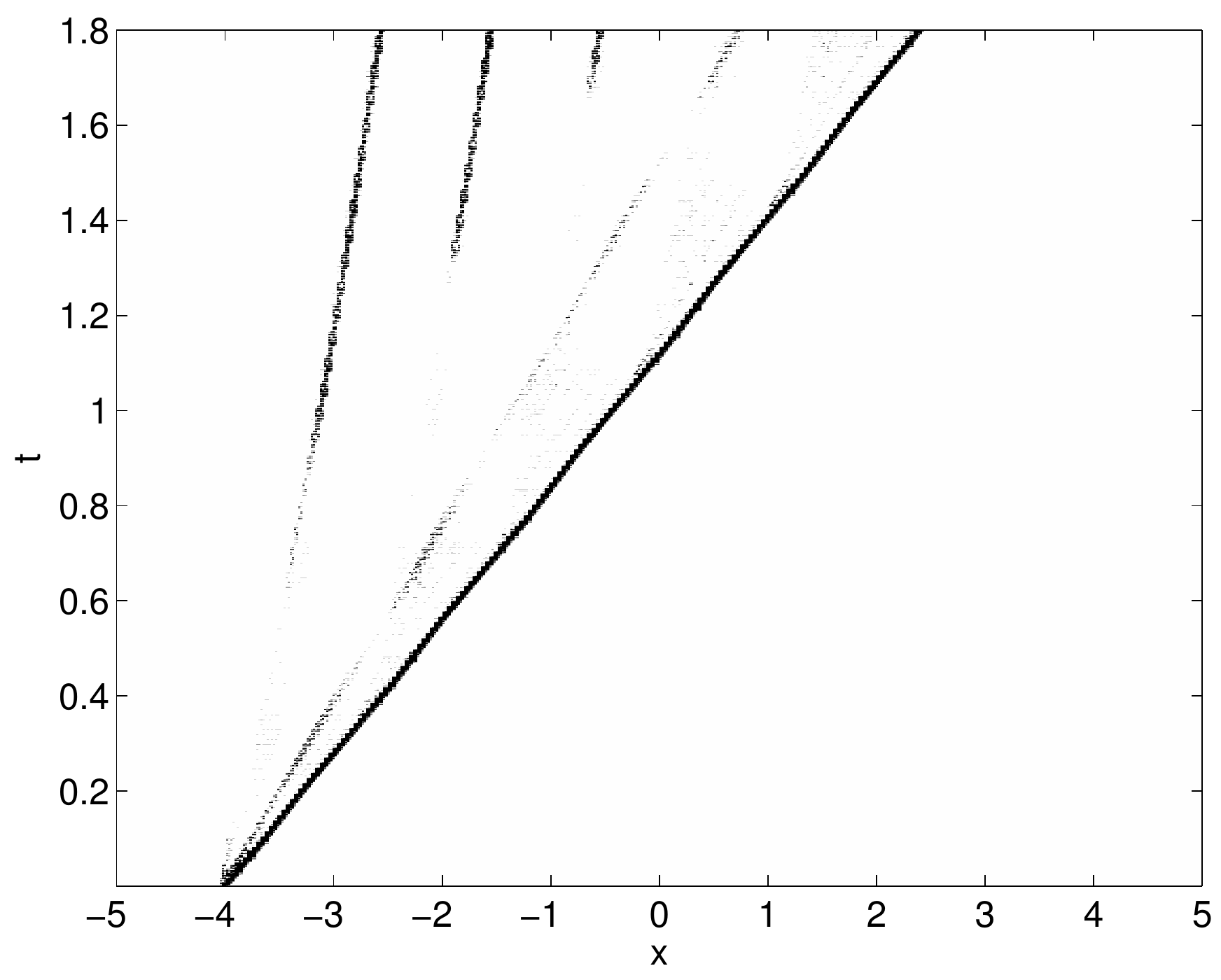}} \\
\caption{Time history plot of detected troubled cells, using multiwavelets (density), Shu-Osher problem, 512 elements. First row: $k=1$, second row: $k=2$.}\label{fig:SineC}
\end{figure}

\begin{figure}[ht!]
 \centering
\subfigure[KXRCF, $k=1$]{\includegraphics[scale = 0.28]{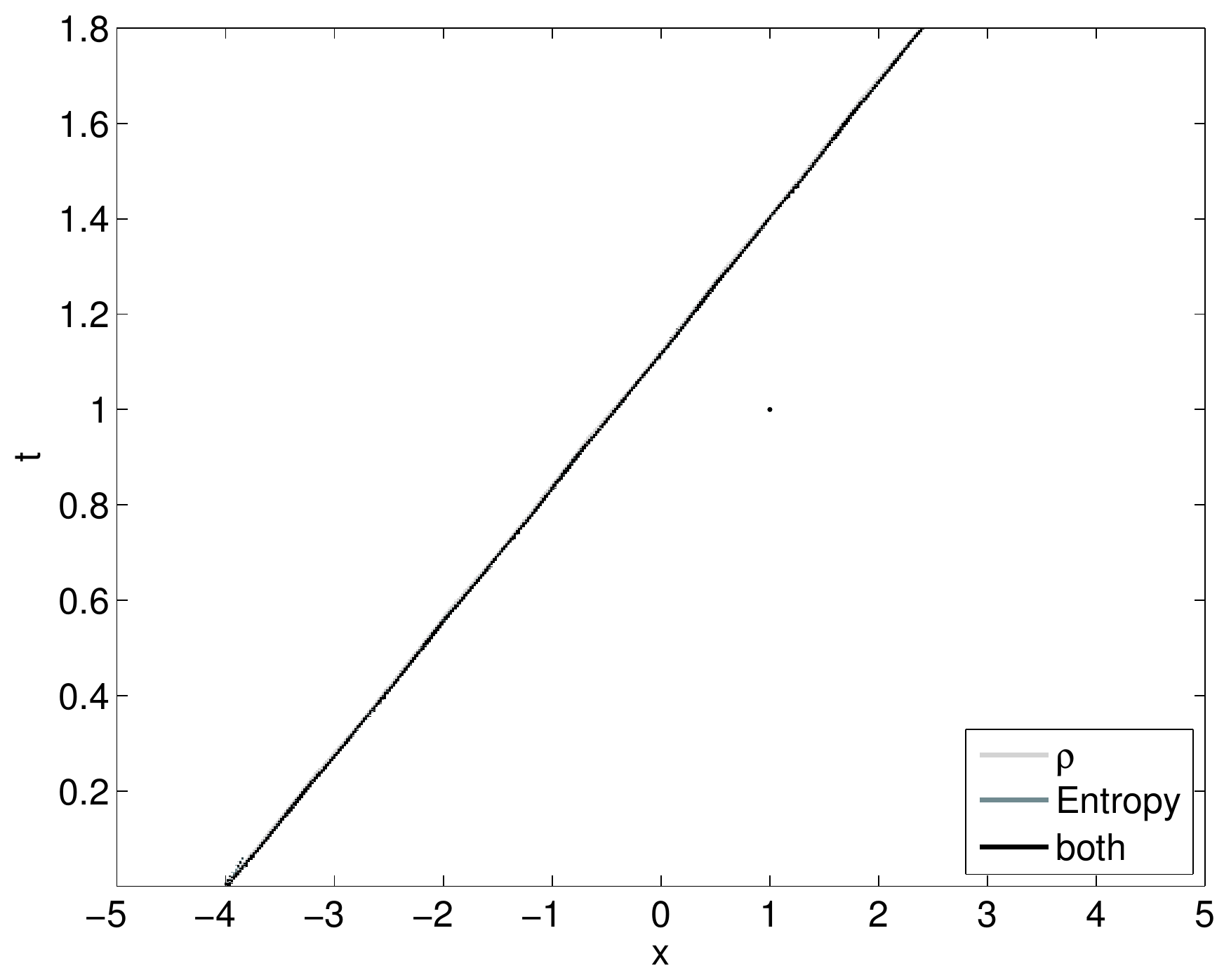}}
\subfigure[Harten, $k=1$]{\includegraphics[scale = 0.28]{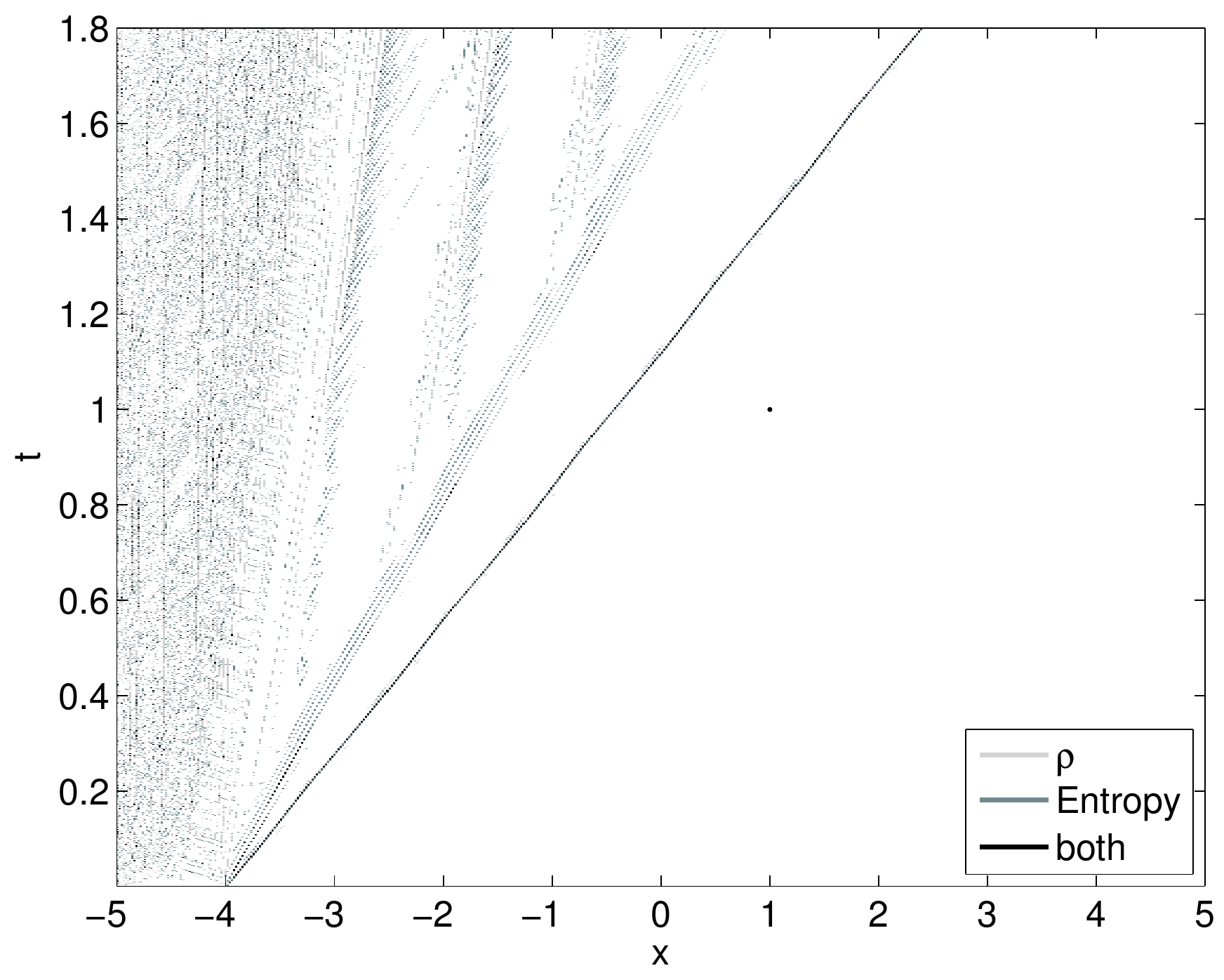}} \\
\vspace{-0.3cm}
\subfigure[KXRCF, $k=2$]{\includegraphics[scale = 0.28]{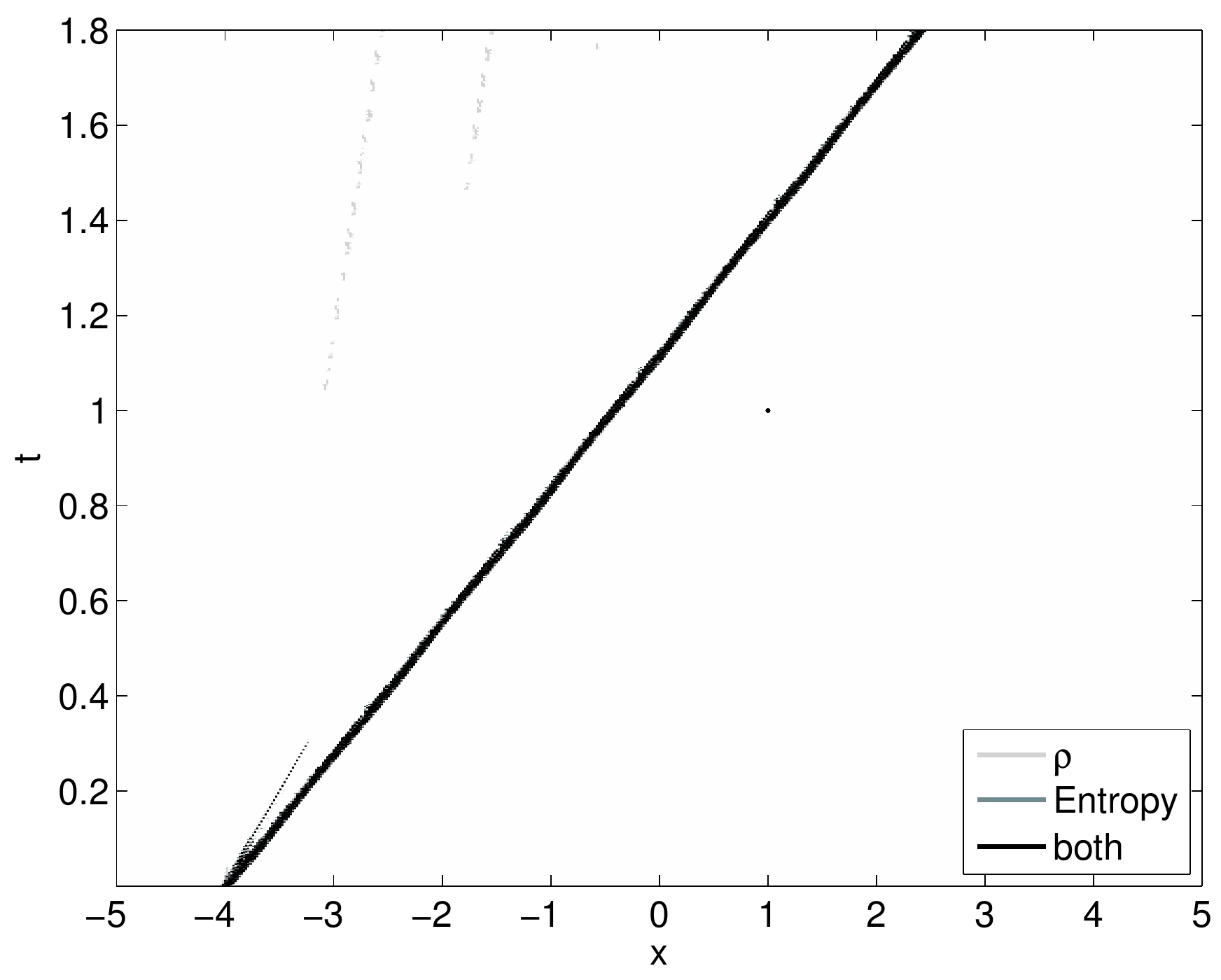}}
\subfigure[Harten, $k=2$]{\includegraphics[scale = 0.28]{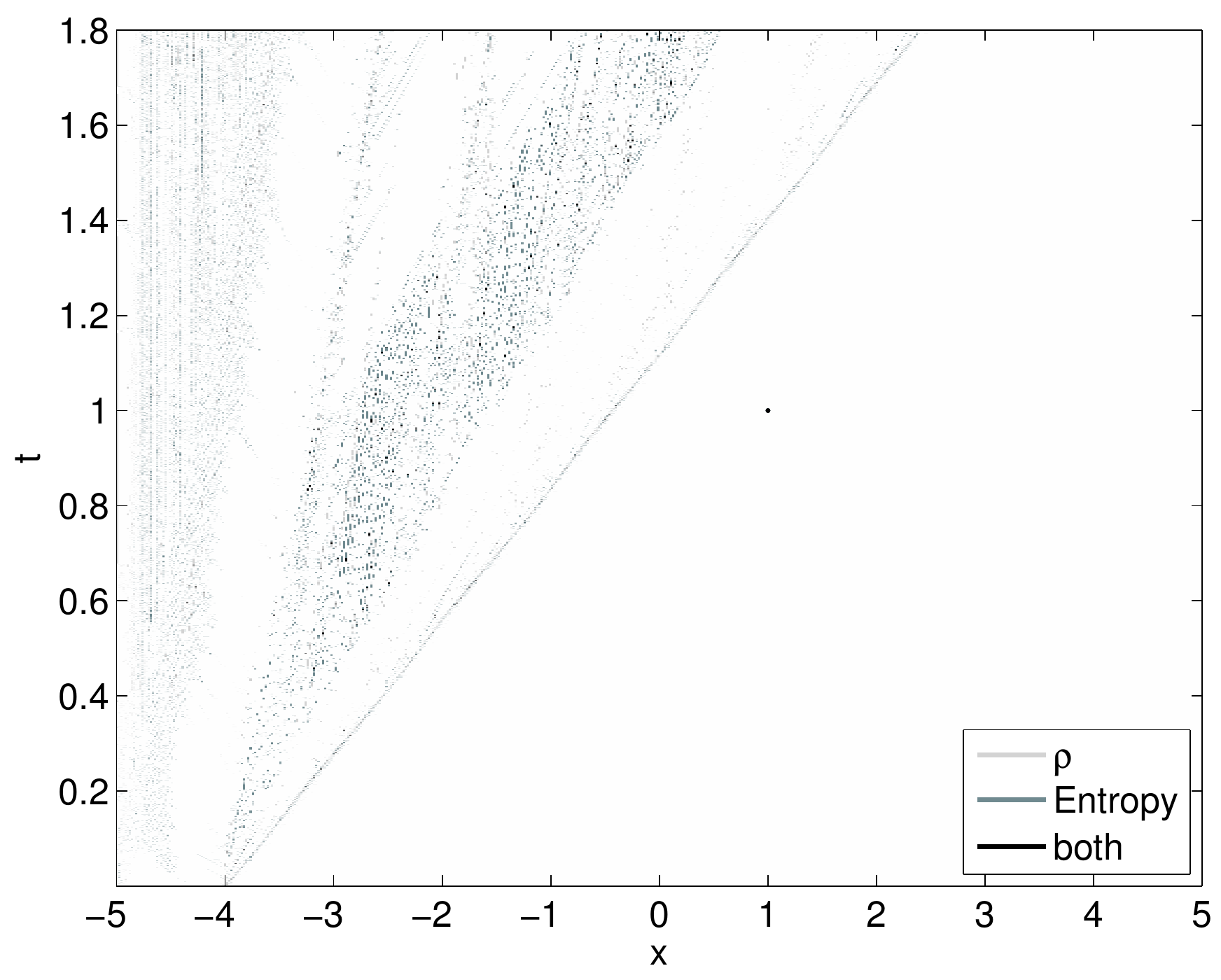}} \\
\vspace{-0.3cm}
\caption{Time history plot of detected troubled cells, using KXRCF or Harten ($\alpha = 1.5$), using density and entropy, Shu-Osher problem, 512 elements.}\label{fig:SineKH}
\end{figure}

\begin{figure}[ht!]
\centering
 \subfigure[$C=0.5$]{\includegraphics[scale = 0.22]{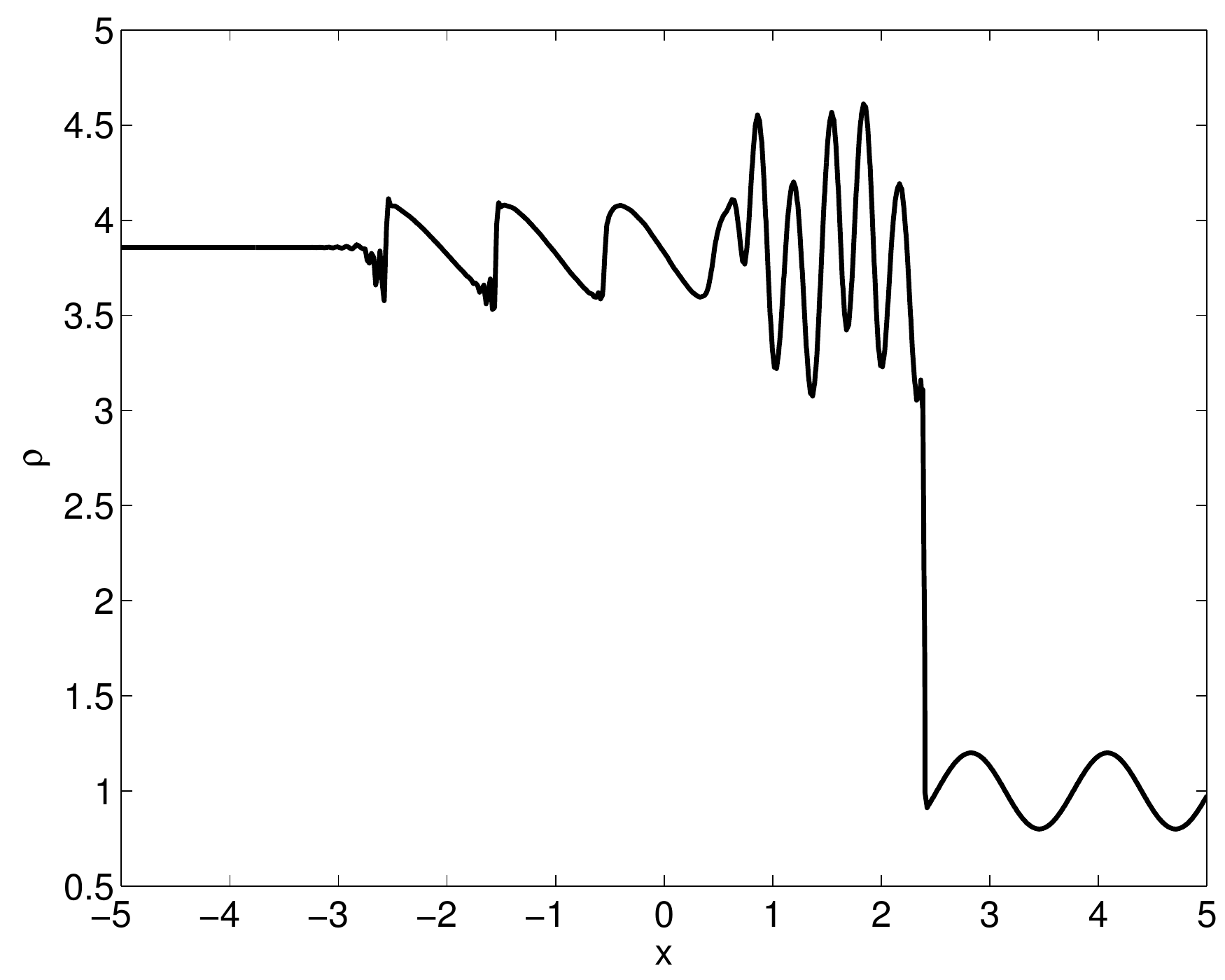}}
 \subfigure[$C=0.1$]{\includegraphics[scale = 0.22]{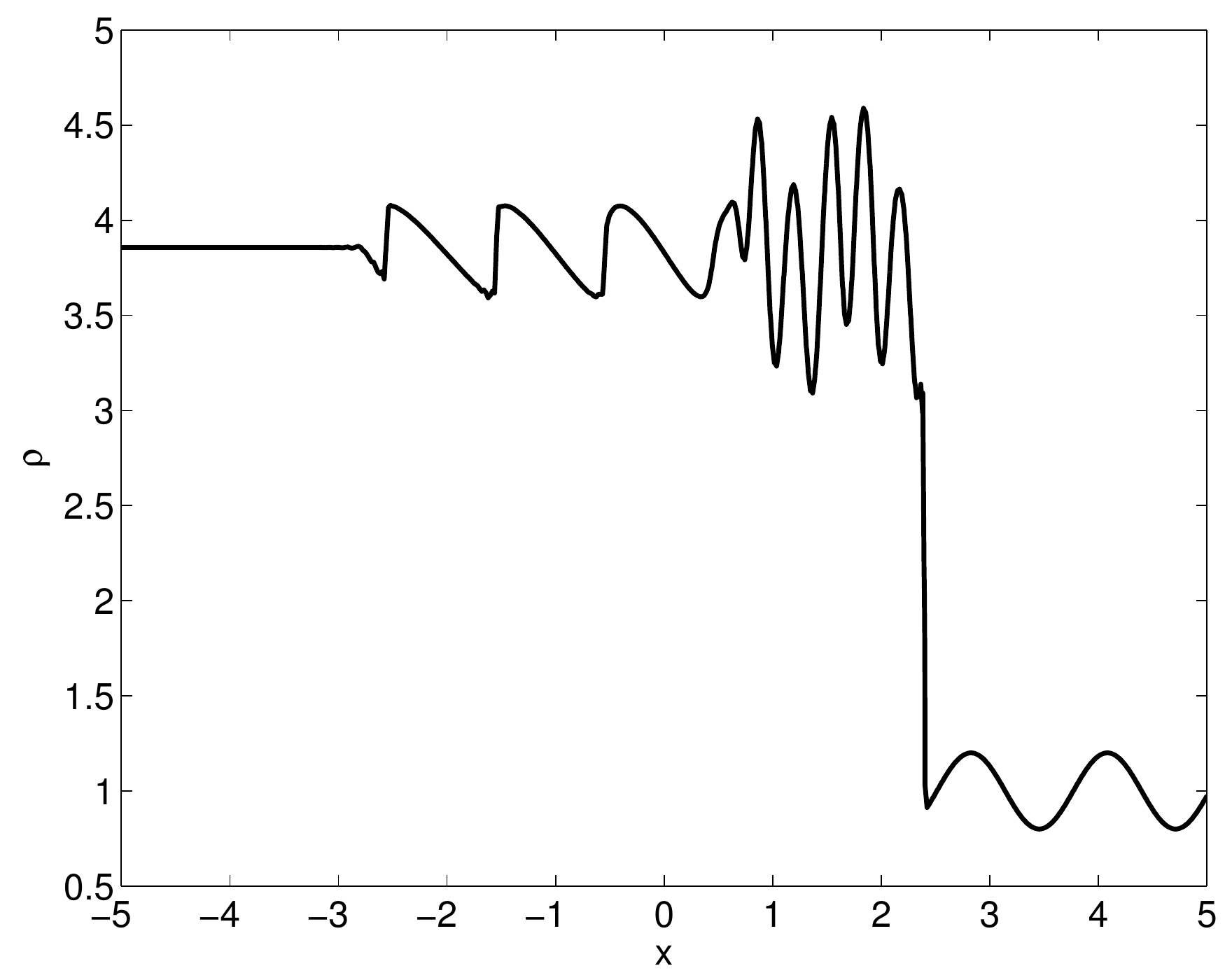}}
 \subfigure[$C=0.05$]{\includegraphics[scale = 0.22]{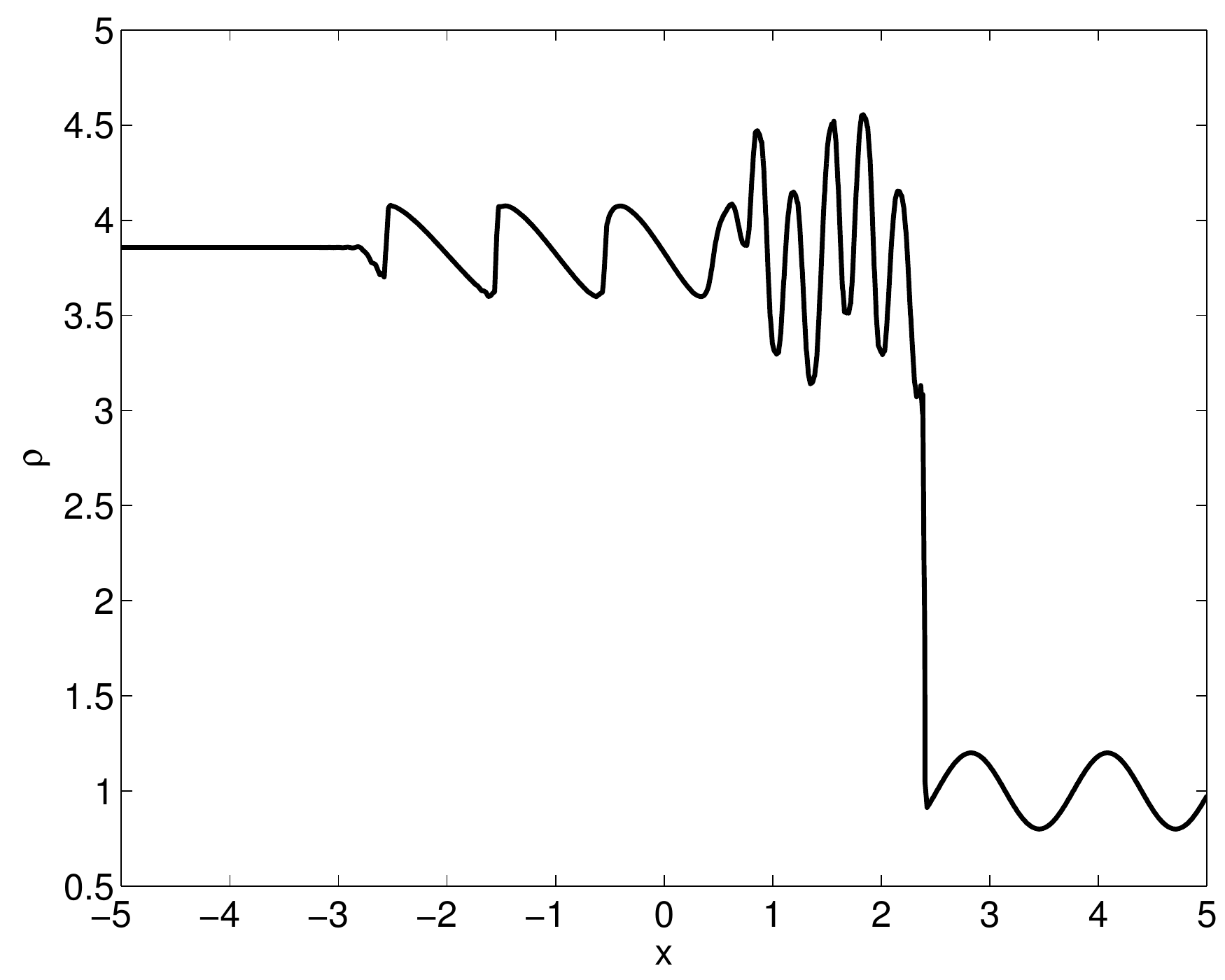}} \\
\vspace{-0.3cm}
 \subfigure[$C=0.5$]{\includegraphics[scale = 0.22]{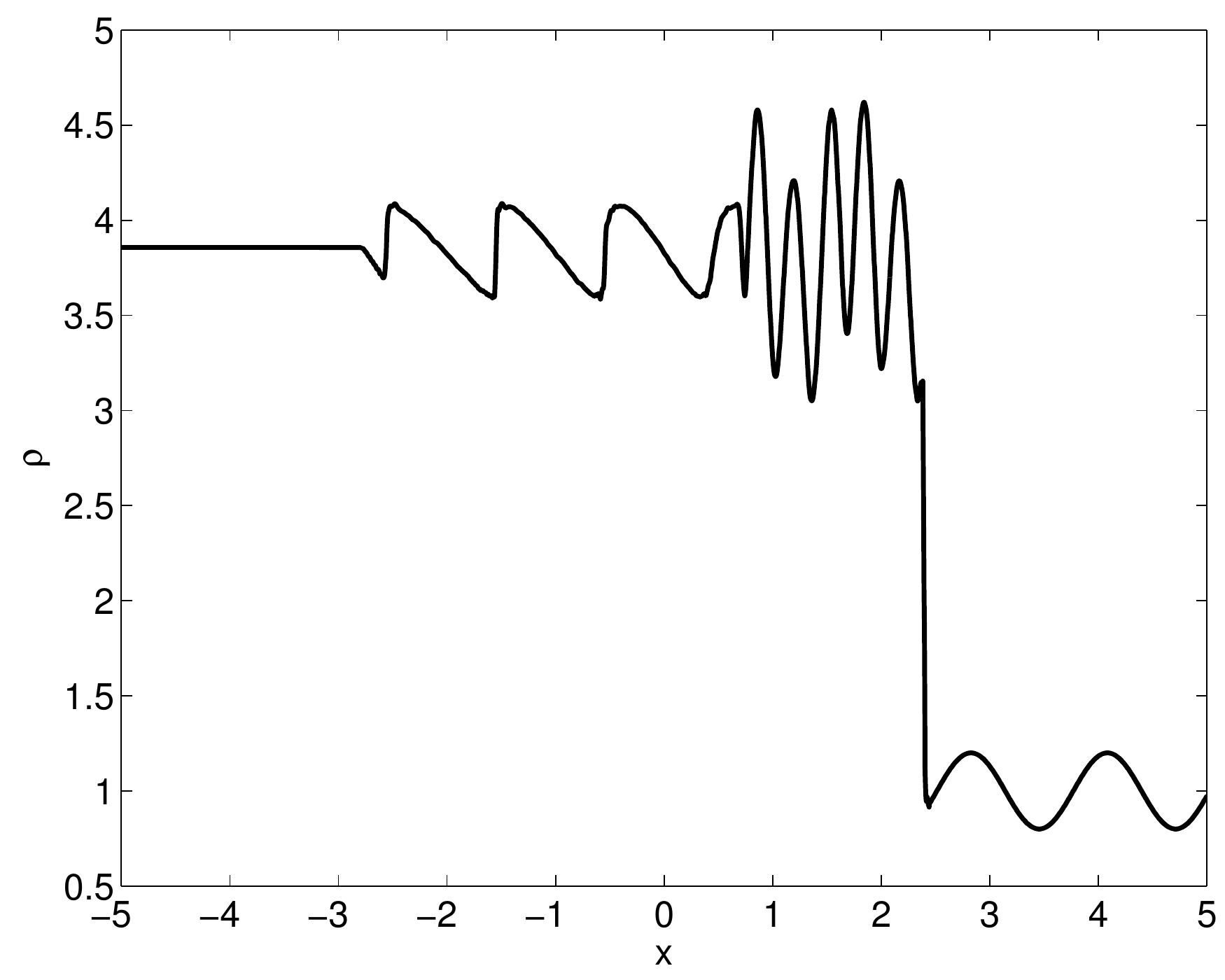}} 
 \subfigure[$C=0.1$]{\includegraphics[scale = 0.22]{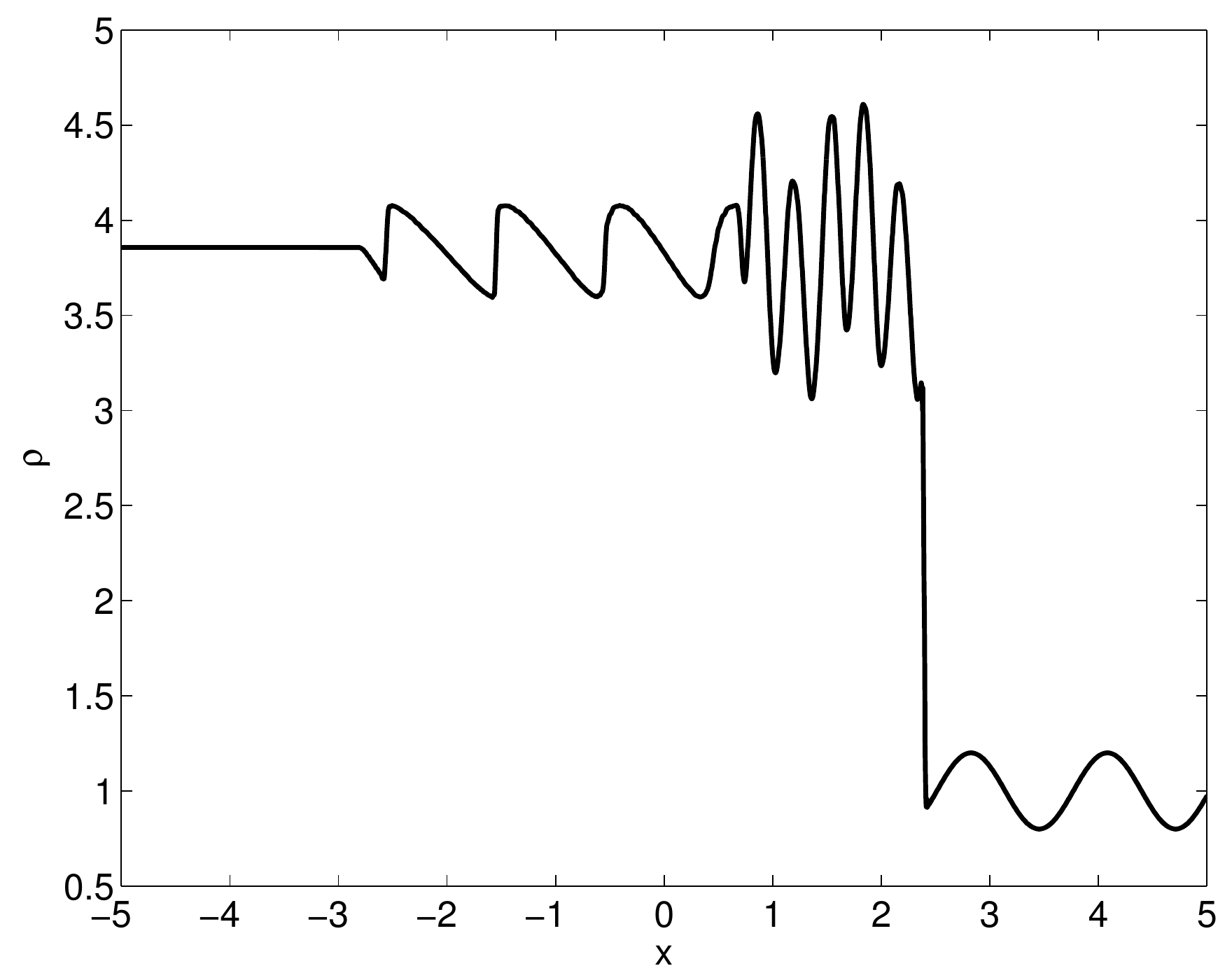}} 
 \subfigure[$C=0.05$]{\includegraphics[scale = 0.22]{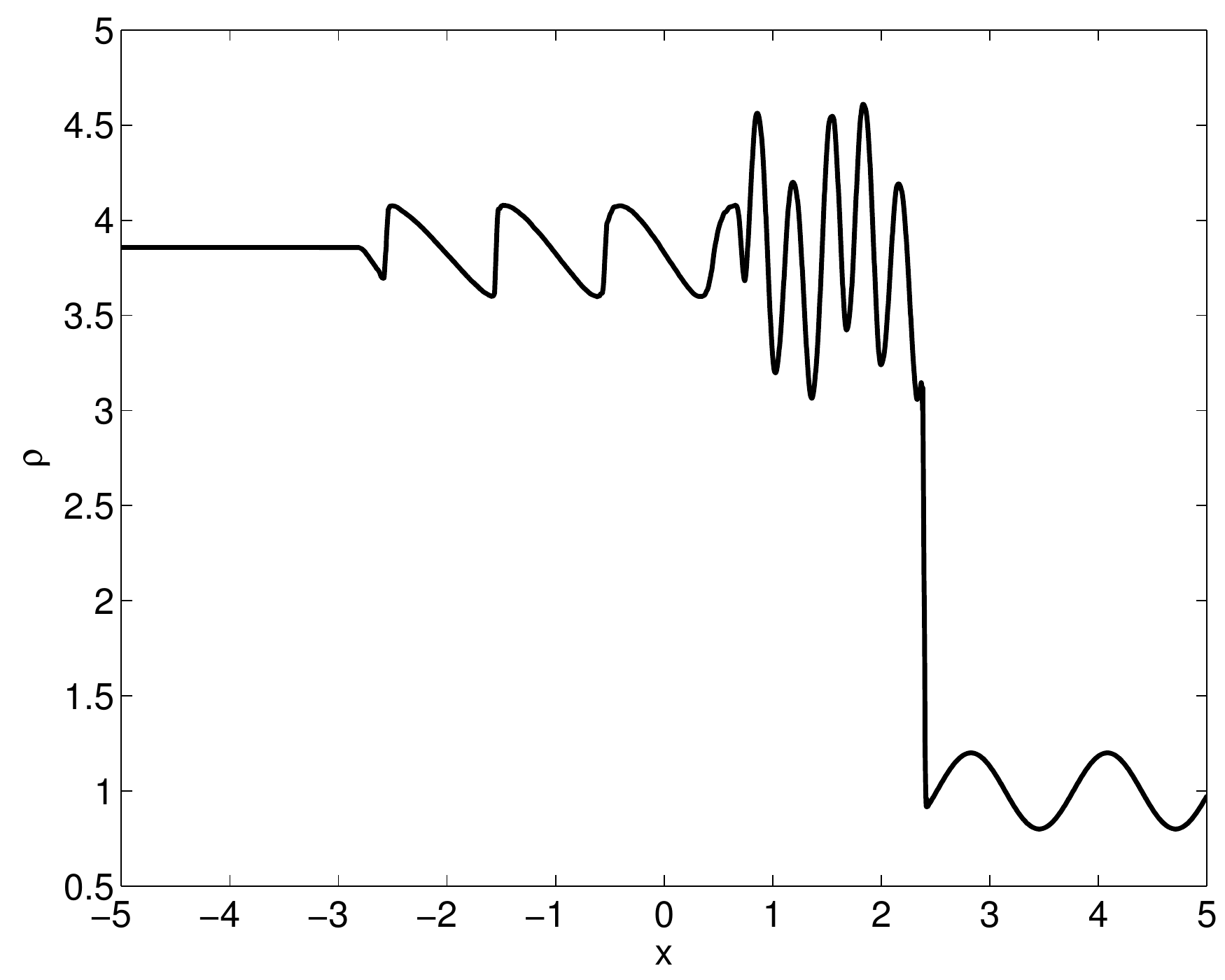}} \\
\vspace{-0.3cm}
\caption{Approximation at $T=1.8$, multiwavelet indicator (density), Shu-Osher problem, 512 elements. First row: $k=1$, second row: $k=2$.}\label{fig:SineCsol}
\end{figure}

\begin{figure}[ht!]
 \centering
\subfigure[KXRCF, $k=1$]{\includegraphics[scale = 0.28]{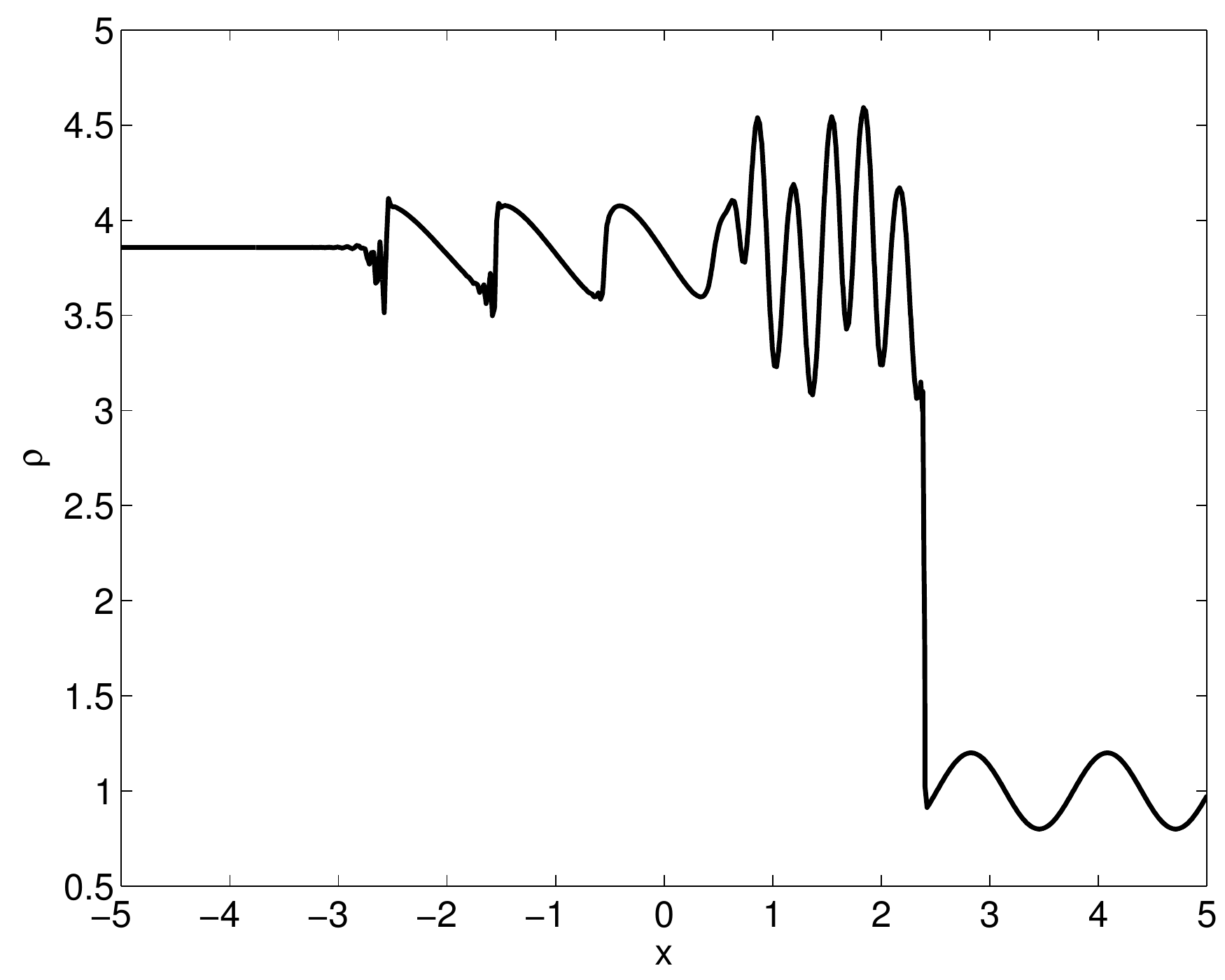}}
\subfigure[Harten, $k=1$]{\includegraphics[scale = 0.28]{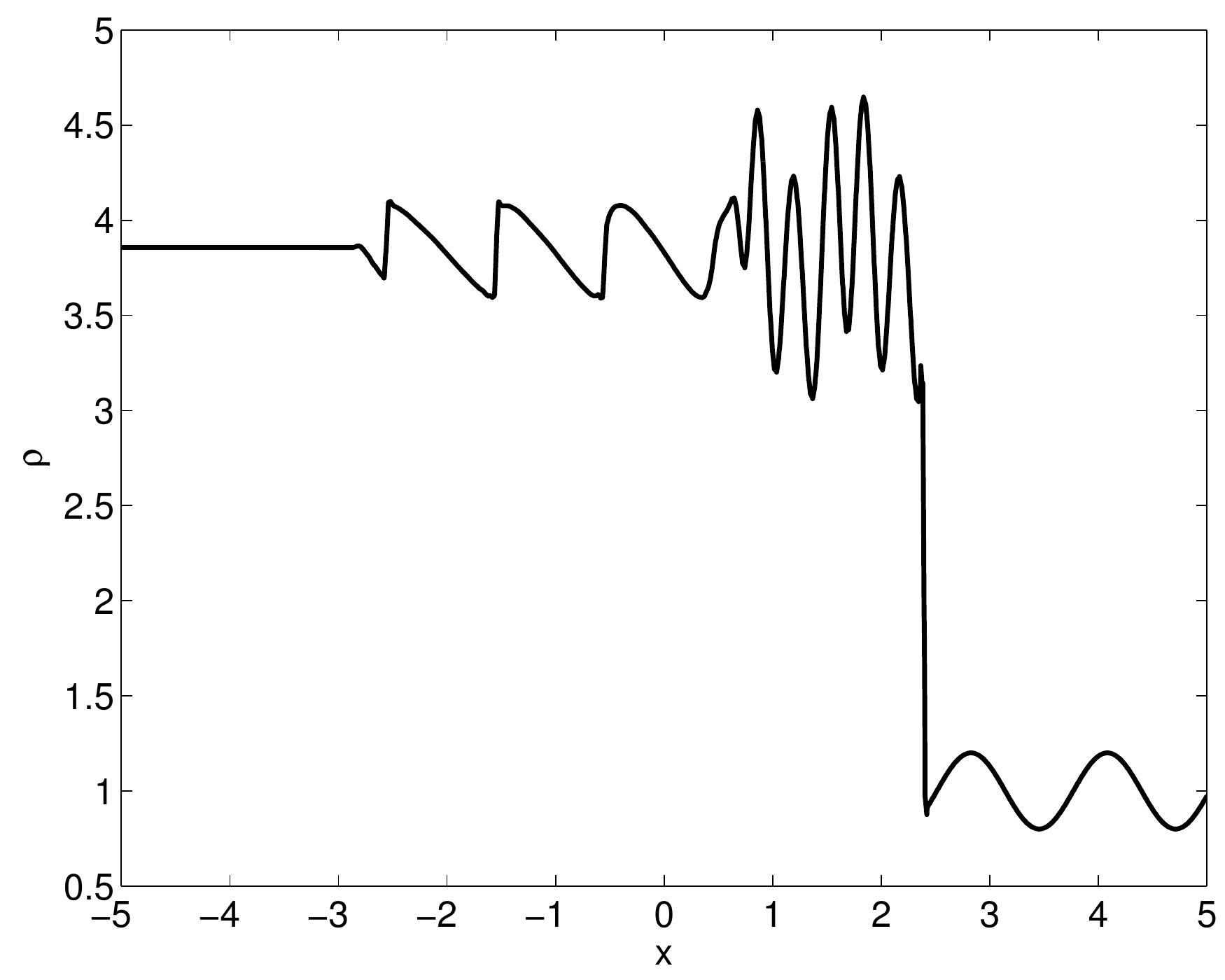}} \\
\subfigure[KXRCF, $k=2$]{\includegraphics[scale = 0.28]{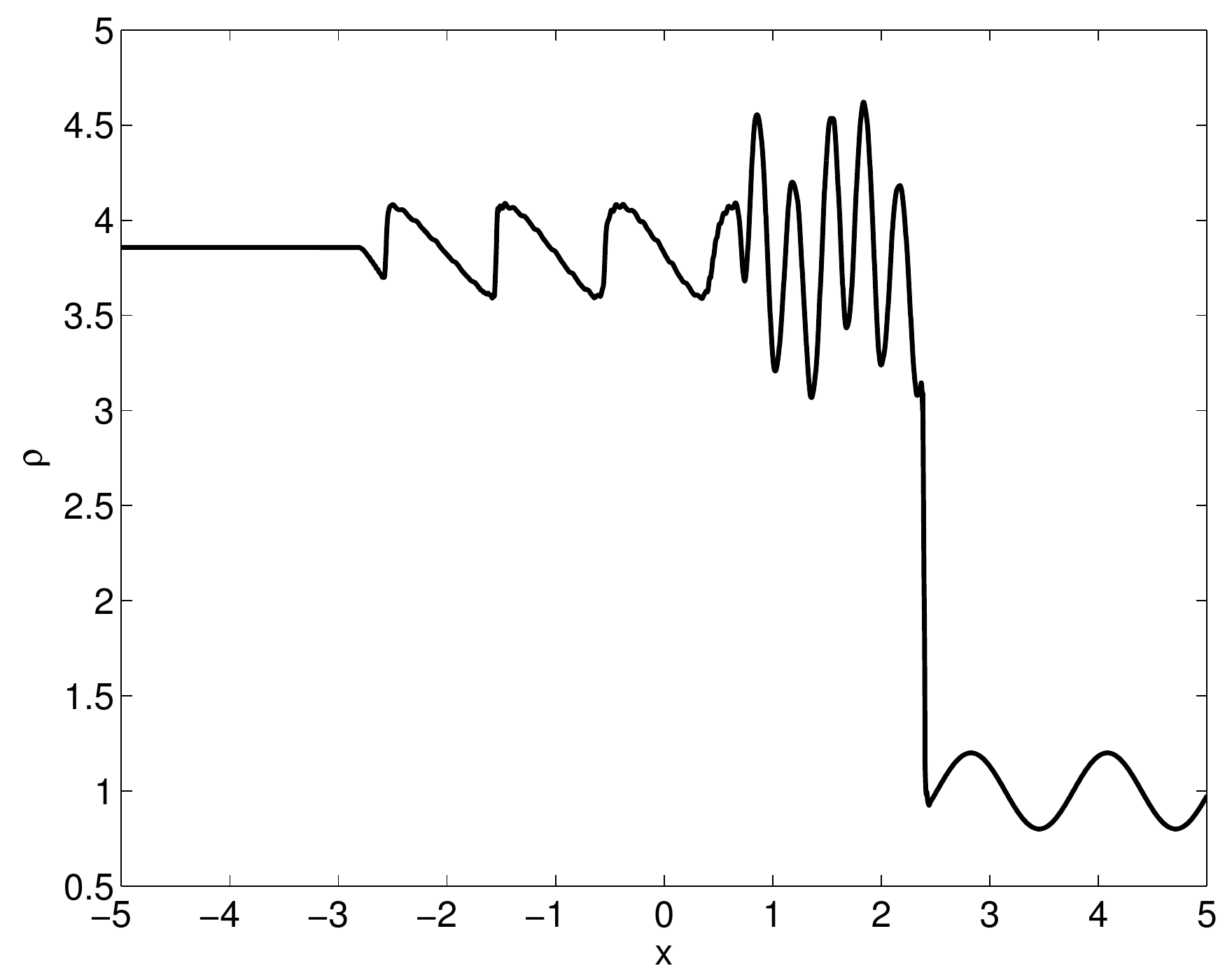}}
\subfigure[Harten, $k=2$]{\includegraphics[scale = 0.28]{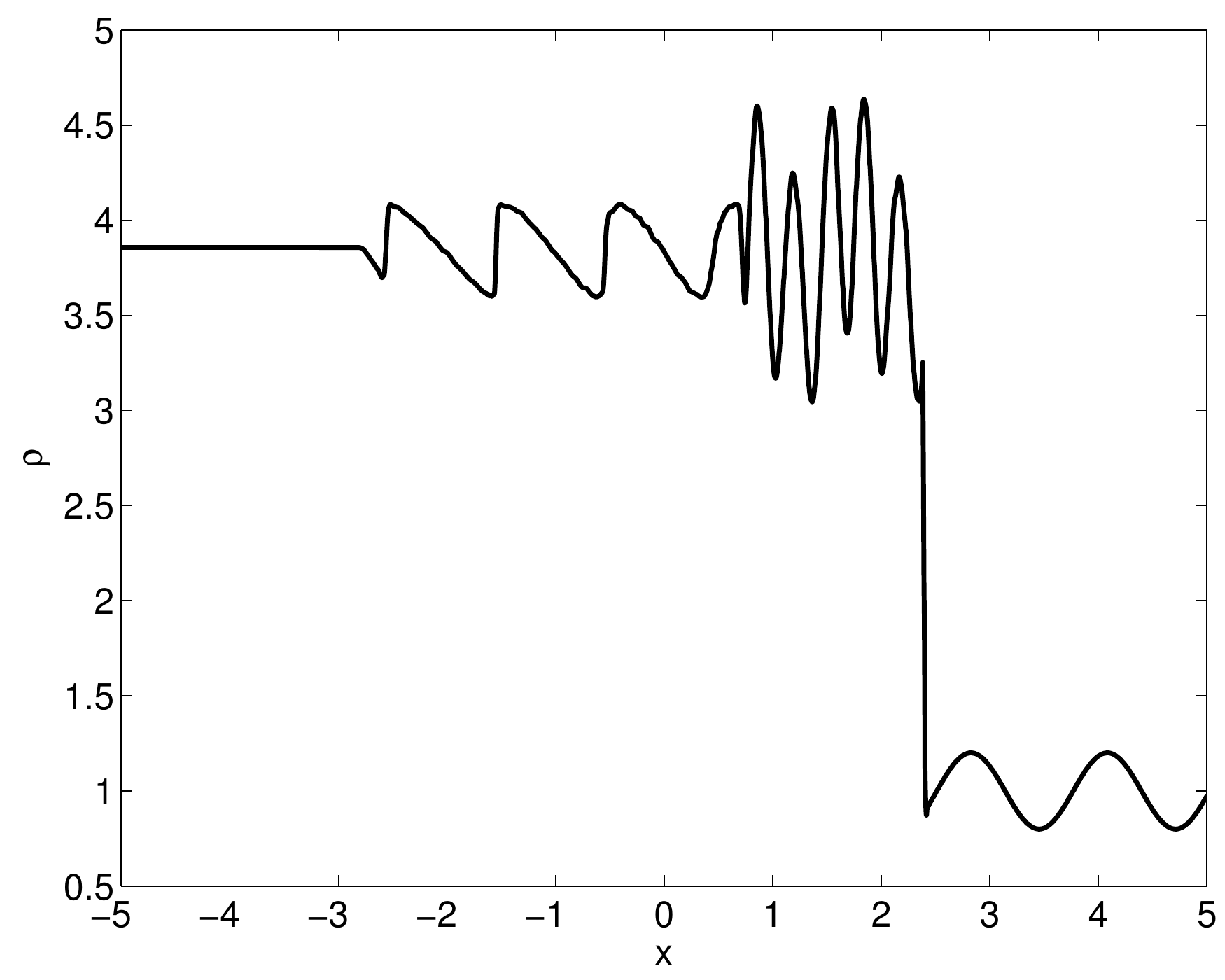}}
\caption{Approximation at $T=1.8$, KXRCF or Harten indicator ($\alpha = 1.5$), using density and entropy, Shu-Osher problem, 512 elements.}\label{fig:SineKHsol}
\end{figure}

For this example, the unmodified moment limited results are given in Figure \ref{fig:Sinemoment}.  It is clearly visible that almost every element is limited, and the peaks in the oscillating region at the left side of the initial discontinuity ($0.5 \leq x \leq 2$) have been lowered. Furthermore, the local smooth extrema at the right side of the initial discontinuity are limited.  The multiwavelet troubled-cell indicator was applied using the values of $C$ to  $0.5, 0.1$, and $0.05$ (Figures \ref{fig:SineC} and \ref{fig:SineCsol}, final time is $T=1.8$). Using $C=0.5$, the left shocks are not captured, and the only detected discontinuity is the strong shock in the initial condition. This is easily visible in Figure \ref{fig:SineC}. However, $C=0.1$ is perhaps more useful: in Figure \ref{fig:SineC}, we recognize the newly formed shocks. The solution looks much better in this region. The value $C=0.05$ may be too small: the continuous oscillating region is detected as well.  The results using the KXRCF or Harten's indicator can be inspected in Figures \ref{fig:SineKH} and \ref{fig:SineKHsol}. The KXRCF  indicator is very poor: in the linear case, the initial discontinuity is detected only. Therefore, the solution is very oscillatory. For the quadratic case, density selects part of the two left shocks, which leads to better results.  Harten's indicator works well in both the linear and the quadratic case, but the indicated troubled cells are much more scattered. The multiwavelet indicator using a convenient value of $C$ performs well for this case. 

\newpage

\newpage

\newpage

\subsubsection{Discussion}

 In Table \ref{tab:compare}, the average and maximum percentages of troubled elements in time are compared for each test problem, as in \cite{Qiu05S_2}. Generally, a troubled-cell indicator is said to be more accurate if smaller percentages of troubled cells are found. However, the approximation allows for more oscillations when fewer elements or incorrect regions are detected. This can, for example, be seen in the linear cases of Sod's problem (KXRCF), Lax's problem ($C=0.9$), and the Shu-Osher problem ($C=0.5$), where, although the smallest percentage is found, the resulting approximation oscillates.  Marked in bold are the smallest averages percentages that give rise to a nonoscillatory solution. It seems that the multiwavelet indicator leads to the best results, thereby detecting the smallest possible percentages.

\begin{table}[ht!] \scriptsize 
\begin{center}\resizebox{\columnwidth}{!}{
 \begin{tabular}{|c|c|c|c|c|c|c|c|c|c|}  \hline 
\multicolumn{10}{|c|}{Sod, 64 elements, $k=1$ and $k=2$} \\
\hline 
 \multicolumn{2}{|c|}{$C=0.9$} &  \multicolumn{2}{c|}{$C=0.5$} &  \multicolumn{2}{c|}{$C=0.1$} &  \multicolumn{2}{c|}{KXRCF} &  \multicolumn{2}{c|}{Harten, $\alpha = 1.5$}\\
\hline
Ave & Max & Ave & Max & Ave & Max& Ave & Max & Ave & Max \\
\hline
\textbf{1.8342} &  4.6875 &  5.0272 & 10.9375 &  17.1188 &  26.5625 & 1.4096  & 4.6875 & 8.4659 & 15.6250 \\
\textbf{1.8587} &  4.6875 &  3.4539 & 14.0625 &  13.4046 &  23.4375 & 2.9337 & 6.2500  & 4.4768 & 12.5000  \\
  \hline
\hline
 \multicolumn{10}{|c|}{Lax, 128 elements, $k=1$ and $k=2$} \\
\hline
 \multicolumn{2}{|c|}{$C=0.9$} &  \multicolumn{2}{c|}{$C=0.5$} &  \multicolumn{2}{c|}{$C=0.1$} &  \multicolumn{2}{c|}{KXRCF} &  \multicolumn{2}{c|}{Harten, $\alpha = 1.5$}\\
\hline
Ave & Max & Ave & Max & Ave & Max& Ave & Max & Ave & Max \\
\hline
 0.9176 & 2.3438  & 1.8977  &  6.2500  &   \textbf{5.5439} &  7.8125 & 2.3946 & 3.9062 & 9.2553 & 17.1875 \\
 \textbf{0.9962} & 3.1250  &  1.8974 &   6.2500 &  5.3720 & 9.3750 & 3.3744 & 5.4688 & 2.4032 &  9.3750\\
  \hline
\hline
 \multicolumn{10}{|c|}{Blast, 512 elements, $k=1$ and $k=2$} \\
\hline
 \multicolumn{2}{|c|}{$C=0.25$} &  \multicolumn{2}{c|}{$C=0.1$} & \multicolumn{2}{|c|}{$C=0.05$}& \multicolumn{2}{c|}{KXRCF}&  \multicolumn{2}{c|}{} \\
\hline
Ave & Max & Ave & Max & Ave & Max& Ave & Max & Ave & Max \\
\hline
\textbf{1.0758} & 3.1250 &  1.9260 & 4.1016 & 2.6178 & 6.2500 & 6.9520 & 9.3750 &  &  \\
\hline
\multicolumn{2}{|c|}{$C=0.1$} &  \multicolumn{2}{c|}{$C=0.05$} &  \multicolumn{2}{c|}{$C=0.01$} &  \multicolumn{2}{c|}{KXRCF}& \multicolumn{2}{c|}{}  \\
\hline
Ave & Max & Ave & Max & Ave & Max& Ave & Max & Ave & Max \\
\hline
\textbf{1.6792} & 5.2734 &  2.3354 & 7.0312 & 3.9640 &  11.9141 &  12.3097 &  20.1172 &  &  \\
  \hline
\hline
 \multicolumn{10}{|c|}{Shu-Osher, 512 elements, $k=1$ and $k=2$} \\
\hline
 \multicolumn{2}{|c|}{$C=0.5$} &  \multicolumn{2}{c|}{$C=0.1$} &  \multicolumn{2}{c|}{$C=0.05$} &  \multicolumn{2}{c|}{KXRCF} &  \multicolumn{2}{c|}{Harten, $\alpha = 1.5$} \\
\hline
Ave & Max & Ave & Max & Ave & Max & Ave & Max & Ave & Max\\
\hline
0.3849 &  0.7812 & \textbf{0.8596}  & 3.3203  & 2.7921 & 14.0625  & 0.6237  &1.1719 & 4.2886  & 9.1797  \\
\textbf{0.3787} & 1.5625 & 0.8093  &  3.9062  & 1.2584 &  8.7891 & 1.2059  & 2.1484 & 2.4105  & 6.2500 \\
  \hline  \end{tabular}}\caption{Average and maximum percentages of cells that are indicated as troubled by our multiwavelet detector, for different $C$, the KXRCF indicator or Harten's indicator. For the Blast wave example, Harten's indicator did not work correctly, and different values of $C$ are used for $k=1$ and $k=2$. Marked in bold are the smallest average percentages that belong to a nonoscillatory solution.}\label{tab:compare}
\end{center}
\end{table}

A useful property of troubled-cell indicators is the decrease of percentages if the resolution is increased \cite{Qiu05S_2}. In all examples, $C=0.1$ is a good choice for detecting troubled cells. We therefore keep it fixed for each example, and double the number of elements in our discretization. The percentages of troubled cells are approximately halved, which can be seen in Table \ref{tab:resolution}. This nice behavior is due to the multiwavelet approach. The KXRCF and Harten's troubled-cell indicator have the same property, although the rate of decrease is smaller, \cite{Qiu05S_2}. 

\begin{table}[ht!] \scriptsize
\begin{center}
 \begin{tabular}{|c|c|c|c|c|c|c|c|c|}  \hline 
&\multicolumn{2}{|c|}{Sod} & \multicolumn{2}{c|}{Lax} &\multicolumn{2}{|c|}{Blast} & \multicolumn{2}{c|}{Shu-Osher}\\
&\multicolumn{2}{|c|}{128 elements} & \multicolumn{2}{c|}{256 elements} &\multicolumn{2}{|c|}{1024 elements} & \multicolumn{2}{c|}{1024 elements} \\
\hline
& Ave & Max & Ave & Max &Ave & Max & Ave & Max\\
\hline
$k=1$ & 7.9687 & 13.2812 & 2.9076 & 3.9062 &  0.9002 & 2.5391& 0.3857& 1.1719\\
$k=2$ &  6.7245 & 14.0625 &2.8696 &  5.0781 &  0.7948 & 2.6367&  0.3877&  1.4648\\
\hline
 \end{tabular}
 \caption{Average and maximum percentages of cells that are indicated as troubled by our multiwavelet detector ($C=0.1$), using twice as many elements as in Table \ref{tab:compare}. 
 \label{tab:resolution}}
\end{center}
\end{table}

\subsection{Two-dimensions:  Double Mach reflection}
The performance of the multiwavelet troubled-cell indicator is now considered in two-dimensions for the double Mach reflection problem \cite{Woo84C}.  Here, it is compared with Krivodonova et al.'s indicator \cite{Kri04XRCF}.  

Using the vector ${\bf u} = (\rho, \rho u, \rho v, E)^\top$, the two-dimensional Euler equations are given by
\begin{subequations}
 \begin{alignat}{5}
   {\bf u}_t + {\bf f}({\bf u})_x + {\bf g}({\bf u})_y & = {\bf 0}, & {\bf x} \in \Omega, t \geq 0, \\
   {\bf u}({\bf x}, 0) &= {\bf u}_0({\bf x}), & {\bf x} \in \Omega.
 \end{alignat}
where
\begin{align}
 {\bf f}({\bf u}) &= \left(\rho u, \rho u^2 + p, \rho uv, (E+p)u \right)^\top; \\
 {\bf g}({\bf u}) &= \left(\rho v, \rho uv, \rho v^2 +p, (E+p)v \right)^\top,
\end{align}
\end{subequations}
and the equation of state is given by $p = (\gamma - 1)(E-\rho(u^2+v^2)/2)$. 

The computational domain of this problem is $[0,4] \times [0,1]$. At $t=0$, this domain is divided into two regions that are separated by 
$y(x) = \sqrt{3}(x-1/6)$ \cite{Tan05W}. The following initial conditions are used:
\begin{subequations}
\begin{align}
  {\bf u}_L &= (8, 8.25 \cos (30^\circ),-8.25 \sin (30^\circ), 563.5)^\top; \\
  {\bf u}_R &= (1.4,0,0, 2.5)^\top.
\end{align}
\end{subequations}

At the left boundary, ${\bf u}_L$ is used as a boundary condition, and at the right, ${\bf u}_R$.  The top boundary is divided into two regions: for $x < 1/6 + (1+20t)/\sqrt{3}$, ${\bf u}_L$ is used, whereas ${\bf u}_R$ is used to the right.   At the bottom boundary, ${\bf u}_L$ is used for $x < 1/6$, and a reflecting wall is used for $x \geq 1/6$. 

\begin{figure}[ht!]
\centering
 \subfigure[$\alpha$ mode, $C=0.05$]{\includegraphics[scale = 0.34]{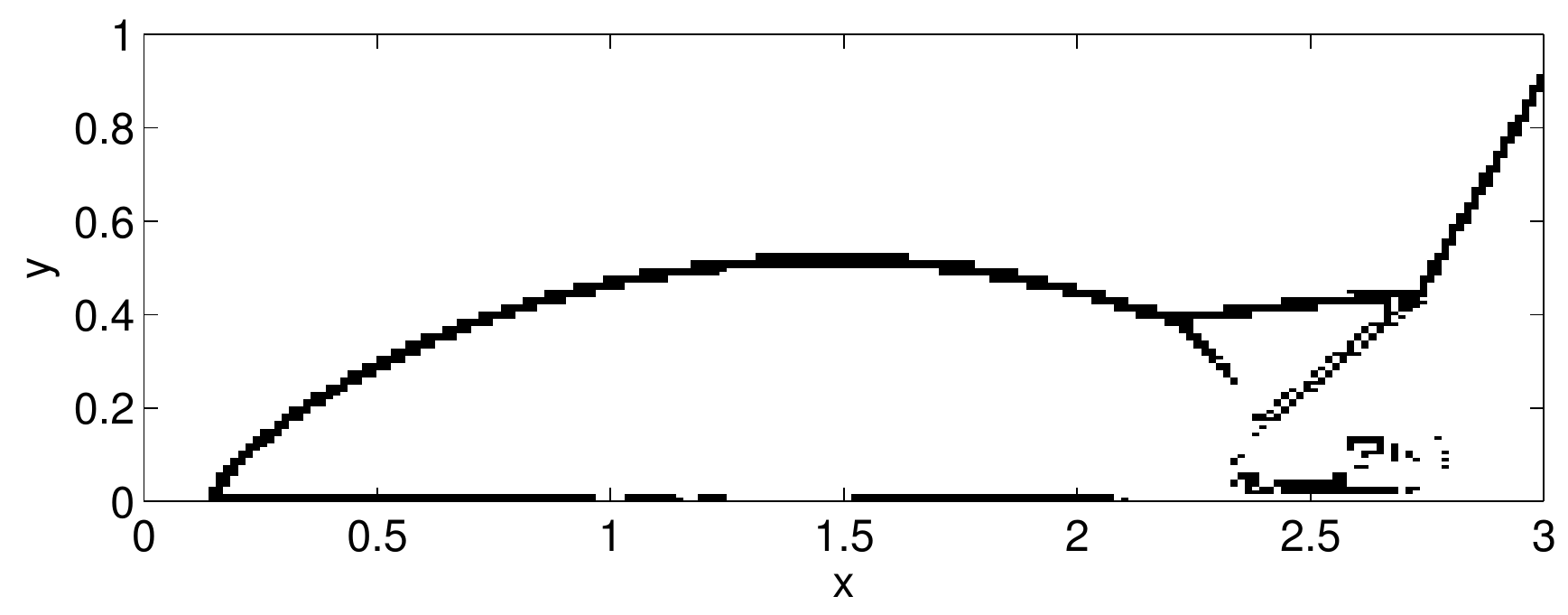}}
 \subfigure[$\beta$ mode, $C=0.05$]{\includegraphics[scale = 0.34]{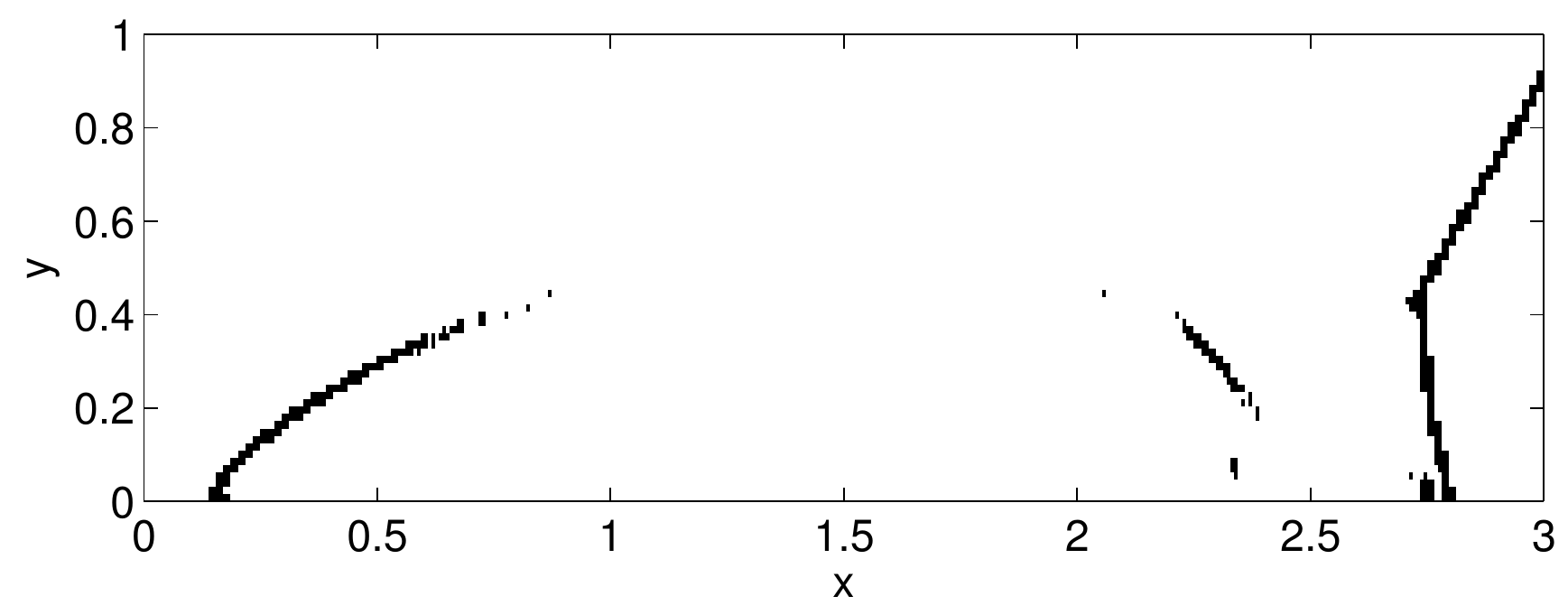}} \\
 \subfigure[$\gamma$ mode, $C=0.05$]{\includegraphics[scale = 0.34]{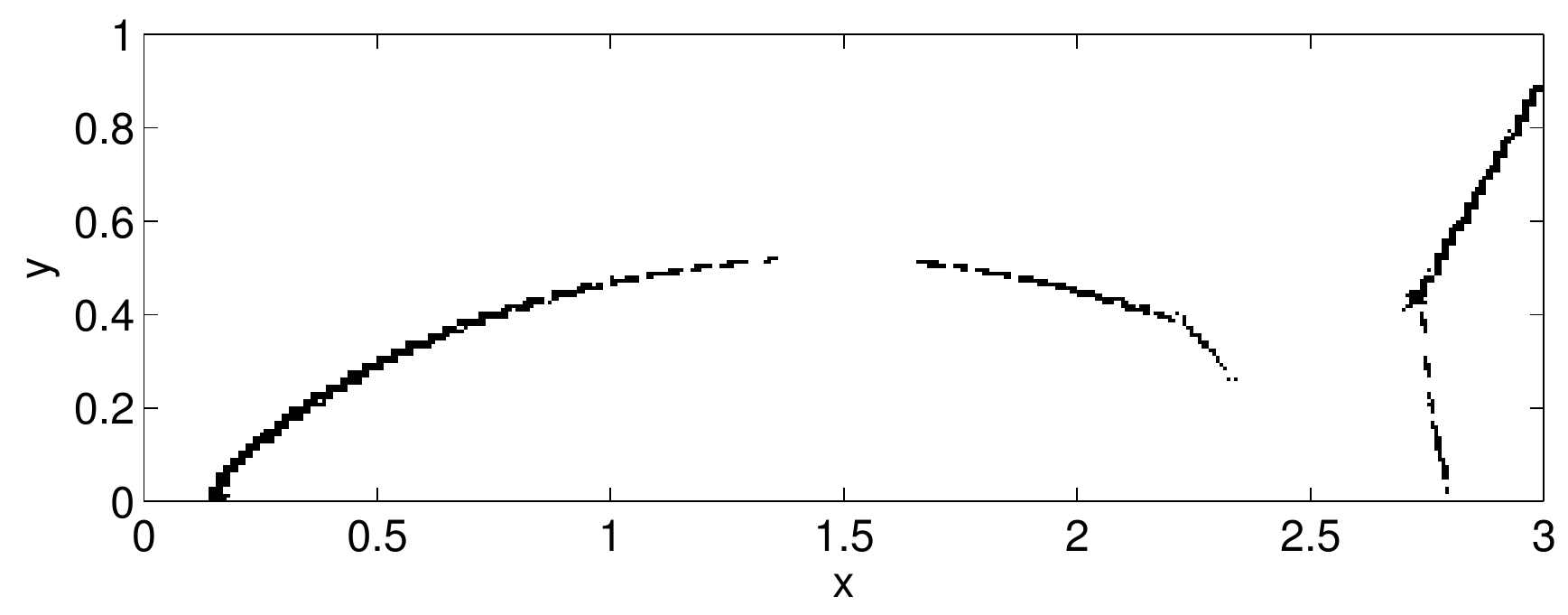}} 
 \subfigure[Combination, $C=0.05$]{\includegraphics[scale = 0.34]{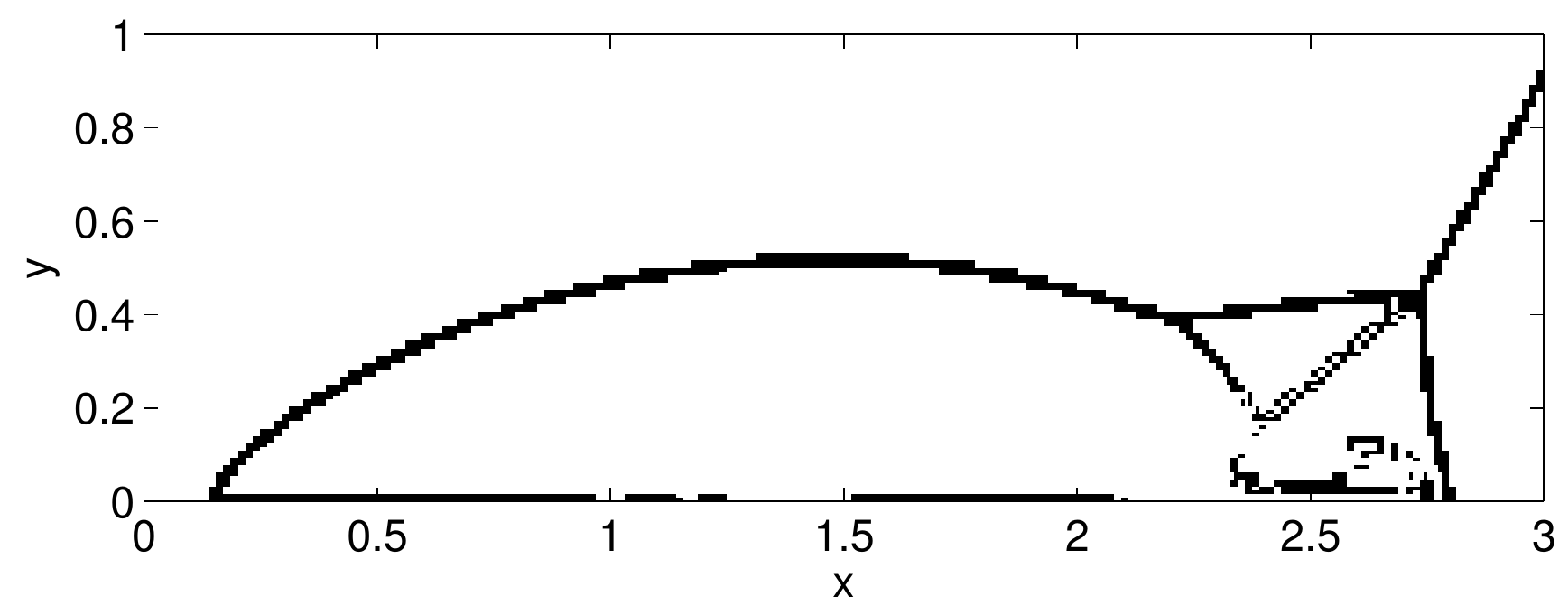}} \\
 \subfigure[KXRCF]{\includegraphics[scale = 0.34]{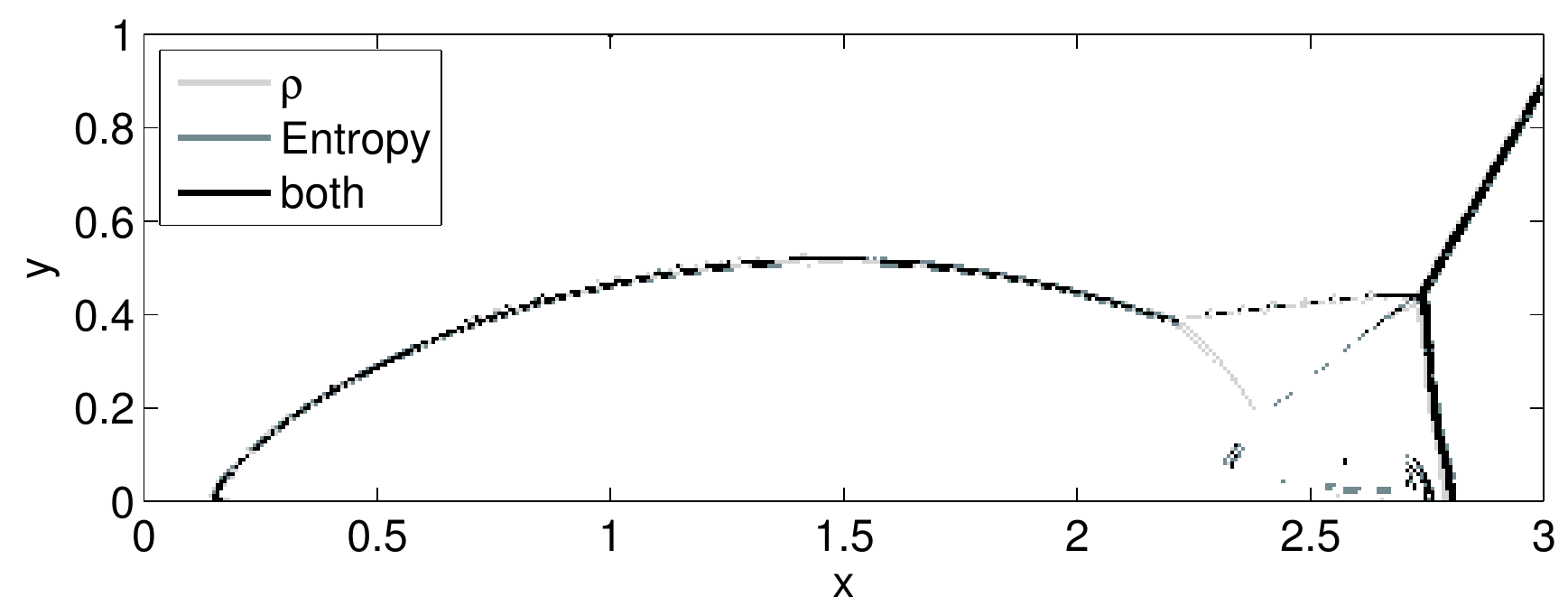}}
\caption{Detected troubled cells at $T=0.2$, double Mach reflection problem, $k=1$, $\Delta x = \Delta y = 1/128$.}\label{fig:DoubleMachk1C}
\end{figure}

\begin{figure}[ht!]
\centering
 \subfigure[$\alpha$ mode, $C=0.05$]{\includegraphics[scale = 0.34]{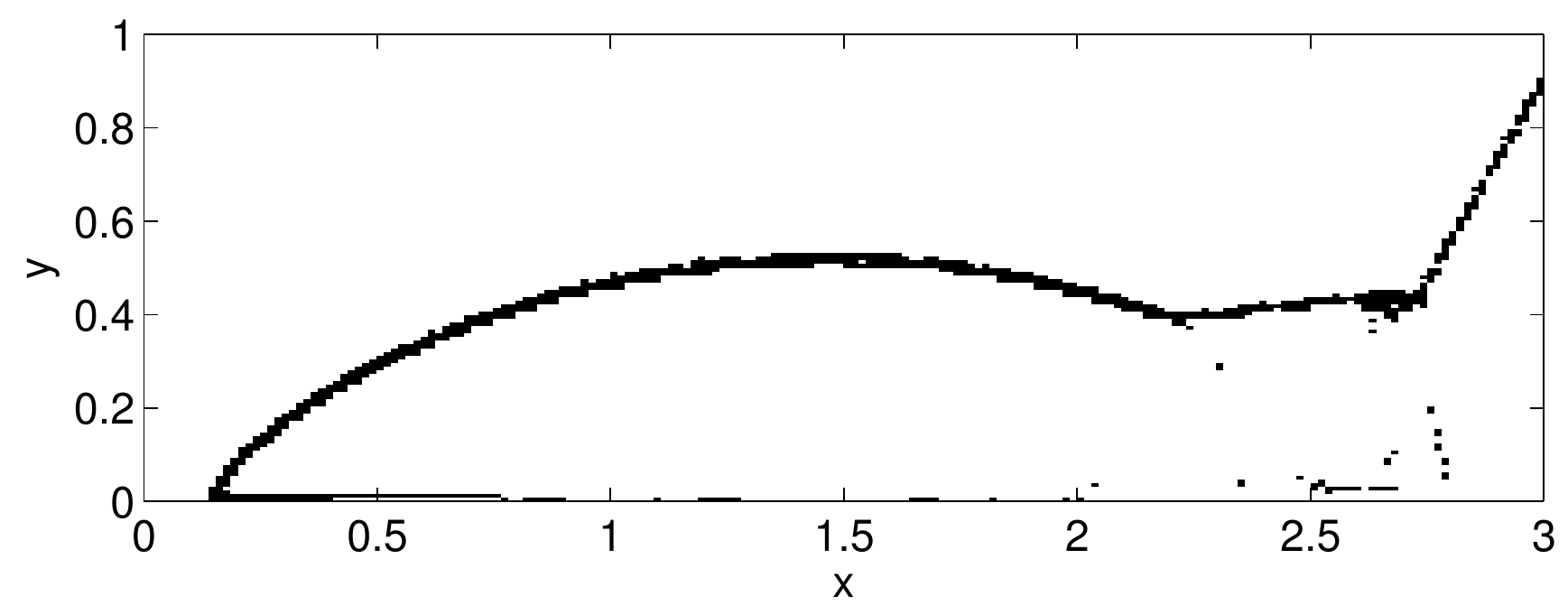}}
 \subfigure[$\beta$ mode, $C=0.05$]{\includegraphics[scale = 0.34]{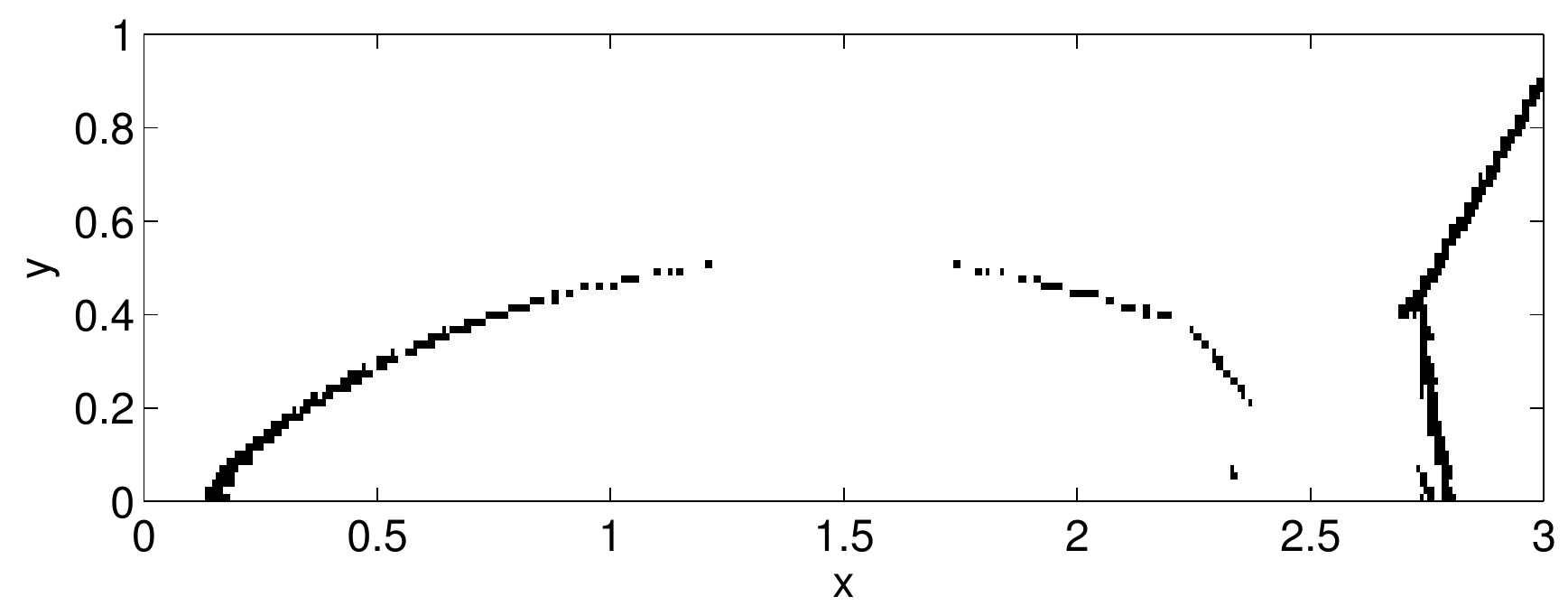}} \\
 \subfigure[$\gamma$ mode, $C=0.05$]{\includegraphics[scale = 0.34]{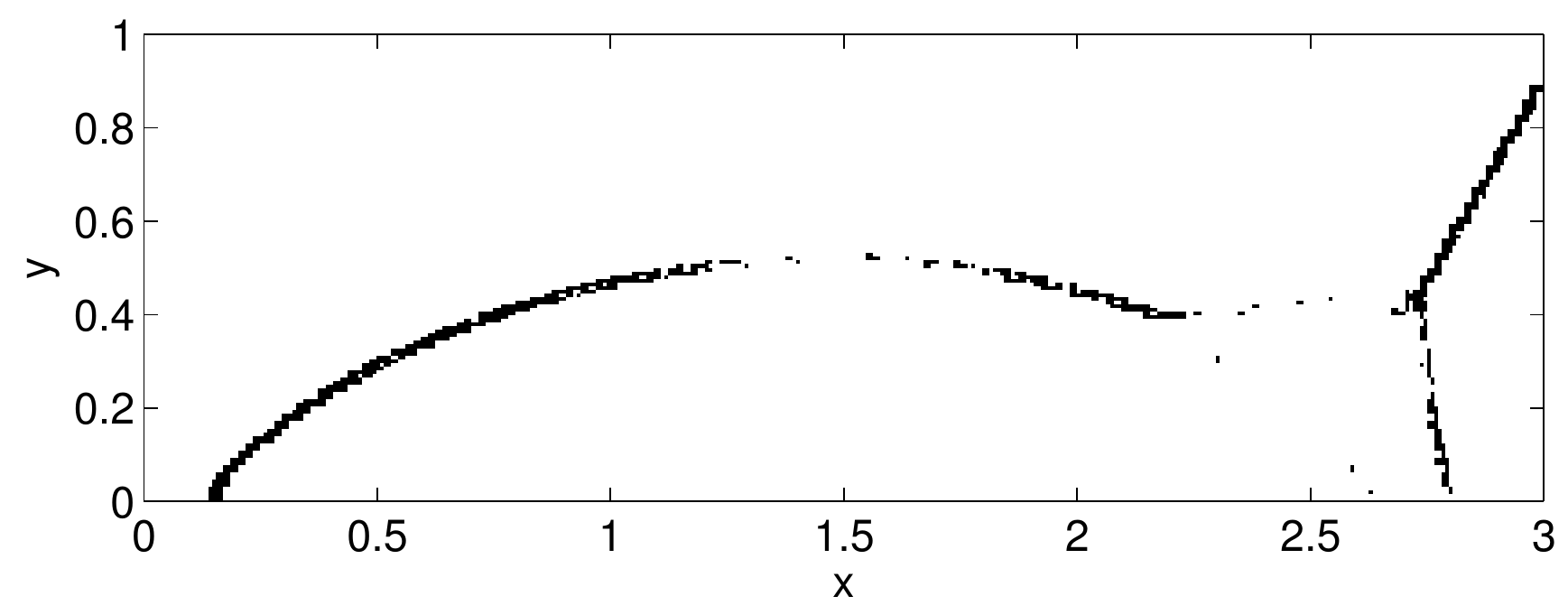}}
 \subfigure[Combination, $C=0.05$]{\includegraphics[scale = 0.34]{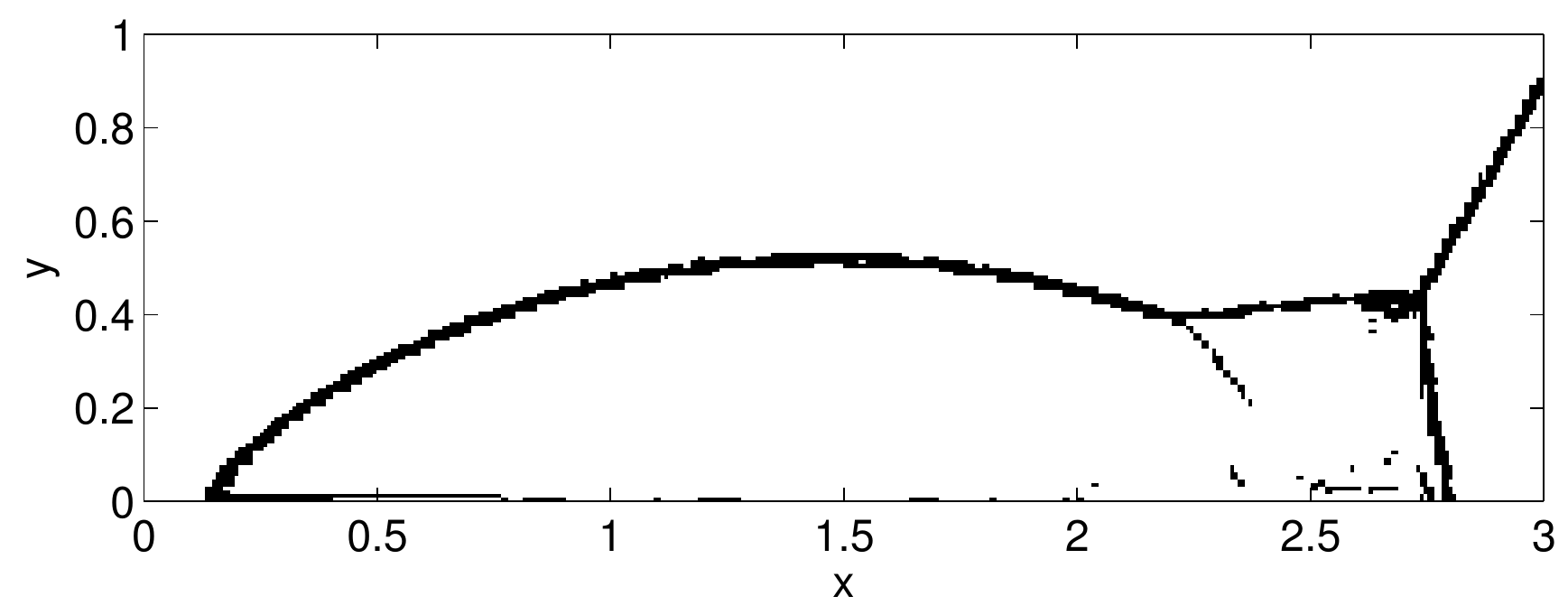}} \\
 \subfigure[KXRCF]{\includegraphics[scale = 0.34]{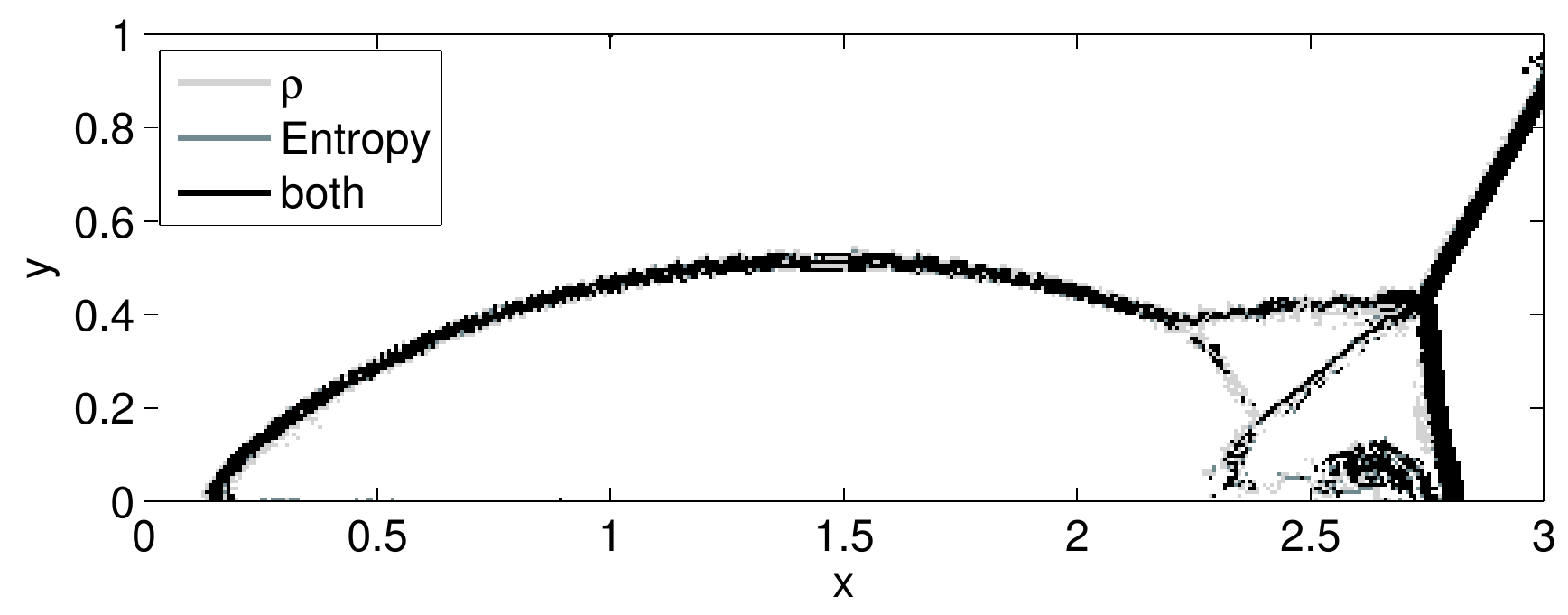}}
\caption{Detected troubled cells at $T=0.2$, double Mach reflection problem, $k=2$, $\Delta x = \Delta y = 1/128$.}\label{fig:DoubleMachk2C}
\end{figure}

\begin{figure}[ht!]
\centering
 \subfigure[$C=0.05$]{\includegraphics[scale = 0.4]{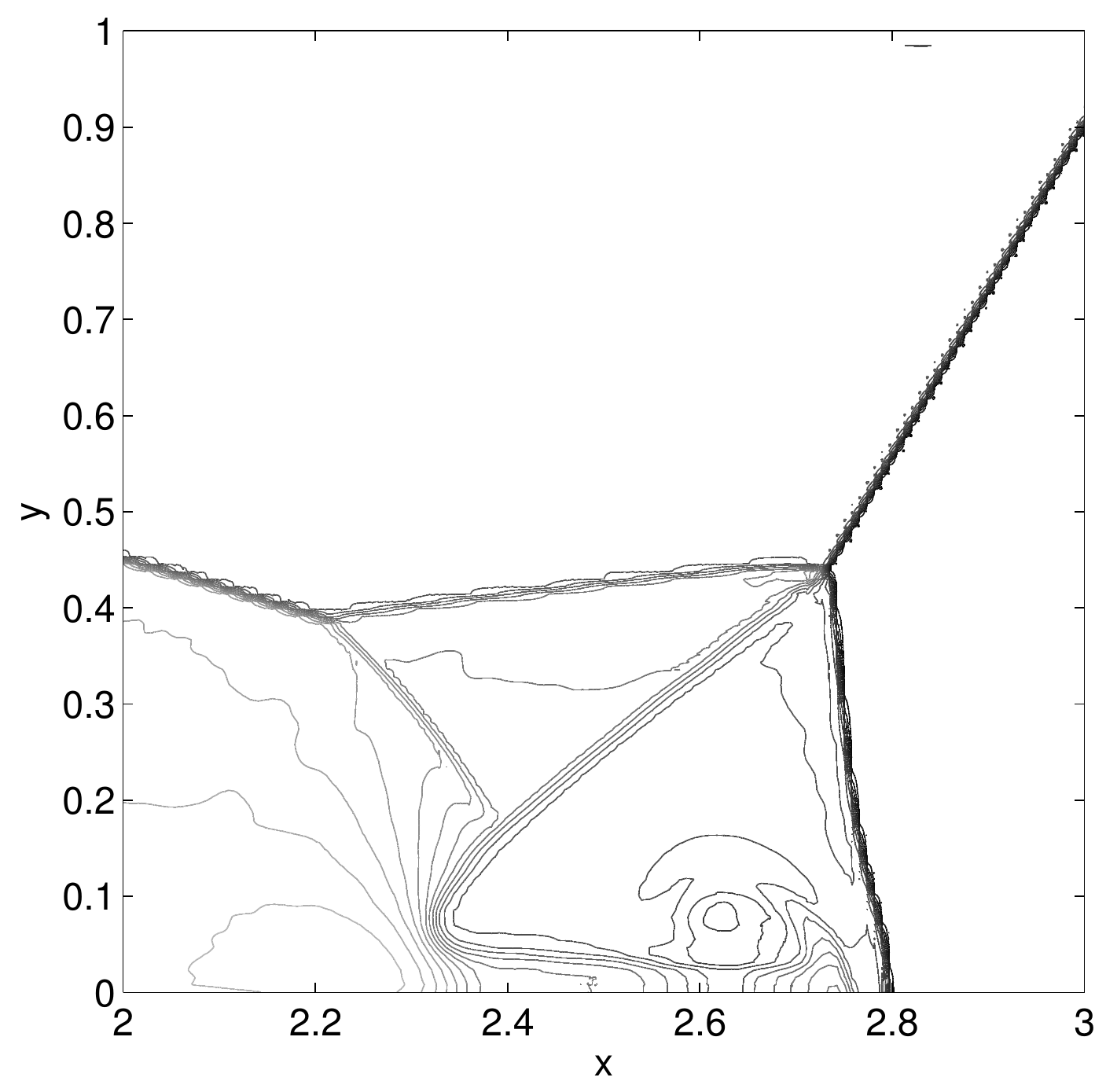}}
 \subfigure[KXRCF]{\includegraphics[scale = 0.4]{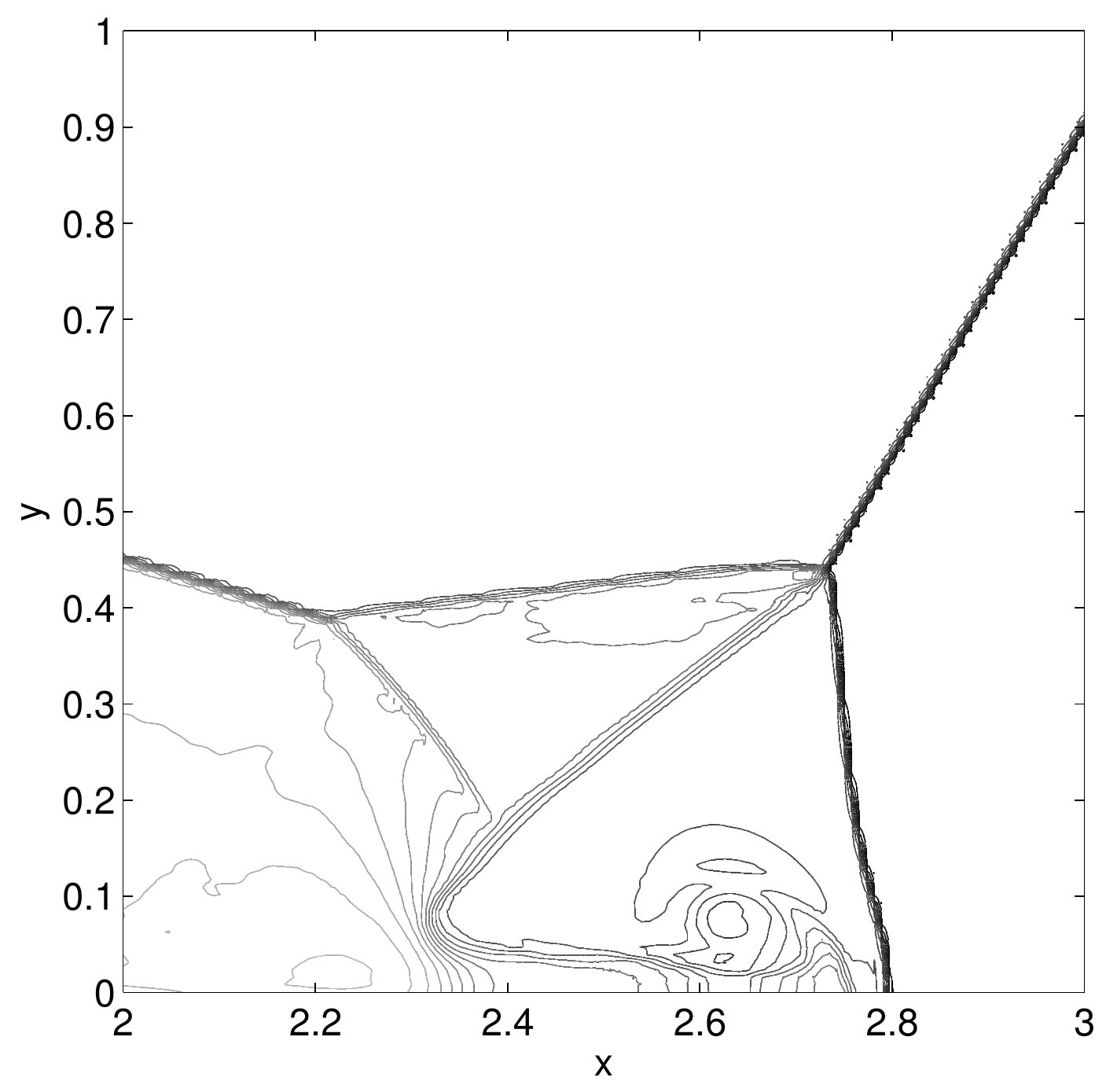}} \\

 \subfigure[$C=0.05$]{\includegraphics[scale = 0.4]{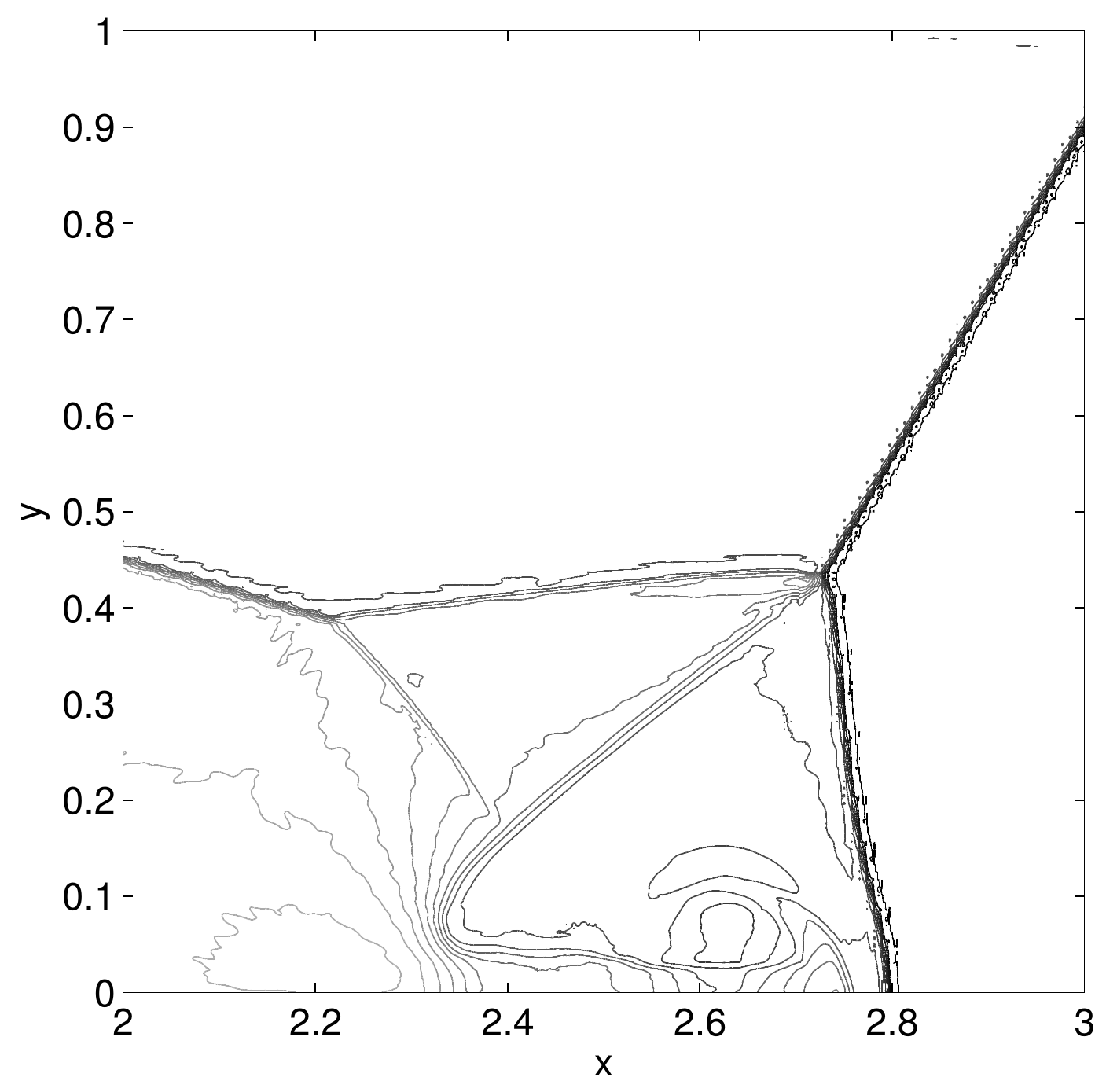}}
 \subfigure[KXRCF]{\includegraphics[scale = 0.4]{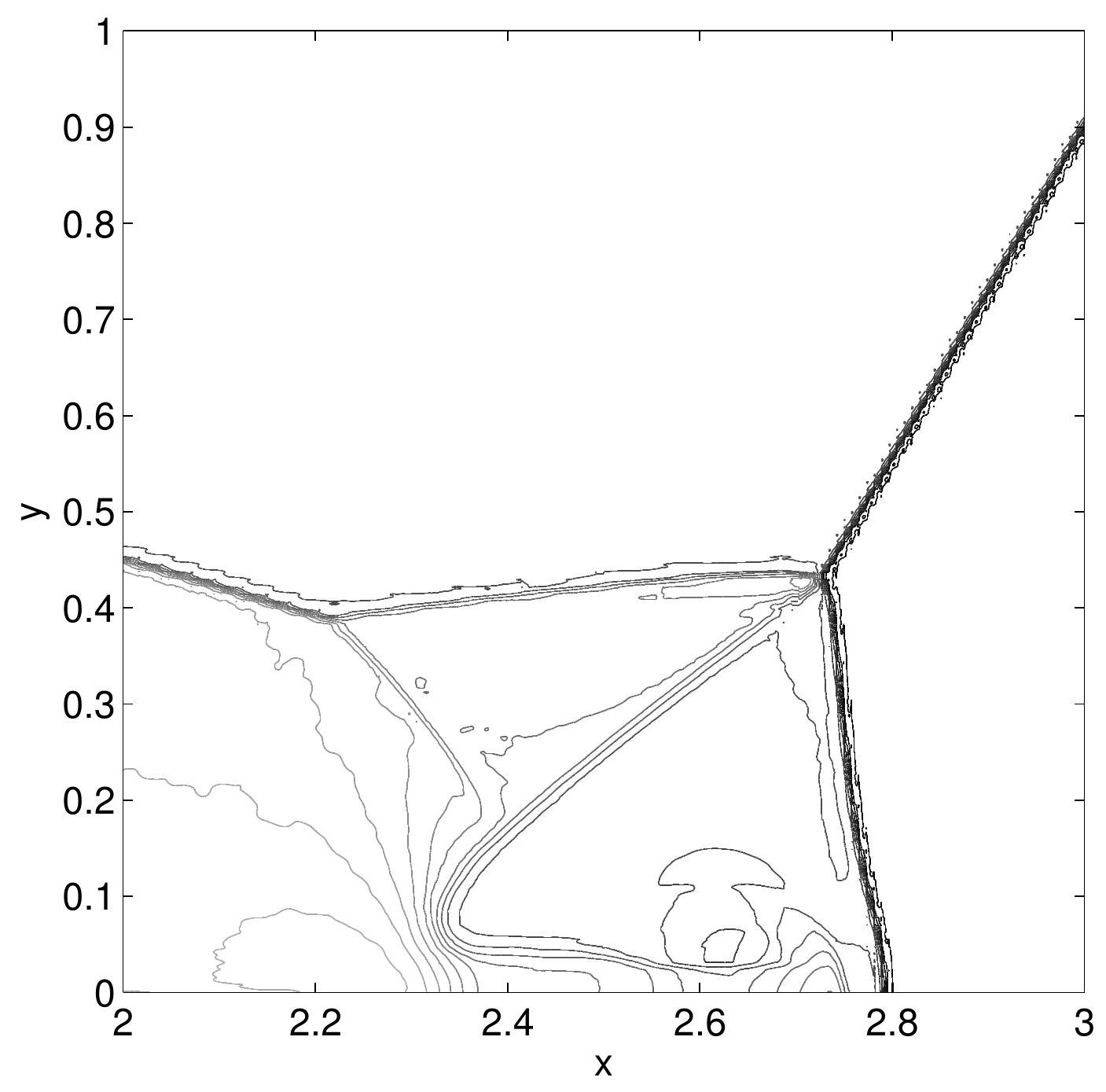}} \\
\caption{Contour lines of approximation, double Mach reflection problem at $T=0.2$, $\Delta x = \Delta y = 1/128$. First row: $k=1$, second row: $k=2$.}\label{fig:DoubleMachk2KXRCF}
\end{figure}

The results at $T=0.2$ using $\Delta x = \Delta y = 1/128$ are given in Figures \ref{fig:DoubleMachk1C} to \ref{fig:DoubleMachk2KXRCF} along with the identified troubled cells. The percentages of detected troubled cells are given in Table \ref{tab:percentages}. As in \cite{Mal98}, we can see that the $\alpha, \beta$ and $\gamma$ modes detect different troubled cells based on direction. Although Qiu et al. use both density and entropy or density and energy to compute troubled cells \cite{Qiu05S_2}, for our multiwavelet indicator using only density is enough to detect troubled cells. The use of entropy does not produce significant changes for the detected troubled-cell regions. The KXRCF indicator using density and entropy as indicator variables works very good, detecting exactly the discontinuous regions in the solution. Note that for $k=2$, more elements are detected in the turbulent region than if the multiwavelet indicator ($C=0.05$) is used. The approximate solutions of the multiwavelet and the KXRCF approach look quite similar. In the turbulent region, more details of the DG simulation can be seen because we allow the solution to oscillate in continuous regions. Similar to the Blast wave problem, the combination of Harten's troubled-cell indicator and the moment limiter was found to be unstable and therefore we do not include this comparison. 

Using a troubled-cell indicator, the moment limiter is applied only in a small portion of the elements, whereas the unmodified moment limiter limits every element. Therefore, the total computation time decreases by using a troubled-cell indicator.  In Table \ref{tab:computationtime}, the total computation times using one of the different indicators can be compared. It is clear that the multiwavelet indicator is faster than the KXRCF approach. 

For high resolution computations, the same behavior is found: the indicator perfectly finds the troubled regions, and the moment limiter is applied only in these elements. Because the element size decreases, the limited region itself becomes smaller.

\begin{table}[ht!]
\begin{center}
 \begin{tabular}{|c|c|c|}  \hline 
$k$ & $C=0.05$ & KXRCF \\
\hline
1 & 50 & 85 \\
2 & 214 & 335 \\
  \hline \end{tabular}\caption{Total computation time in minutes for double Mach, $\Delta x = \Delta y = 1/128$.}\label{tab:computationtime}
\end{center}
\end{table}

\begin{table}[ht!] 
\begin{center}
 \begin{tabular}{|c|c|c|c|c|}  \hline 
 & \multicolumn{2}{c|}{$C=0.05$} & \multicolumn{2}{c|}{KXRCF}\\
\hline
$k$ & Ave & Max & Ave & Max \\
\hline
1 & 2.2916 & 4.0115 & 1.5190 & 2.2629\\
2 & 2.0978 & 3.0106 &  3.3784 & 5.3650 \\
  \hline \end{tabular}\caption{Average and maximum percentages of cells that are indicated as troubled for the double Mach reflection problem, $\Delta x = \Delta y = 1/128$. }\label{tab:percentages}
\end{center}
\end{table}

\clearpage
\section{Conclusions}
\label{sec:conclude}
In this paper we have introduced a {\it global} multiwavelet troubled-cell indicator.  This technique relates the DG approximation to the multiwavelet expansion and uses information from the multiwavelet expansion in order to identify troubled cells.  In the numerical results, we demonstrated that this technique performs well, even in the vicinity of a strong shock with weaker local shocks and has a robust performance compared with other methods.  Furthermore, our results showed that, because of the choice in how the multiwavelet expansion is implemented, it performs faster than the currently used troubled-cell indicators. Future work will be to see if we can improve upon the performance in detecting local structures, to decide in advance which value of the parameter we should use, and to extend this to unstructured meshes.

\noindent {\bf Acknowledgements:}  The authors gratefully wish to acknowledge the useful comments provided by Lilia Krivodonova, Jianxian Qiu, Chi-Wang Shu and Arnold Heemink that helped to shape this work.

\bibliographystyle{plain}
\bibliography{References}

\end{document}